\documentclass[preprint,12pt,sort,compress,numbers]{elsarticle}

\usepackage{soul}
\usepackage{lineno}
\input{packages.tex}
\def\Dpartial#1#2{ \frac{\partial #1}{\partial #2} }
\def\Dparttwo#1#2{ \frac{\partial^2 #1 }{ \partial #2^2} }
\def\DpartN#1#2#3{ \frac{\partial^{#3} #1 }{\partial #2^{#3}} }

\newif\ifrefcolor
\refcolorfalse

\def\refA#1{\ifrefcolor\textcolor{black!30!green}{#1}\else#1\fi}
\def\refB#1{\ifrefcolor{\color{purple}#1}\else#1\fi}

\newcommand{\LHS}{\mathrm{LHS}}

\newcommand{\tW}{\tilde{W}}
\newcommand{\tOmega}{\tilde{\Omega}}

\newcommand{\cO}{\mathcal{O}}

\newcommand{\cJ}{\mathcal{J}}

\newcommand{\rB}{\mathbb{B}}
\newcommand{\rR}{\mathbb{R}}
\newcommand{\rI}{\mathbb{I}}

\newcommand{\bA}{\mathbf{A}}

\newcommand{\bn}{\mathbf{n}}
\newcommand{\bc}{\mathbf{c}}

\newcommand{\bx}{\mathbf{x}}

\def\jump#1{\left\llbracket #1 \right\rrbracket}
\def\avg#1{\left\{\!\!\left\{ #1 \right\}\!\!\right\}}

\newcommand{\logLogSlopeTriangle}[6]
{

    \pgfplotsextra
    {
        \pgfkeysgetvalue{/pgfplots/xmin}{\xmin}
        \pgfkeysgetvalue{/pgfplots/xmax}{\xmax}
        \pgfkeysgetvalue{/pgfplots/ymin}{\ymin}
        \pgfkeysgetvalue{/pgfplots/ymax}{\ymax}

        \pgfmathsetmacro{\xArel}{#1}
        \pgfmathsetmacro{\yArel}{#3}
        \pgfmathsetmacro{\xBrel}{#1-#2}
        \pgfmathsetmacro{\yBrel}{\yArel}
        \pgfmathsetmacro{\xCrel}{\xArel}

        \pgfmathsetmacro{\lnxB}{\xmin*(1-(#1-#2))+\xmax*(#1-#2)} 
        \pgfmathsetmacro{\lnxA}{\xmin*(1-#1)+\xmax*#1} 
        \pgfmathsetmacro{\lnyA}{\ymin*(1-#3)+\ymax*#3} 
        \pgfmathsetmacro{\lnyC}{\lnyA+#4*(\lnxA-\lnxB)}
        \pgfmathsetmacro{\yCrel}{(\lnyC-\ymin)/(\ymax-\ymin)} 

        \coordinate (A) at (rel axis cs:\xArel,\yArel);
        \coordinate (B) at (rel axis cs:\xBrel,\yBrel);
        \coordinate (C) at (rel axis cs:\xCrel,\yCrel);

        \draw[#5]   (A)-- node[pos=0.5,anchor=#6] {1}
                    (B)-- 
                    (C)-- node[pos=0.5,anchor=west] {#4}
                    cycle;
    }
}

\pgfplotsset{select coords between index/.style 2 args={
    x filter/.code={
        \ifnum\coordindex<#1\fi
        \ifnum\coordindex>#2\fi
    }
}}

\begin{document}

\begin{frontmatter}



\title{
Hard-constrained Physics-informed Neural Networks for Interface Problems
}


\author[casc]{Seung Whan Chung\corref{cor1}}
\cortext[cor1]{corresponding authors}
\ead{chung28@llnl.gov}
\author[ced]{Stephen T. Castonguay}
\author[aeed]{Sumanta Roy}
\author[sandia]{Michael S. Penwarden}
\author[pnnl]{Yucheng Fu}
\author[aeed]{Pratanu Roy\corref{cor1}}
\ead{roy23@llnl.gov}
\affiliation[casc]{organization={Center for Applied Scientific Computing, Lawrence Livermore National Laboratory},
  city={Livermore},
  state={CA},
  postcode={94550}, 
  country={USA}}

\affiliation[aeed]{organization={Atmospheric, Earth and Energy Division, Lawrence Livermore National Laboratory},
city={Livermore},
state={CA},
postcode={94550}, 
country={USA}}
\affiliation[sandia]{organization={Sandia National Laboratories},
city={Albuquerque},
state={NM},
postcode={87123},
country={USA}}
\affiliation[ced]{organization={Computational Engineering Division, Lawrence Livermore National Laboratory},
city={Livermore},
state={CA},
postcode={94550}, 
country={USA}}
\affiliation[pnnl]{organization={Physical and Computational Sciences Directorate, Pacific Northwest National Laboratory},
city={Richland},
state={WA},
postcode={99354},
country={USA}}

\begin{abstract}
Physics-informed neural networks (PINNs) have emerged as a flexible framework for solving partial differential equations, but their performance on interface problems remains challenging because continuity and flux conditions are typically imposed through soft penalty terms. The standard soft-constraint formulation leads to imperfect interface enforcement and degraded accuracy near interfaces. We introduce two ansatz-based hard-constrained PINN formulations  for interface problems that embed the interface physics into the solution representation and thereby decouple interface enforcement from PDE residual minimization. The first, termed the \textbf{windowing approach}, constructs the trial space from compactly supported windowed subnetworks so that interface continuity and flux balance are satisfied by design. The second, called the \textbf{buffer approach}, \refB{discretely hard-constrains subnetworks with additive} buffer functions that enforce boundary and interface constraints at \refB{sampled} points through a lightweight correction. We study these formulations on one- and two-dimensional elliptic interface benchmarks and compare them with soft-constrained baselines.
\refB{A statistical study on the one-dimensional cases shows that} the windowing approach attains very high accuracy (as low as $\mathcal{O}(10^{-9})$) \refB{in occasional runs} on simple structured cases\refB{, but its rigid ansatz makes training stiff, so its median errors are often worse than the soft-constrained baselines. In contrast,}
the buffer approach \refB{achieves the smallest median error (as low as $\mathcal{O}(10^{-5})$)} across a wider range of source terms and interface configurations. In two dimensions, the buffer formulation is shown to be more robust \refB{against general geometries with discrete hard-constraining, while} the windowing construction becomes more sensitive to overlap and corner effects and over-constrains the problem. 
This positions the buffer method as a straightforward and \refB{promising approach for interface problems, with the search for an optimal buffer function form left for future work.}
\end{abstract}




\end{frontmatter}

 \nolinenumbers

\section{Introduction}

Physics-informed neural networks (PINNs) embed partial differential equations (PDEs) into the training objective of neural networks. PINNs have emerged as a flexible, mesh-free alternative for forward and inverse problems across science and engineering \cite{karniadakis2021physics}. Since their introduction in \cite{lagaris1998artificial,raissi2019physics}, a large body of work has refined and advanced their approximation, optimization, and generalization properties \cite{wang2023expert}. Despite this progress, the accurate enforcement of boundary and interface conditions remains a central challenge---particularly for problems with discontinuous coefficients or nonsmooth solution structure. In such settings, standard PINN formulations often lose accuracy near interfaces, even when the PDE residual is fitted well in the interior.

Most existing PINN formulations for interface problems enforce continuity or jump conditions through penalty terms added to the loss \cite{jagtap2020extended,jagtap2020conservative,pmlr-v202-li23w}. This soft-constrained strategy is attractive because it is simple to implement and compatible with a wide range of architectures, including domain-decomposed and interface-aware PINN variants. However, it converts the training problem into a multi-objective optimization in which PDE residuals, boundary conditions, and interface conditions must be balanced against one another. In practice, this often leads to sensitivity to loss weights, imperfect satisfaction of interface physics, and degraded solution quality near discontinuities \cite{cuomo2022scientific, lu2021physics}, especially for high-contrast coefficients, multiple interfaces, or mixed boundary conditions. 
Several architectures have been proposed to improve soft-constrained interface learning, including domain-decomposed multi-domain PINN (M-PINN)~\cite{zhang2022multi}, interface PINNs (I-PINNs)~\cite{sarma2024interface} and their adaptive variants such as AdaI-PINNs~\cite{roy2024adaptive}, as well as discontinuity-capturing networks that augment the input space to represent piecewise solutions more effectively~\cite{hu2022discontinuity, tseng2023discontinuity, yao2023deep, fan2025novel, bi2025extended}.
While these approaches are flexible and easy to implement,
they are still subject to the inherent weakness of multi-objective optimization, leading to
inexact enforcement of interface conditions, sensitivity to penalty weights, imbalance among competing loss terms, and difficulty representing sharp jumps---especially in high-contrast or multi-interface problems~\cite{roy2025adaptive,sarkar2025adaptive}.

Hard-constrained PINN strategies overcome this issue by constructing solution ansatzes that exactly satisfy Dirichlet, Neumann, and Robin conditions, thus decoupling boundary enforcement from PDE residual enforcement. A common recipe is to multiply a learnable “free” network by a signed distance function from the boundary and add a particular solution that encodes the boundary data; the resulting ansatz renders the boundary loss unnecessary and improves stability~\cite{moseley2023finite}. Recent studies have formalized this approach, including exact Dirichlet enforcement via auxiliary networks for the distance function and boundary data and unified frameworks that extend hard enforcement to Neumann and Robin conditions~\cite{sukumar2022exact,wang2023exact, liu2022unified, roy2024exact}. Similar ideas have also appeared in variational PINN formulations, where structure-preserving treatments of essential boundary data reduce optimization stiffness by avoiding penalty terms \cite{liao2019deep, yu2025natural}.
\par
There are still important challenges for multiplicative hard-constraining approaches.
First, designing a multiplier that exactly satisfies Neumann conditions or boundary
conditions involving higher-order derivatives is not straightforward, especially for
inhomogeneous boundaries where multiple boundary conditions intersect
~\cite{sukumar2022exact,goschel2025enforcing}. One alternative is to
adopt a first-order formulation by introducing auxiliary variables that represent
derivatives of the solution~\cite{liu2023gradient,peng2025information}. In that
setting, boundary conditions on higher-order derivatives can be reformulated as
Dirichlet conditions on the auxiliary variables and then imposed exactly. Peng and
Tang~\cite{peng2025information} extended this idea to interface problems. However,
training still becomes a multi-objective optimization over the primary and auxiliary
variables and therefore remains sensitive to penalty weights and imbalance among the
competing loss terms. Interface problems introduce an additional difficulty because
field values or their gradients may be discontinuous across material or domain
boundaries. Tresckow \textit{et al.}~\cite{von2025multi} proposed a hard-constrained PINN for PDEs on multi-patch domains; however, their construction enforces only $C^0$ continuity of the solution across interfaces.
Second, defining a distance function near corners of the domain, such as polygonal vertices, is itself delicate and can lead to singular higher-order derivatives ~\cite{sukumar2022exact}. As a result, these regions are often either excluded or
down-weighted in the physics loss, which can compromise accuracy in the corner regions.
G\"oschel \textit{et al.}~\cite{goschel2025enforcing}
generalized the approach of Sukumar and Srivastava~\cite{sukumar2022exact} to
handle non-$C^1$ boundaries, but its extension to interface problems has not yet
been demonstrated.
\par

In this work, we investigate two hard-constrained PINN formulations for elliptic interface problems that represent fundamentally different ansatz strategies.
The first is a \textbf{windowing} formulation, which follows a multiplicative
hard-constraining strategy: interior, boundary, and interface subnetworks are
multiplied by simple, compactly-supported window functions so that the interface
conditions are embedded directly into the trial space.
Contrastingly, we propose a \textbf{buffer} formulation, which enforces these constraints through \textbf{additive} correction terms while leaving the neural network component largely unrestricted.
These additive corrections are determined from boundary and interface mismatches at discretely sampled points. \refB{The buffer method is therefore exact in one dimension and \emph{discretely} hard-constrained in higher dimensions}, enabling straightforward implementation on general geometries. Both hard-constrained formulations remove interface penalty terms from the loss and thereby decouple interface enforcement from PDE residual minimization, eliminating multi-objective loss balancing.
However, the buffer formulation exhibits optimization and approximation behavior
that differs fundamentally from that of typical multiplicative hard-constraining
approaches: instead of restricting the neural networks through window functions,
it preserves their flexibility by enforcing the constraints through a lightweight
additive correction.

Closely related additive hard-constraining approaches have been proposed very
recently~\cite{sukumar2026wachspress,dong2026novel}.
Sukumar and Roy~\cite{sukumar2026wachspress} enforced Dirichlet conditions
along the entire boundary, rather than only at discrete sample points, using
Wachspress-based transfinite coordinates.
Building on this idea, Dong and Zhang~\cite{dong2026novel} employed
TFC-transfinite constructions on mapped curved quadrilateral domains
to additionally enforce Neumann and Robin conditions.
Both formulations, however, are currently restricted to two-dimensional
domains, and extending them to higher dimensions with interface conditions
remains an open challenge.

The main contributions of this paper are therefore threefold. First, we formulate two ansatz-based hard-constrained PINN methods for elliptic interface problems: a windowing approach and a buffer approach.
Second, we show and analyze how multiplicative hard-constraining, represented
here by the windowing approach, can restrict neural-network expressivity through
the imposed window functions and thereby induce sensitivity to geometric overlap
and sharp corners, whereas the novel buffer approach is more robust because it does not restrict the interior subnetworks.
Third, we present comparative one- and two-dimensional benchmarks against soft-constrained PINN baselines and demonstrate that
the buffer formulation can provide a clear path toward robust and geometrically flexible practical implementations.

The remainder of the paper is organized as follows. Section \ref{sec:methodologies} introduces the model problems and the two hard-constrained formulations. Section \ref{sec:examples} presents one- and two-dimensional numerical experiments, compares the windowing and buffer approaches, and benchmarks them against soft-constrained methods. Section \ref{sec:conclusion} concludes with a summary of the main findings and directions for future work.
\section{Methodologies}\label{sec:methodologies}

\subsection{Problem Setup}
We consider a generic spatial domain composed of multiple disjoint subdomains,
\begin{equation}
\Omega = \bigcup_m^M \Omega_m,
\end{equation}
where each subdomain $\Omega_m$ has distinct physical properties. We denote the subdomain boundary of $\Omega_b$ as $\Gamma_b=\partial\Omega_b\cap\partial\Omega$,
and the interface between $\Omega_i$ and $\Omega_j$ as $\Gamma_{ij} = \partial\Omega_i\cap\partial\Omega_j$.
We also denote the set of the boundary of the domain
as $\rB=\{b | \partial\Omega_b\cap\partial\Omega\ne\varnothing\}$,
and the set of the interfaces between all subdomains as $\rI=\{(i,j) | \partial\Omega_i\cap\partial\Omega_j\ne\varnothing\}$.
Throughout this work, the following Poisson equation serves as a representative problem to illustrate the proposed methods
\begin{equation}\label{eq:poisson}
\LHS[u] \equiv -\nabla\cdot\kappa(\bx)\nabla u(\bx) = f(\bx)
\qquad
\forall\bx\in\Omega,
\end{equation}
where the diffusivity $\kappa(\bx)$ has discontinuities
at the subdomain interfaces.
We consider either Dirichlet or Neumann boundary conditions, partitioning
$\Gamma_b$ into $\Gamma_b^{(D)}$ and $\Gamma_b^{(N)}$, respectively.
On $\Gamma_b^{(D)}$, we impose
\begin{subequations}
\begin{equation}
u(\bx) = u_b(\bx)
\qquad
\forall\bx\in\Gamma_b^{(D)},
\end{equation}
whereas on $\Gamma_b^{(N)}$, we impose
\begin{equation}
\bn\cdot\nabla u(\bx) = \bn\cdot\nabla u_b(\bx)
\qquad
\forall\bx\in\Gamma_b^{(N)},
\end{equation}
\end{subequations}
where $\bn$ denotes the outward-facing normal vector on $\partial\Omega_b$.
\refB{Throughout this work we write $\bn_m$ for the outward-facing normal of subdomain
$\Omega_m$, so that $\bn_j=-\bn_i$ on a shared interface $\Gamma_{ij}$, and all interface
quantities are expressed in this orientation. For one-dimensional domains, $\bn_m$
reduces to a scalar $n_m=\pm1$.}
For the interface $\Gamma_{ij}$,
the physics state $u(\bx)$ and its flux are continuous across the interfaces,
\begin{subequations}\label{eq:itf_cond}
\begin{equation}
\jump{u}_{ij} \equiv u\big|_{\Omega_i} - u\big|_{\Omega_j} = 0
\qquad
\forall\bx\in\Gamma_{ij}
\end{equation}
\begin{equation}
\avg{\bn\cdot\kappa\nabla u}_{ij}
\equiv \refB{\frac{\bn_i\cdot\kappa\nabla u\big|_{\Omega_i} + \bn_j\cdot\kappa\nabla u\big|_{\Omega_j}}{2}}
= 0,
\qquad
\forall\bx\in\Gamma_{ij}.
\end{equation}
\end{subequations}
\par
While we consider here the case of continuous solution and flux across the interface, the formulations developed in this work are not limited to zero-jump conditions and can be extended naturally to interfaces with prescribed nonzero jumps. Throughout this study, hard-constrained PINN models are trained against the residual loss
for the physics governing equation,
\begin{equation}\label{eq:J-int}
\cJ_{physics} = \sum_{k}^K ( -\nabla\cdot\kappa(\bx_k)\nabla u(\bx_k) - f(\bx_k) )^2,
\end{equation}
where $\{\bx_k\}_{k=1}^K\subset\Omega$ are interior collocation points.
For soft-constrained methods, residual losses against the boundary and
interface conditions are evaluated at $K_{bnd}$ boundary collocation points, in
addition to the physics loss.
\begin{equation}\label{eq:J-dbc_1}
\cJ_{dbc} = \sum_{b\in\rB}\sum_{\bx_k\in\Gamma_b^{(D)}}( u(\bx_k) - u_b(\bx_k) )^2,
\end{equation}
\begin{equation}\label{eq:J-dbc_2}
\cJ_{nbc} = \sum_{b\in\rB}\sum_{\bx_k\in\Gamma_b^{(N)}} ( \bn\cdot \nabla u(\bx_k) - \bn\cdot\nabla u_b(\bx_k) )^2,
\end{equation}
\begin{equation}\label{eq:J-dbc_3}
\cJ_{int} = \sum_{(i,j)\in\rI}\sum_{\bx_k\in\Gamma_{ij}} ( u(\bx_k)\big|_{\Omega_i}  - u(\bx_k)\big|_{\Omega_j}  )^2,
\end{equation}
\begin{equation}\label{eq:J-dbc_4}
\cJ_{fint} = \sum_{(i,j)\in\rI}\sum_{\bx_k\in\Gamma_{ij}} \refB{( \bn_i\cdot \kappa \nabla u(\bx_k)\big|_{\Omega_i} + \bn_j\cdot \kappa \nabla u(\bx_k)\big|_{\Omega_j} )^2},
\end{equation}
We study two alternatives to penalty-based interface enforcement. The first is an exact hard-constrained windowing ansatz. The second is a buffer-correction strategy that is exact in 1D and discretely hard-constrained in higher dimensions.
\subsection{Windowing approach}\label{subsec:windowing}
In this section, we describe the \textit{windowing approach}, whose objective is to automatically enforce interface and boundary constraints within the PINN framework.
This is achieved through a solution ansatz written as a sum of interior, boundary, and interface contributions,
\begin{subequations}\label{eq:window-ansatz}
\begin{equation}
u_{NN}(\bx;\theta)
=
\sum_m^M u^{(S)}_m(\bx;\theta_m)
+
\sum_{b\in\rB} u^{(B)}_b(\bx;\theta_b)
+
\sum_{(i,j)\in\rI} u^{(I)}_{ij}(\bx;\theta_{ij}),
\label{eq:widnow_approach}
\end{equation}
where $\theta_m$, $\theta_b$, and $\theta_{ij}$ denote the neural-network
parameters associated with the interior, boundary, and interface terms,
respectively, and $\theta$ denotes their union.
The interior contribution $u^{(S)}_m$ provides the degrees of freedom used to represent the PDE solution within each subdomain.
Here, each contribution is a product of neural networks with so-called window functions,
\begin{equation}\label{eq:uD-m}
u^{(S)}_m(\bx;\theta_m) = W_m(\bx)NN_m(\bx;\theta_m).
\end{equation}
The boundary contribution $u^{(B)}_b$ encodes the prescribed boundary value or normal derivative,
also being products of a window function with a trainable function,
\begin{equation}\label{eq:uB-m}
u^{(B)}_b(\bx;\theta_b) = W_{b,d}(\bx) g_{b,d}(\bx;\theta_{b,d}) + W_{b,n}(\bx) g_{b,n}(\bx;\theta_{b,n}).
\end{equation}
Similarly, the interface contribution $u^{(I)}_{ij}$ carries the shared interface value and flux information so that continuity and flux balance are enforced by construction,
\begin{equation}\label{eq:uI-m}
u^{(I)}_{ij}(\bx;\theta_{ij}) = W_{ij,d}(\bx) g_{ij,d}(\bx;\theta_{ij,d}) + W_{ij,n}(\bx) \kappa^{-1}(\bx)g_{ij,n}(\bx;\theta_{ij,n}).
\end{equation}
\end{subequations}
The trainable functions $g_{b,d}$, $g_{b,n}$, $g_{ij,d}$ and $g_{ij,n}$
either prescribe the boundary/interface conditions
or neural networks defined only on the boundary/interface,
which will be introduced subsequently.
The purpose of window functions is to confine neural networks
to the neighborhood of the corresponding interior/boundary/interface of the domains.
In theory, any architecture can be used for neural network functions $NN(\bx;\theta)$. In this study, for simplicity, we only consider simple fully connected multi-layer perceptrons. Each term with the corresponding window is first introduced using a one-dimensional space example,
and extension to higher dimensions will be discussed subsequently.

\subsubsection{Window functions for one-dimensional problems}
The window functions are designed so that each component of the ansatz controls only the constraint it is intended to impose. Interior windows vanish at the edges of each subdomain, ensuring that the interior subnetworks do not interfere with boundary or interface conditions. Boundary windows are then chosen according to the type of prescribed data: Dirichlet windows enforce the solution value while having no constraining effect on the normal derivative, whereas Neumann windows enforce the normal derivative while leaving the solution value unconstrained. Interface windows follow the same principle, allowing the interface terms to carry the shared solution and flux information across neighboring subdomains. As a result, continuity and flux balance are embedded directly into the ansatz rather than enforced through penalty terms.

The interior neural network, $NN_m: \Omega_m\to\rR$, is restricted by the interior window function $W_m(\bx)$.
For a one-dimensional domain, each interior window function is located at the subdomain center $x_m$ with its size $\Delta x_m$,
\begin{equation}\label{eq:tW-int}
W_m(x) = \tW_{int}\!\left(\frac{\vert x - x_{m}\vert}{\Delta x_{m}}\right),
\end{equation}
where $\tW_{int}(\tau)$ represents the window function on the normalized spatial variable $\tau = \frac{|x-x_m|}{\Delta x_m}$.
This idea of constructing a solution ansatz from multiple neural networks, each multiplied by an auxiliary window function, was originally introduced by Moseley et al.~\cite{moseley2023finite} in the Finite Basis PINN (FB-PINN) framework. Their work focused on capturing high frequency solutions, and employed smoothly decaying window functions.
In the present work, however, our goal is not frequency localization but the exact embedding of boundary and interface constraints into the ansatz. For this reason, we instead propose to use polynomial window functions, since they can be constructed to satisfy prescribed value and derivative conditions in a simple and explicit manner.
\par
For the 1D problem, the window functions have to satisfy the following boundary conditions,
\begin{equation}\label{eq:wint-bc}
\begin{split}
\tW_{int}(0) = 1
\qquad
\tW_{int}(1) = 0\\
\tW'_{int}(0) = 0
\qquad
\tW'_{int}(1) = 0.
\end{split}
\end{equation}
The function and first derivatives of $\tW_{int}(\tau)$ vanish at its boundary $\tau=1$.
This choice ensures that the interior neural networks $NN_m$ do not contribute to the solution value or gradient at the boundary or interface,
thereby allowing the boundary and interface conditions to be imposed exactly through the remaining ansatz terms.
For higher-order PDEs where boundary or interface conditions involve derivatives beyond first-order,
the window function must be constructed such that its higher-order derivatives $\frac{d^n}{d\tau^n}\tW_{int}$ also vanish at $\tau=1$.
\par
Technically, $\tW'_{int}(0) = 0$ is not required by the boundary or interface
conditions, but by the physics equation itself.
With $\tW'_{int}(0) \ne 0$, the second-order derivative
$\frac{d^2}{dx^2}\tW(|x-x_m|)$ becomes discontinuous at the subdomain center,
thereby introducing an unintended discontinuity into the Poisson equation.
By the same reasoning, its higher-order derivatives
$\frac{d^n}{d\tau^n}\tW_{int}(0)$ should vanish for higher-order PDEs.
\par
For the present second-order setting, the lowest-degree polynomial that satisfies (\ref{eq:wint-bc}) is a cubic polynomial~\cite{roy2024exact},
\begin{equation}
\tW_{int}(\tau) = 1 - 3\tau^2 + 2\tau^3.
\end{equation}
Although the cubic polynomial starts as a simple baseline construction, higher-order polynomial windows may, however, offer greater flexibility and improved training behavior for some problems. Their relative performance is investigated in ~\ref{subsec:comparison-window-order}.\todo{Kevin: section for polynomial order comparison}
\par
Along the external boundary, the ansatz must incorporate the prescribed boundary value while retaining flexibility in the normal derivative. To achieve this, the boundary value of the solution ansatz $u_{NN}$ is represented using two
boundary functions $g_{b,d}, g_{b,n}:\partial\Omega_b\cap\partial\Omega\to\rR$.
For Dirichlet boundary conditions,
$g_{b,d}$ enforces the prescribed boundary condition
while $g_{b,n}$ allows a free trainable scalar parameter $\theta_{b,n}$ for the boundary slope,
\begin{equation}\label{eq:gbd}
g_{b,d}(x) = u_b, \;
g_{b,n}(x) = \theta_{b,n}.
\end{equation}
Similarly, for Neumann conditions,
$g_{b,n}$ enforces the prescribed boundary condition
while $g_{b,d}$ allows a free trainable scalar parameter $\theta_{b,d}$ for the boundary value,
\begin{equation}\label{eq:gbn}
g_{b,d}(\bx) = \theta_{b,d}, \;
g_{b,n}(\bx) = \Dpartial{u_b}{x}.
\end{equation}
\par
To ensure that $g_{b,d}$ and $g_{b,n}$ enforce Dirichlet and Neumann conditions independently, each is multiplied by the corresponding boundary window function similar to (\ref{eq:tW-int}),
\begin{subequations}\label{eq:Wb}
\begin{equation}
W_{b,d}(x) = \tW_{d}\!\left(\frac{\vert x - x_{b}\vert}{\Delta x_{b}}\right), 
\end{equation}
\begin{equation}
W_{b,n}(x) = -\mathrm{sign}(\bn)\Delta x_{b}\tW_{n}\!\left(\frac{\vert x - x_{b}\vert}{\Delta x_{b}}\right).
\end{equation}
\end{subequations}
Along the normal direction, $x_{b}$ denotes the boundary location itself rather than the center of the subdomain. The prefactor $(-\mathrm{sign}(\bn)\Delta x_{b})$ in $W_{b,n}$
is included to cancel the chain-rule factor arising from the derivative of $\tW_{n}$ and preserve the correct sign with respect to the boundary normal. To ensure that $g_{b,d}$ only enforces the Dirichlet condition and does not affect the Neumann condition, the associated window function
$\tW_d$ must satisfy $\tW_d(0)=1$ and $\tW'_d(0) = 0$.
Additionally, both its value and first derivative must vanish at $\tau=1$, so that this term has no influence on other boundaries or interfaces. These requirements on $\tW_d$ are identical to those in (\ref{eq:wint-bc}) for $\tW_{int}$.
By contrast, the Neumann window function must enforce the derivative condition without altering the boundary value, and is therefore required to satisfy
\begin{equation}\label{eq:wn-bc}
\begin{split}
\tW_{n}(0) = 0
\qquad
\tW_{n}(1) = 0\\
\tW'_{n}(0) = 1
\qquad
\tW'_{n}(1) = 0,
\end{split}
\end{equation}
so that it only contributes to the Neumann condition but not to the Dirichlet condition.
\par
Having defined the boundary windows so that value and derivative constraints are cleanly separated, we use the same principle at internal interfaces. Unlike the boundary case, at the interface $\Gamma_{ij}$,
the values of the solution $u$ and flux $\kappa\nabla u$ are not prescribed individually, but rather their differences across the interface as described in (\ref{eq:itf_cond}). The
interface functions $g_{ij,d}, g_{ij,n}:\Gamma_{ij}\to\rR$
therefore represent the shared interface value and flux information,
with free trainable scalar parameters $\theta_{ij,d}, \theta_{ij,n}$ controlling their contribution to the ansatz:
\begin{subequations}\label{eq:gij}
\begin{equation}
g_{ij,d}(x) \equiv u_{NN}(x)\bigg|_{\Gamma_{ij}} = \theta_{ij,d}
\end{equation}
\begin{equation}\label{eq:g-ij-n}
g_{ij,n}(x) \equiv \kappa(x)\Dpartial{u_{NN}}{x}(x)\bigg|_{\Gamma_{ij}} = \theta_{ij,n}.
\end{equation}
\end{subequations}
Since $g_{ij,n}$ in (\ref{eq:g-ij-n}) determines the solution flux $\kappa\Dpartial{u}{x_{k_{ij}}}$,
not the slope $\Dpartial{u}{x_{k_{ij}}}$,
a factor $\kappa^{-1}$ is further multiplied in the solution ansatz (\ref{eq:window-ansatz}).
This way, the window functions for the Dirichlet and Neumann conditions can be reused
for the interface window functions,
\begin{subequations}\label{eq:Wij}
\begin{equation}
W_{ij,d}(x) = \tW_{d}\!\left(\frac{\vert x - x_{ij}\vert}{\Delta x_{ij}}\right)
\end{equation}
\begin{equation}
W_{ij,n}(x) = -\mathrm{sign}(\bn)\Delta x_{ij}\tW_{n}\!\left(\frac{\vert x - x_{ij}\vert}{\Delta x_{ij}}\right),
\end{equation}
\end{subequations}
with the interface subdomain center $x_{ij}$ located right on the interface and its size $\Delta x_{ij}$.
\par
While the interface condition (\ref{eq:itf_cond}) only describes a continuous solution and flux,
discontinuities in the solution and flux can also be handled with
discontinuous interface functions.
For example, for the following discontinuous interface condition,
\begin{subequations}\label{eq:itf_cond_jump}
\begin{equation}
\jump{u}_{ij} = h_{ij,d},
\qquad
x=x_{ij}
\end{equation}
\begin{equation}
\avg{\bn\cdot\kappa\nabla u}_{ij} = h_{ij,n},
\qquad
x=x_{ij},
\end{equation}
\end{subequations}
the interface functions can be redefined as,
\begin{subequations}
\begin{equation}
g_{ij,d}(\bx) =
\begin{cases}
\theta_{ij,d} & x < x_{ij} \\
\theta_{ij,d} + h_{ij,d} & x > x_{ij}
\end{cases}
\end{equation}
\begin{equation}
g_{ij,n}(x) =
\begin{cases}
\theta_{ij,n} & x < x_{ij} \\
\theta_{ij,n} + h_{ij,n} & x > x_{ij}.
\end{cases}
\end{equation}
\end{subequations}
\par
The equations (\ref{eq:tW-int}), (\ref{eq:gbd}), (\ref{eq:gbn}), (\ref{eq:Wb}), (\ref{eq:gij}) and (\ref{eq:Wij})
constitute the solution ansatz (\ref{eq:window-ansatz}) for the windowing approach.

\subsubsection{Admissible polynomial window functions}\label{sec:window_functions}
The boundary and interface constraints introduced above admit multiple polynomial window functions. Although all such choices satisfy the required hard-constraint conditions, they differ in their higher-order derivative structure, which can influence optimization and residual representation. We therefore consider a small family of admissible polynomial windows with increasing smoothness and compare their behavior in the numerical experiments.
\begin{figure}[tbph]
    \begin{tikzpicture}[
]
    \begin{groupplot}[
        group style={
            group name = my plots,
            group size= 3 by 1,
            xlabels at =edge bottom,
            horizontal sep=2cm,
            vertical sep=2.2cm,
        },
        height = 0.35\textwidth,
        width = 0.35\textwidth,
        enlarge x limits={false, abs value = 5mm},
        enlarge y limits={false, abs value = 5mm},
        xtick={0, 0.5, 1.},
        name=chung,
    ]    
\pgfplotsset{set layers=standard}%

        \nextgroupplot[
            ymax=1.1,
            xlabel={$\tau$},
            ylabel={$\tW^{(k)}_{int}(\tau)$},
            tick scale binop ={\times},
            legend style={
                draw=none, fill=none,
                at={(rel axis cs: 0, 1.0)},
                anchor=south west,
                nodes={scale=1.0},
                legend cell align={left},
                legend columns=3,
                /tikz/every even column/.append style={column sep=0.5cm},
            },
        ]

            \addplot+ [
                line width=1.0,
                solid,
                mark=none,
                blue,
            ]
            table [
                x index=0, y index=1,
            ]{data_interior_window.txt};

            \addplot+ [
                line width=1.0,
                solid,
                mark=none,
                brown,
            ]
            table [
                x index=0, y index=2,
            ]{data_interior_window.txt};

            \addplot+ [
                line width=1.0,
                solid,
                mark=none,
                red,
            ]
            table [
                x index=0, y index=3,
            ]{data_interior_window.txt};

            \legend{$k=1$, $k=2$, $k=3$};

        \nextgroupplot[
            ymax=1.1,
            xlabel={$\tau$},
            ylabel={$\tW^{(k)}_{d}(\tau)$},
            tick scale binop ={\times},
        ]

            \addplot+ [
                line width=1.0,
                solid,
                mark=none,
                blue,
            ]
            table [
                x index=0, y index=1,
            ]{data_dirichlet_window.txt};

            \addplot+ [
                line width=1.0,
                solid,
                mark=none,
                brown,
            ]
            table [
                x index=0, y index=2,
            ]{data_dirichlet_window.txt};

            \addplot+ [
                line width=1.0,
                solid,
                mark=none,
                red,
            ]
            table [
                x index=0, y index=3,
            ]{data_dirichlet_window.txt};

        \nextgroupplot[
            ymax=0.5,
            xmin=0., xmax=1.,
            xlabel={$\tau$},
            ylabel={$\tW^{(k)}_{n}(\tau)$},
            tick scale binop ={\times},
        ]

            \addplot+ [
                dashed,
                line width=1.0,
                gray,
                mark=none,
            ] {x};

            \addplot+ [
                line width=1.0,
                solid,
                mark=none,
                blue,
            ]
            table [
                x index=0, y index=1,
            ]{data_neumann_window.txt};

            \addplot+ [
                line width=1.0,
                solid,
                mark=none,
                brown,
            ]
            table [
                x index=0, y index=2,
            ]{data_neumann_window.txt};

            \addplot+ [
                line width=1.0,
                solid,
                mark=none,
                red,
            ]
            table [
                x index=0, y index=3,
            ]{data_neumann_window.txt};
            \logLogSlopeTriangle{0.45}{0.1}{0.6}{1}{}{north};

  \end{groupplot}
\node[below = 1.5cm of my plots c1r1.south west,
    anchor=west,
] {(a)};
\node[below = 1.5cm of my plots c2r1.south west,
    anchor=west,
] {(b)};
\node[below = 1.5cm of my plots c3r1.south west,
    anchor=west,
] {(c)};
\end{tikzpicture}
%
    \caption{Admissible polynomial window functions considered in this work: (a) \refB{interior} window functions, (b) Dirichlet window functions, (c) Neumann window functions. Here, $k$ denotes the maximum vanishing derivative order.
    }
    \label{fig:windows}
\end{figure}
\begin{table}[tbh]
\begin{tabular}{|c|c|c|c|}
\hline
$\tW_{int}$ & $\tau=0$ & $\tau=1$ & analytic form \\
\hline
$\tW^{(1)}_{int}$ & $\tW^{(1)}_{int}=1$, $\partial_{\tau}\tW^{(1)}_{int}=0$ & \multirow{3}{6em}{$\partial^n_{\tau}\tW_{int}=0$ for $n=0,1$} & $1 - 3\tau^2 + 2\tau^3$ \\
\cline{1-2} \cline{4-4}
$\tW^{(2)}_{int}$ & $\tW^{(2)}_{int}=1$, $\partial^n_{\tau}\tW^{(2)}_{int}=0$ for $n=1,2$ &  & $1 - 4\tau^3 + 3\tau^4$ \\
\cline{1-2} \cline{4-4}
$\tW^{(3)}_{int}$ & $\tW^{(3)}_{int}=1$, $\partial^n_{\tau}\tW^{(3)}_{int}=0$ for $n=1,2,3$ &  & $1 - 5\tau^4 + 4\tau^5$ \\
\hline
\end{tabular}
\caption{Admissible interior window functions and their constraints.}
\label{tab:wint}
\end{table}
\begin{table}[tbh]
\begin{tabular}{|c|c|c|c|}
\hline
$\tW_{d}$ & $\tau=0$ & $\tau=1$ & analytic form \\
\hline
$\tW^{(1)}_{d}$ & \multirow{3}{5em}{$\tW_{d}=1$, $\partial_{\tau}\tW_{d}=0$} & $\partial^n_{\tau}\tW_{d}=0$ for $n=0,1$ & $1 - 3\tau^2 + 2\tau^3$ \\
\cline{1-1} \cline{3-4}
$\tW^{(2)}_{d}$ & & $\partial^n_{\tau}\tW_{d}=0$ for $n=0,1,2$ & $1 - 6\tau^2 + 8\tau^3 - 3\tau^4$ \\
\cline{1-1} \cline{3-4}
$\tW^{(3)}_{d}$ & & $\partial^n_{\tau}\tW_{d}=0$ for $n=0,1,2,3$ & $1 - 10\tau^2 + 20\tau^3 - 15\tau^4 + 4\tau^5$ \\
\hline
\end{tabular}
\caption{Admissible Dirichlet window functions and their constraints.}
\label{tab:wd}
\end{table}
\begin{table}[tbh]
\begin{tabular}{|c|c|c|c|}
\hline
$\tW_{n}$ & $\tau=0$ & $\tau=1$ & analytic form \\
\hline
$\tW^{(1)}_{n}$ & \multirow{3}{5em}{$\tW_{n}=0$, $\partial_{\tau}\tW_{n}=1$} & $\partial^n_{\tau}\tW_{n}=0$ for $n=0,1$ & $\tau(1-\tau)^2$ \\
\cline{1-1} \cline{3-4}
$\tW^{(2)}_{n}$ & & $\partial^n_{\tau}\tW_{n}=0$ for $n=0,1,2$ & $\tau(1-\tau)^3$ \\
\cline{1-1} \cline{3-4}
$\tW^{(3)}_{n}$ & & $\partial^n_{\tau}\tW_{n}=0$ for $n=0,1,2,3$ & $\tau(1-\tau)^4$ \\
\hline
\end{tabular}
\caption{Admissible Neumann window functions and their constraints.}
\label{tab:wn}
\end{table}

Figure~\ref{fig:windows} shows the candidate polynomial window functions investigated in this study.
Their specific forms and constraints are also summarized in Tables~\ref{tab:wint}-\ref{tab:wn}.
\par
First, for the interior window function $\tW^{(k)}_{int}$,
vanishing derivatives in the center of the subdomain $\tau=0$ are added with increasing derivative order,
\begin{equation}
\DpartN{\tW_{int}^{(k)}}{\tau}{n}\bigg|_{\tau=0} = 0,
\qquad
\forall n=1, 2, \ldots, k,
\end{equation}
with the superscript $(k)$ indicating the maximum enforced vanishing derivative order.
With increasing $k$, the interior window $\tW^{(k)}_{int}$ becomes smoother in the center of the subdomain, where no discontinuity is expected in the solution or in the physics equation.
On the other hand, at the subdomain boundary $\tau=1$,
the derivatives higher than the first order remain non-zero.
This is necessary to keep the second derivative of the interior solution term from vanishing,
\begin{equation}
\Dparttwo{}{\tau}\left(\tW_{int}(\tau)NN(\tau;\theta)\right)
=
\Dparttwo{\tW_{int}}{\tau}NN
+
\Dpartial{\tW_{int}}{\tau}\Dpartial{NN}{\tau}
+
\tW_{int}\Dparttwo{NN}{\tau}
\ne 0,
\end{equation}
as $\tau\to1$,
so that the solution ansatz \eqref{eq:window-ansatz} can resolve the physics equation \eqref{eq:poisson} near the interface or the boundary.
\par
For the boundary/interface subdomains, unlike the interior window,
$\tau=1$ is located in the middle of the physical domain
where no discontinuity is expected in the solution or in the physics equation.
Thus, for the Dirichlet window function $\tW^{(k)}_d$
and the Neumann window function $\tW^{(k)}_n$,
the vanishing derivative conditions are added at the subdomain edge $\tau=1$,
\begin{equation}
\DpartN{\tW_{d}^{(k)}}{\tau}{n}\bigg|_{\tau=1} =
\DpartN{\tW_{n}^{(k)}}{\tau}{n}\bigg|_{\tau=1} = 0,
\qquad
\forall n=1, 2, \ldots, k.
\end{equation}
\par
We note that other forms can be used for the window function,
such as the sigmoid function in Mosely \textit{et al.}~\cite{moseley2023finite}.
This function can be considered as (negligibly) vanishing derivatives of infinite order.
However, a smoother function does not necessarily perform better in training convergence.
This will be shown in ~\ref{subsec:comparison-window-order}.

\subsubsection{Extension to higher dimensions}
\label{subsec:window-multi-dim}

The same design principles used in one dimension can be carried over to higher-dimensional settings. As a simple and natural extension, we consider a 
$D$-dimensional hyper-rectangular subdomain.
The corresponding window function is constructed as the product of one-dimensional window functions,
\begin{equation}
W_m(\bx) = \prod_{k=1}^D\tW_{int}\!\left(\frac{\vert x_k - x_{m,k}\vert}{\Delta x_{m,k}}\right),
\end{equation}
where the subdomain center is $\bx_m =(x_{m,1}, \ldots, x_{m,D})$
and the subdomain size is $\Delta \bx_m =(\Delta x_{m,1}, \ldots \Delta x_{m,D})$.
Similarly, the boundary window functions are defined as
\begin{equation}
W_{b,d}(\bx) = \tW_{d}\!\left(\frac{\vert x_{k_b} - x_{b,k_b}\vert}{\Delta x_{b,k_b}}\right)
\prod_{k\ne k_b}\tW_{int}\!\left(\frac{\vert x_k - x_{b,k}\vert}{\Delta x_{b,k}}\right)
\end{equation}
\begin{equation}
W_{b,n}(\bx) = -\mathrm{sign}(n_{k_b})\Delta x_{b,k_b}\tW_{n}\!\left(\frac{\vert x_{k_b} - x_{b,k_b}\vert}{\Delta x_{b,k_b}}\right)
\prod_{k\ne k_b}\tW_{int}\!\left(\frac{\vert x_k - x_{b,k}\vert}{\Delta x_{b,k}}\right),
\end{equation}
where $k_b$ denotes the coordinate direction normal to the boundary.
The interface window functions can be constructed in an analogous manner.
\par
This outer-product form of the window enforces the boundary condition along one
edge and vanishes at corners where two boundary or interface surfaces meet.
Accordingly, additional boundary or interface enforcement is required at the
domain corners.
As an illustrative example,
consider a two-dimensional corner at $\bx_c=(x_{c1}, x_{c2})$,
where the two boundaries $\Gamma_1=\{\bx| x_1=x_{c1}\}$ and $\Gamma_2=\{\bx| x_2=x_{c2}\}$ meet.
In this case, we augment the solution ansatz \eqref{eq:window-ansatz} with a corner contribution,
\begin{equation}\label{eq:window-corner}
\begin{split}
u^{(C)}(\bx;\theta)
&= W^{(C)}_{1,d}(\bx)g_{1,d}(\bx;\theta_{1,d}) + W^{(C)}_{1,n}(\bx) g_{1,n}(\bx;\theta_{1,n})\\
&+ W^{(C)}_{2,d}(\bx)g_{2,d}(\bx;\theta_{2,d}) + W^{(C)}_{2,n}(\bx) g_{2,n}(\bx;\theta_{2,n}),
\end{split}
\end{equation}
where $g_{1,d}$ and $g_{1,n}$ correspond to the boundary functions for $\Gamma_1$ ($b=1$),
and $g_{2,d}$ and $g_{2,n}$ correspond to those for $\Gamma_2$ ($b=2$).
\par
For the corner window functions $W^{(C)}_{\cdot,\cdot}$,
the Cartesian outer-product construction
is generally insufficient to satisfy the boundary/interface conditions at both sides,
particularly with inhomogeneous boundary conditions.
Instead, in this study, the corner windows are defined in polar coordinates centered at the corner,
\begin{subequations}\label{eq:wc}
\begin{equation}
W^{(C)}_{b,d}(\bx)
= \tW_d\!\left(\frac{|\alpha - \alpha_b|}{\Delta \alpha}\right)\tW_{int}\!\left(\frac{r}{\Delta r}\right)
\end{equation}
\begin{equation}
W^{(C)}_{b,n}(\bx)
= -\mathrm{sign}(\bn)\Delta\alpha\,\tW_n\!\left(\frac{|\alpha - \alpha_b|}{\Delta \alpha}\right)\tW_{int}\!\left(\frac{r}{\Delta r}\right),
\end{equation}
where $(r, \alpha)$ is the polar coordinates relative to the corner, defined by
\begin{equation}
\bx = \bx_c + (r\cos\alpha,r\sin\alpha),
\end{equation}
\end{subequations}
and $\alpha_b$ denotes the polar angle associated with boundary $\Gamma_b$, with $\Delta \alpha = |\alpha_2 - \alpha_1|$.
The sign of the outward normal vector $\mathrm{sign}(\bn)$ is defined along the azimuthal direction.
\par
In one-dimensional settings,
the boundary and interface functions defined in \eqref{eq:gbn}, \eqref{eq:gbd} and \eqref{eq:gij}
are parameterized by scalar trainable variables.
This is sufficient because boundaries and interfaces reduce to points,
and the associated boundary/interface conditions are scalar quantities
defined at those points, without spatial variation.
\par
In higher dimensions, however,
boundaries and interfaces form curves (in 2D) or surfaces (in 3D),
and the boundary/interface conditions generally exhibit spatial dependence along these manifolds.
Consequently, the scalar trainable parameters must be replaced by neural networks capable of representing spatial variation along the boundary or interface.
\par
Specifically, for \eqref{eq:gbd}, we replace the scalar parameter with
\begin{equation}
g_{b,n}(\bx) = NN_{b,n}(\bx;\theta_{b,n}),
\end{equation}
and for \eqref{eq:gbn},
\begin{equation}\label{eq:gbd-highd}
g_{b,d}(\bx) = NN_{b,d}(\bx;\theta_{b,d}).
\end{equation}
Similarly, for the interface functions in \eqref{eq:gij}, we define
\begin{subequations}
\begin{equation}
g_{ij,d}(\bx) =
NN_{ij,d}(\bx;\theta_{ij,d})
\end{equation}
\begin{equation}
g_{ij,n}(\bx) =
NN_{ij,n}(\bx;\theta_{ij,n})
\end{equation}
\end{subequations}
These neural networks are defined only on the boundary or interface manifold and are independent of the normal coordinate direction $x_{k_b}$.
In other words, they vary tangentially along the boundary/interface but remain constant in the normal direction.
We note that many real systems will involve geometries much more complicated than what is described in this section. The windowing approach can be applied to non-orthogonal geometries (see, e.g. \ref{app:2dpoisson-full-window}), but the formation of window functions for complicated geometries is not a trivial task and presents a challenge for practical usage. 
\subsubsection{Discussion on different forms of windowing approach}
\todo[inline]{Kevin: the optimal location of this section may not be here. Pratanu: One option is to move it to the Results/Discussion section.}
We recognize that different forms of windowing approaches have been widely used
for hard-constraining boundary conditions.
For example, Moseley \textit{et al.}~\cite{moseley2023finite} proposed to use the form
\begin{equation}
u_{NN}(\bx;\theta)
=
g(\bx) + W_{int}(\bx)NN(\bx;\theta),
\end{equation}
where $g(\bx):\Omega\to\rR$ has the prescribed boundary condition.
The major difference in this approach is that
the interior neural networks determine both the interior solution
and the non-prescribed boundary value.
For example, for the one-dimensional Poisson problem with Dirichlet boundary conditions,
\begin{subequations}
\begin{equation}
\Dparttwo{u}{x} = f(x)
\qquad
\forall x\in[0, 1]
\end{equation}
\begin{equation}
u(0) = u(1) = 0,
\end{equation}
\end{subequations}
the window function $W_{int}(x) = x(1-x)$ can be used,
\begin{equation}
u_{NN}(x;\theta) = x(1-x)NN(x;\theta).
\end{equation}
In such a case, there is no need for a free trainable function for the boundary slope $\Dpartial{u}{x}$.
\par
However, this ansatz---constructed as a product of a window function and a neural network---is not applicable for general boundary conditions,
particularly for Neumann conditions.
For such an ansatz,
the derivative along the normal direction (say, the $k_b$-th axis) becomes
\begin{equation}\label{eq:n-dot-grad-WNN}
\bn\cdot\nabla (W(\bx)NN(\bx;\theta))
=\mathrm{sign}(n_{k_b})\left[ \Dpartial{W}{x_{k_b}}(\bx)NN(\bx;\theta) + W(\bx)\Dpartial{NN}{x_{k_b}}(\bx;\theta) \right].
\end{equation}
To enforce Neumann boundary conditions in a hard-constrained manner,
the required boundary slope must be prescribed by the function $g(\bx)$, while the neural network contribution to the boundary slope must vanish.
However, since the boundary value itself is not prescribed, $W(\bx)NN(\bx;\theta)$ should remain nonzero at the boundary to represent the free boundary solution value.
This requires $W(\bx)\ne0$ at the boundary.
However, if $W(\bx)\ne 0$,
the second term in \eqref{eq:n-dot-grad-WNN} introduces a nonzero contribution from the neural network to the boundary slope, thereby violating the hard constraint on the Neumann condition.
\par
The only case in which the boundary slope in \eqref{eq:n-dot-grad-WNN} vanishes while $W(\bx)\ne0$
is when $NN(\bx;\theta)$ is independent of $x_{k_b}$,
so that $\partial NN/\partial x_{k_b}=0$.
In that case, the neural network effectively reduces to a free trainable boundary function \eqref{eq:gbd-highd} proposed in this study.
Otherwise, a different form of solution ansatz is required to hard-constrain the Neumann condition.
Similar reasoning applies to more general conditions, including Robin boundary conditions and interface conditions.
\par
Notably, Sukumar and Srivastava proposed an alternative formulation for Neumann condition,
in which a differential operator is applied to neural networks~\cite[Eq. 24]{sukumar2022exact},
\begin{equation}
u_{NN}(\bx;\theta) =
[1 - \nabla\phi\cdot\nabla]NN_1(\bx;\theta)
+\phi\Dpartial{u}{\bn}\bigg|_{\partial\Omega}
+\phi^2NN_2(\bx;\theta),
\end{equation}
where $\phi(\bx)$ is an approximate distance function, playing a role analogous to a window function. However, this method suffers from singular behavior of the distance function associated with boundary corners.

These considerations show that, while alternative windowing constructions are possible, designing a fully hard-constrained masking ansatz that simultaneously preserves flexibility and enforces boundary and interface conditions in a clean and localized manner is not straightforward. The product-window formulation adopted here provides a systematic way to separate interior, boundary, and interface roles, but it also introduces geometric and derivative-coupling constraints that can complicate the construction, especially in higher dimensions. This motivates the introduction of an alternative formulation that retains strong control of boundary and interface conditions while avoiding direct restrictions on the interior neural-network representation.
\subsection{Buffer approach}
We next introduce a buffer-based formulation, in which the neural-network solution is left unrestricted and boundary/interface conditions are enforced via the addition of auxiliary correction terms. In this approach, each subdomain neural network is augmented by a corresponding boundary buffer function $g_m(\bx)$,
\begin{equation}\label{eq:buffer}
u_{NN}(\bx;\theta)
=
\sum_m^M\left[ NN_m(\bx;\theta_m) + g_m(\bx;\bc(\theta_m)) \right].
\end{equation}
The key difference from the windowing approach is
that the neural networks are not restrained by any window function, and therefore generally retain non-zero contributions at boundaries and interfaces.
Instead, the buffer function $g_m$
can adjust its degrees of freedom $\bc(\theta_m)$ (DOFs) for a given network parameter set $\theta_m$,
so as to cancel the boundary/interface mismatch introduced by the neural networks
and satisfy the required boundary/interface conditions.
In this respect, it is conceptually similar to lifting-based methods, which enforce boundary conditions through additive corrections rather than multiplicative restriction of the neural network.
\par
We introduce this approach with a few examples of gradually increasing complexity in boundary and interface conditions.

\subsubsection{One-dimensional, continuous problem with boundary conditions}
For example, for a continuous, one-dimensional Poisson problem (\ref{eq:poisson})
with Neumann-Dirichlet boundary conditions
\begin{equation}\label{eq:1d-ND-bc}
\partial_x u(0) = \partial_x u_0, 
\qquad
u(1) = u_1,
\end{equation}
the one subdomain solution ansatz is
\begin{equation}\label{eq:buffer1}
u_{NN}(x;\theta) = NN(x;\theta) + g(x).
\end{equation}
The buffer function $g(x)$ then must satisfy the boundary conditions at the given parameter $\theta$,
\begin{subequations}
\begin{equation}
\partial_xg(0) = \partial_xu_0 - \partial_xNN(0;\theta)
\end{equation}
\begin{equation}
g(1) = u_1 - NN(1;\theta).
\end{equation}
\end{subequations}
A linear function is the simplest form for $g(x)$ to satisfy these boundary conditions,
\begin{equation}
g(x) = c_0 + c_1x,
\end{equation}
whose coefficients $\bc=\{c_0, c_1\}$ are the DOFs for $g(x)$.
These coefficients are determined by the following system describing the boundary conditions,
\begin{equation}
\begin{pmatrix}
0 & 1 \\
1 & 1
\end{pmatrix}
\begin{pmatrix}
c_0 \\ c_1
\end{pmatrix}
=
\begin{pmatrix}
\partial_xu_0 - \partial_xNN(0;\theta)\\
u_1 - NN(1;\theta)
\end{pmatrix}.
\end{equation}
These DOFs are therefore determined by $\theta$. In this way, the solution ansatz (\ref{eq:buffer1}) hard-constrains the boundary condition for the current neural-network parameters $\theta$ at each training iteration.

\subsubsection{One-dimensional problem with interfaces}\label{subsec:buffer-1d-interface}
This approach can be extended to problems with interfaces, such as the target problem in Section~\ref{subsec:1interface}.
For the example of a one-dimensional Poisson problem (Eq. \ref{eq:poisson}) with discontinuous diffusivity
\begin{equation}
\kappa(x) = \kappa_m
\qquad
\forall x\in\Omega_m,
\end{equation}
with the same Neumann-Dirichlet boundary condition (\ref{eq:1d-ND-bc}) and the interface condition (\ref{eq:itf_cond}).
As in the previous section, with the solution ansatz (\ref{eq:buffer}),
we obtain the same boundary condition of $g_1$ and $g_M$ for the left-most and right-most subdomains,
\begin{subequations}\label{eq:g-bc}
\begin{equation}
g'_1(0) = \partial_xu_0 - NN'_1(0;\theta_1)
\end{equation}
\begin{equation}
g_M(1) = u_1 - NN_M(1;\theta_M).
\end{equation}
\end{subequations}
On the other hand, (\ref{eq:itf_cond}a) provides a jump condition for $g_i$ and $g_j$ at the interface $x_{ij}$,
\begin{equation}\label{eq:g-jump}
\jump{NN}_{ij} + \jump{g}_{ij} = 0,
\end{equation}
which provides only one equation for two unknowns, $g_i(x_{ij})$ and $g_j(x_{ij})$, and an additional arbitrary constraint is needed to determine their values.
\refA{
Likewise, the flux condition (\ref{eq:itf_cond}b) provides only one equation
for the two unknowns $\kappa_ig'_i(x_{ij})$ and $\kappa_jg'_j(x_{ij})$,
\begin{equation}\label{eq:g-flux}
\avg{\bn\cdot\kappa\nabla NN}_{ij} + \avg{\bn\cdot\kappa\nabla g}_{ij} = 0.
\end{equation}
Each of (\ref{eq:g-jump}) and (\ref{eq:g-flux}) must therefore be supplemented by one
additional condition, which can be arbitrary as long as it is linearly independent of the
condition it accompanies.
Here we prescribe a fixed ratio between the two unknowns of each pair,
\begin{subequations}\label{eq:g-itf-ratio}
\begin{equation}
g_i(x_{ij}) = - \gamma_0 g_j(x_{ij})
\end{equation}
\begin{equation}
n_i\kappa_ig'_i(x_{ij}) = \gamma_1 n_j\kappa_jg'_j(x_{ij}).
\end{equation}
\end{subequations}
\refB{Note that the two
buffer values carry opposite signs while the two normal fluxes carry the same sign,
since the latter are already expressed in each subdomain's own outward normal.}
While $\gamma_0=\gamma_1=1$ is the intuitive choice, the two ratios are otherwise
unconstrained and are left as hyperparameters to be tuned.
The sensitivity of the trained solution to them is examined in
Section~\ref{subsubsec:buffer2d-sensitivity}.
Solving (\ref{eq:g-jump}) and (\ref{eq:g-flux}) together with (\ref{eq:g-itf-ratio})
yields the interface conditions for the buffer functions,
\begin{subequations}\label{eq:g-itf}
\begin{equation}
\begin{split}
g_i(x_{ij}) &= -\frac{\gamma_0}{1 + \gamma_0}\jump{NN}_{ij}\\
g_j(x_{ij}) &= \frac{1}{1 + \gamma_0}\jump{NN}_{ij},
\end{split}
\end{equation}
with $\jump{NN}_{ij}\equiv NN_i(x_{ij};\theta_i) - NN_j(x_{ij};\theta_j)$, and
\begin{equation}
\begin{split}
n_i\kappa_i g'_i(x_{ij}) &= -\frac{2\gamma_1}{1 + \gamma_1}\avg{\bn\cdot\kappa\nabla NN}_{ij}\\
n_j\kappa_j g'_j(x_{ij}) &= -\frac{2}{1 + \gamma_1}\avg{\bn\cdot\kappa\nabla NN}_{ij},
\end{split}
\end{equation}
\refB{with $\avg{\bn\cdot\kappa\nabla NN}_{ij}\equiv(n_i\kappa_iNN'_i(x_{ij};\theta_i) + n_j\kappa_jNN'_j(x_{ij};\theta_j))/2$.}
\end{subequations}
When $\gamma_0=1$ or $\gamma_1=1$, the corresponding pair of unknowns becomes equal in
magnitude, recovering the even split between the two subdomains.
}
Equations (\ref{eq:g-bc}) and (\ref{eq:g-itf}) together specify the boundary and interface conditions that must be satisfied by each buffer function $g_m(x)$.
Since the leftmost and rightmost subdomains each share one external boundary and one internal interface, their corresponding buffer functions $g_1$ and $g_M$ are subject to three conditions: one boundary condition and two interface conditions. In contrast, each interior subdomain is bounded by two interfaces and therefore requires four interface conditions on $g_m$ with $m=2,\ldots, M-1$.
We therefore choose quadratic polynomials for $g_1$ and $g_M$, and cubic polynomials for the interior buffer functions.
\begin{equation}\label{eq:g-poly}
g_m(x) =
\begin{cases}
c_{m,0} + c_{m,1}x + c_{m,2}x^2 & m = 1, M \\
c_{m,0} + c_{m,1}x + c_{m,2}x^2 + c_{m,3}x^3 & m = 2, \ldots, M-1.
\end{cases}
\end{equation}
Substituting (\ref{eq:g-poly}) into (\ref{eq:g-bc}) and (\ref{eq:g-itf})
leads to the systems of equations to determine DOFs of $g_m$,
\begin{subequations}\label{eq:gcoeffs}
\begin{equation}
\begin{pmatrix}
0 & 1 & 0 \\
1 & x_{12} & x_{12}^2 \\
0 & 1 & 2x_{12}
\end{pmatrix}
\begin{pmatrix}
c_{1,0} \\ c_{1,1} \\ c_{1,2}
\end{pmatrix}
=
\begin{pmatrix}
\partial_xu_0 - \partial_xNN_1(0;\theta_1)\\
-\frac{1}{2}\jump{NN}_{12}\\
-\frac{1}{2}\avg{\bn\cdot\kappa\nabla NN}_{12}\\
\end{pmatrix}
\end{equation}
\begin{equation}
\begin{pmatrix}
1 & 1 & 1 \\
1 & x_{M-1,M} & x_{M-1,M}^2 \\
0 & 1 & 2x_{M-1,M}
\end{pmatrix}
\begin{pmatrix}
c_{M,0} \\ c_{M,1} \\ c_{M,2}
\end{pmatrix}
=
\begin{pmatrix}
u_1 - NN_M(1;\theta_M)\\
\frac{1}{2}\jump{NN}_{M-1,M}\\
\frac{1}{2}\avg{\bn\cdot\kappa\nabla NN}_{M-1,M}\\
\end{pmatrix},
\end{equation}
and for $m=2,\ldots,M-1$,
\begin{equation}
\begin{pmatrix}
1 & x_{m-1,m} & x_{m-1,m}^2 & x_{m-1,m}^3 \\
0 & 1 & 2x_{m-1,m} & 3x_{m-1,m}^2 \\
1 & x_{m,m+1} & x_{m,m+1}^2 & x_{m,m+1}^3 \\
0 & 1 & 2x_{m,m+1} & 3x_{m,m+1}^2
\end{pmatrix}
\begin{pmatrix}
c_{m,0} \\ c_{m,1} \\ c_{m,2} \\ c_{m,3}
\end{pmatrix}
=
\begin{pmatrix}
\frac{1}{2}\jump{NN}_{m-1,m}\\
\frac{1}{2}\avg{\bn\cdot\kappa\nabla NN}_{m-1,m}\\
-\frac{1}{2}\jump{NN}_{m,m+1}\\
-\frac{1}{2}\avg{\bn\cdot\kappa\nabla NN}_{m,m+1}\\
\end{pmatrix}.
\end{equation}
\end{subequations}
\par
Now, we introduce the procedure to evaluate the solution ansatz (\ref{eq:buffer}) at a given set of parameters $\{\theta_m\}$
in Algorithm~\ref{alg:buffer}. Again, the degrees of freedom of $g_m$ depend on neural-network parameters $\theta$, so the solution ansatz (\ref{eq:buffer}) hard-constrains the boundary and interface conditions for the current value of $\theta$ at each training iteration.
\begin{algorithm}[tbh]
\caption{Function evaluation for buffer approach (\ref{eq:buffer})}
\begin{algorithmic}\label{alg:buffer}
\item[] \textbf{Input:} Collocation points $\{\bx_k\}$ and NN parameters $\{\theta_m\}$
\item[] \textbf{Output:} Solution values $\{u_{NN}(\bx_k)\}$
\STATE Evaluate boundary/interface NN terms in (\ref{eq:gcoeffs})
\STATE Solve (\ref{eq:gcoeffs}) for $g_m(\bx)$ coefficients
\STATE Evaluate (\ref{eq:buffer}) on collocation points $\{\bx_k\}$
\RETURN $\{u_{NN}(\bx_k)\}$
\end{algorithmic}
\end{algorithm}
\refB{
The system matrix in (\ref{eq:gcoeffs}) depends only on the interface locations and not on
$\theta$. The same holds for the higher-dimensional systems introduced in
Section~\ref{subsec:buffer-multi-dim}, whose matrices are fixed once the geometry, sample
points, and radii are chosen. Only the right-hand side changes during training. Each system
is therefore LU-factorized once before training and reused at every iteration, so an
iteration costs only a triangular solve. Because the matrix is constant, differentiating
through the buffer correction requires no derivative of the factorization: the gradient
propagates through the right-hand side alone, at the cost of a second triangular solve with
the same factors. With $n_m$ buffer degrees of freedom in subdomain $m$, the per-iteration
cost is thus $\cO(n_m^2)$, rather than the $\cO(n_m^3)$ that refactorization would incur.
}
\par
\refA{
The system \eqref{eq:gcoeffs} is non-singular for any strictly ordered partition
$0<x_{12}<\cdots<x_{M-1,M}<1$, so the buffer coefficients are uniquely determined
irrespective of where the interfaces are placed.
It degenerates only when an interface collides with a neighboring interface or with a
domain boundary, in which case the conditioning deteriorates with the vanishing subdomain width.
Such a limit reflects a degenerate partition of the posed physics itself rather than a
deficiency of the buffer construction, and would call for a rescaling of the constraints regardless of numerics, presumably through an asymptotic treatment~\cite{hashin2002thin,benveniste2006general}.
We also emphasize that the buffer function admits considerable freedom in its choice of
ansatz, and that polynomials are by no means the only feasible form. In fact, a radial basis
construction is introduced subsequently for higher dimensions.
}

\subsubsection{Extension to higher dimensions}\label{subsec:buffer-multi-dim}
We now  extend the buffer construction described above to higher dimensional settings.
For higher-dimensional spaces, the natural extensions of boundary conditions (\ref{eq:g-bc}) and (\ref{eq:g-itf}) are
\begin{subequations}\label{eq:g-bc-high-dim}
\begin{equation}
g_b(\bx)
=
u(\bx)
- NN_b(\bx)
\qquad
\forall \bx\in\Gamma_b,
\end{equation}
for the Dirichlet condition,
\begin{equation}
\bn\cdot\nabla g_b(\bx)
=
\bn\cdot\nabla u(\bx)
- \bn\cdot\nabla NN_b(\bx)
\qquad
\forall \bx\in\Gamma_b,
\end{equation}
for the Neumann condition, and
\refA{
\begin{equation}
\begin{split}
g_i(\bx) &= -\frac{\gamma_0}{1 + \gamma_0}\jump{NN}_{ij}\\
g_j(\bx) &= \frac{1}{1 + \gamma_0}\jump{NN}_{ij},
\end{split}
\end{equation}
and
\begin{equation}
\begin{split}
\bn_i\cdot\kappa_i\nabla g_i(\bx) &= -\frac{2\gamma_1}{1 + \gamma_1}\avg{\bn\cdot\kappa \nabla NN}_{ij}\\
\bn_j\cdot\kappa_j\nabla g_j(\bx) &= -\frac{2}{1 + \gamma_1}\avg{\bn\cdot\kappa \nabla NN}_{ij},
\end{split}
\qquad \forall\bx\in\Gamma_{ij},
\end{equation}
}
\end{subequations}
for the interface condition.
\par
Since the buffer function $g_m(\bx)$ is parameterized with finite degrees of freedom,
we enforce condition \eqref{eq:g-bc-high-dim} through a carefully designed boundary sampling strategy.
Specifically, we select a finite subset of boundary points $\partial\tOmega_m \subset \partial\Omega_m$ and determine the DOFs of $g_m(\bx)$ so that condition \eqref{eq:g-bc-high-dim} is satisfied exactly at these representative locations.
$g_m(\bx)$ at boundary points outside $\partial\tOmega_m$ are smoothly interpolated between the sample points in $\partial\tOmega_m$,
thus satisfying (\ref{eq:g-bc-high-dim}) only approximately.
In this regard, the buffer approach may be considered an ``approximate'' hard-constraining or ``discrete'' hard-constraining strategy.
\par
We suppose a subdomain boundary $\partial\tOmega_m$
has $D_m$ Dirichlet boundary sample points, $N_m$ Neumann boundary sample points, and $I_m$ interface sample points.
Then $g_m(\bx)$ requires $(D_m + N_m + 2I_m)$ DOFs to exactly satisfy boundary condition (\ref{eq:g-bc-high-dim}) on these sample points.
A choice of $g_m(\bx)$ would be a sum of radial basis functions (RBFs).
For the $d$-th Dirichlet boundary sample point at $\bx_{m, d}$,
\begin{subequations}\label{eq:g-rbf}
\begin{equation}
g_{m, d}^{(D)}(\bx)
= c^{(D)}_{m,d}\exp\left(-\frac{r_{m,d}^2}{r_0^2}\right),
\end{equation}
with radius $r_{m,k} = \Vert \bx - \bx_{m,k}\Vert$ and DOF $c^{(D)}_{m,d}$.
The base radius $r_0$, as a hyper-parameter, determines the overall smoothness of $g_m(\bx)$.
For the $n$-th Neumann boundary sample point at $\bx_{m, n}$,
\begin{equation}
g_{m, n}^{(N)}(\bx)
= c^{(N)}_{m,n}\bn_{m,n}\cdot(\bx - \bx_{m, n})\exp\left(-\frac{r_{m,n}^2}{r_0^2}\right),
\end{equation}
with DOF $c^{(N)}_{m,n}$ and the normal vector $\bn_{m,n}$ located at $\bx_{m,n}$.
For the $i$-th interface sample point at $\bx_{m, i}$,
\begin{equation}
g_{m, i}^{(I)}(\bx)
= c^{(I0)}_{m,i}\exp\left(-\frac{r_{m,i}^2}{r_0^2}\right)
+
c^{(I1)}_{m,i}\bn_{m,i}\cdot(\bx - \bx_{m, i})\exp\left(-\frac{r_{m,i}^2}{r_0^2}\right),
\end{equation}
with DOFs $c^{(I0)}_{m,i}$ and $c^{(I1)}_{m,i}$.
Overall, the buffer function at the subdomain $\Omega_m$ is defined as
\begin{equation}
g_m(\bx) =
\sum_d^{D_m} g_{m, d}^{(D)}(\bx)
+
\sum_n^{N_m} g_{m, n}^{(N)}(\bx)
+
\sum_i^{I_m} g_{m, i}^{(I)}(\bx),
\end{equation}
\end{subequations}
with DOFs $\{c^{(D)}_{m,d}\}_{d=1}^{D_m} + \{c^{(N)}_{m,n}\}_{n=1}^{N_m} + \{c^{(0)}_{m,i}, c^{(I1)}_{m,i}\}_{i=1}^{I_m}$.
Substituting (\ref{eq:g-rbf}) into (\ref{eq:g-bc-high-dim}) at the sample points $\partial\tOmega_m$
yields the system of equations for $g_m$ DOFs, analogous to (\ref{eq:gcoeffs}).
\par
The number of sample points on a boundary surface, their respective radii,
and the choice of interpolation basis functions collectively influence
the convergence of training.
A thorough investigation of these effects and the prescription of optimal values
are extensive and lie beyond the scope of this study.
In practice,
the convergence of the buffer method \eqref{eq:buffer}
depends on the condition number of the linear systems that determine the buffer
coefficients (e.g., Eq.~\ref{eq:gcoeffs}).
For the two-dimensional cases considered here,
the number of sample points is chosen so that the condition number of each
subdomain system remains moderate ($\lesssim\cO(10^2)$)~\cite{schaback1995error}.
The sample points are placed at the roots of Legendre polynomials along the boundary
or interface, in order to suppress interpolation oscillations near the
endpoints~\cite{berrut2004barycentric}.
The radius of each RBF is then determined by the number of sample points $N_s$
and the boundary length $L_s$,
\begin{equation}
r_0 = \frac{L_s}{N_s+1},
\end{equation}
and is shared by all sample points on the same boundary.
This illustrates that the hyperparameters of the buffer functions can be selected
systematically on the basis of desired accuracy and conditioning. Although soft-constrained PINNs need not rely solely on trial-and-error penalty selection, alternatives such as adaptive weighting, gradient-balancing, and augmented-Lagrangian strategies typically introduce auxiliary update rules and associated hyperparameters ~\cite{wang2021understanding, xiang2022self, son2023enhanced}. Consequently, their performance can depend on the design and tuning of these mechanisms while remaining subject to the trade-offs inherent in multi-objective optimization, where constraint satisfaction must be balanced against PDE residual minimization. In contrast, the buffer formulation embeds the constraint correction directly into the trial representation.
\par
Recently, similar additive approaches have been proposed using transfinite
formulations based on generalized barycentric coordinates and TFC-based constrained
expressions on mapped domains~\cite{dong2026novel,sukumar2026wachspress}.
While these formulations construct trial functions that satisfy boundary conditions exactly,
they rely on geometry-specific representations such as coordinate mappings
or specialized basis functions, which are presently restricted to two-dimensional domains;
their extension to interface conditions is also not yet straightforward.
\par
In contrast, the buffer approach enforces boundary and interface conditions through
discrete, local corrections determined from the neural network mismatch at sampled
points, without requiring an explicit global lifting operator.
\refA{The system that determines the buffer coefficients is assembled simply by
substituting the buffer ansatz into the prescribed condition at each boundary or interface
sample point. It can therefore be constructed for virtually any form of boundary or
interface condition, on any geometry for which such sample points can be placed, and the
coefficients follow from a single linear solve at each training iteration.
The radial basis form (\ref{eq:g-rbf}) is adopted here as a convenient realization of this
idea, since its basis functions are centered at the sample points and require no coordinate
mapping or specialized basis construction; the formulation itself does not depend on this
particular choice.}
This makes the method straightforward to implement and applicable to general
high-dimensional geometries with interfaces, at the cost of enforcing the
constraints exactly only at discrete locations.



\section{Results and Discussion}\label{sec:examples}
In this section, we present the numerical results for the proposed formulations on a sequence of benchmark interface problems of increasing complexity. Our goals are to assess their accuracy relative to soft-constrained PINNs, to clarify the distinct behaviors of the windowing and buffer approaches, and to examine their robustness as the source terms and geometries become more challenging. The one-dimensional examples highlight the approximation and optimization mechanisms of the two methods, while the two-dimensional example tests their practical performance in a more complex setting.

\subsection{Model Problems}
We begin by describing the one- and two-dimensional benchmark problems used to evaluate the proposed formulations.


\subsection*{Problem 1: Uniform source with single interface}\label{subsec:1interface}
First, we solve the one-dimensional diffusion equation on a composite domain consisting of two materials characterized by different diffusivities. We consider the one-dimensional domain $\Omega=[0,1]$
composed of two subdomains with their interface at $\Gamma_{12} = \{x_{itf}\equiv 0.5\}$.
The diffusivity $\kappa$ is piecewise constant,
\begin{equation}
\kappa =
\begin{cases}
\kappa_1 \equiv 0.1 & x < x_{itf} \\
\kappa_2 \equiv 1 & x \ge x_{itf},
\end{cases}
\end{equation}
and the source term $f(\bx)=1$ is constant over the domain.
We consider homogeneous Dirichlet boundary conditions at the left and right boundaries,
\begin{equation}\label{eq:homo-diri-bc}
u(0) = u(1) = 0.
\end{equation}
This problem has an analytic solution of the form,
\begin{equation}
u(\bx) =
\begin{cases}
-\frac{x^2}{2\kappa_1} + c_1\frac{x}{\kappa_1} + c_2 & x < x_{itf} \\    
-\frac{x^2}{2\kappa_2} + c_3\frac{x}{\kappa_2} + c_4 & x \ge x_{itf}.
\end{cases}
\end{equation}
The coefficients $c_1$, $c_2$, $c_3$ and $c_4$ are determined
by the boundary and interface conditions,
which constitutes the following system,
\begin{equation}
\begin{pmatrix}
0 & 1 & 0 & 0 \\
0 & 0 & \frac{1}{\kappa_2} & 1 \\
\frac{x_{itf}}{\kappa_1} & 1 & -\frac{x_{itf}}{\kappa_2} & -1 \\
x_{itf} & 0 & -x_{itf} & 0
\end{pmatrix}
\begin{pmatrix}
c_1 \\ c_2 \\ c_3 \\ c_4
\end{pmatrix}
=
\begin{pmatrix}
0 \\ \frac{1}{2\kappa_2} \\ \frac{x_{itf}^2}{2\kappa_1} - \frac{x_{itf}^2}{2\kappa_2} \\ 0
\end{pmatrix}.
\end{equation}
This problem serves as a baseline one-dimensional interface test, designed to assess the accuracy of the proposed methods in the simplest setting with a single diffusivity jump.

\subsection*{Problem 2: Uniform source with three interfaces}\label{subsec:3interface}
We solve the one-dimensional diffusion equation on a domain with three interfaces. This problem is used to evaluate how the methods perform as the number of interfaces increases, without introducing additional complexity in the source term or boundary conditions. The interfaces are located at $\Gamma=\{0.25, 0.5, 0.75\}$~\cite{sarma2024interface},
with the same source term $f(x)=1$ and the same boundary condition (\ref{eq:homo-diri-bc}) as the first problem.
The analytical solution can be found in~\ref{app:prob2-sol}~and~\cite{sarma2024interface}. 

\subsection*{Problem 3: Spatially varying source with one interface}\label{subsec:1int-spatially-varying}
We next consider the one-dimensional diffusion equation with a spatially varying source term across the interface, which provides a more demanding test of residual approximation and training behavior. The interface is located at $\Gamma=\{x_{itf}=0.5\}$, and the source term is defined as,
\begin{equation}
-\nabla\cdot(\kappa\nabla u)
=
\begin{cases}
f_0 & x < x_{itf}\\
f(x) \equiv a\exp\left(-\frac{(x-x_c)^2}{w^2}\right) & x > x_{itf},
\end{cases}
\end{equation}
with $f_0=-0.05$, $a=1$, $x_c=0.75$, and $w=0.1$.
For the boundary condition, we consider a homogeneous Neumann condition at $x=0$ and a homogeneous Dirichlet condition at $x=1$,
\begin{subequations}
\begin{equation}
\Dpartial{u}{x}(0) = 0
\end{equation}
\begin{equation}
u(1) = 0.
\end{equation}
\end{subequations}
This problem has an analytic solution of the form,
\begin{equation}
u(\bx) =
\begin{cases}
-f_0\frac{x^2}{2\kappa_1} + c_1\frac{x}{\kappa_1} + c_2 & x < x_{itf} \\
-F(x) + c_3\frac{x}{\kappa_2} + c_4 & x \ge x_{itf},
\end{cases}
\end{equation}
where $F(x)$ is the double integral of $f(x)$,
\begin{equation}
F(x) = \frac{aw\sqrt{\pi}}{2\kappa_2}(x - x_c)\,\mathrm{erf}\!\left(\frac{x-x_c}{w}\right)
+ \frac{aw^2}{2\kappa_2}\exp\left(-\frac{(x-x_c)^2}{w^2}\right).
\end{equation}
The coefficients are determined by the following system that describes the boundary and interface conditions,
\begin{equation}
\begin{pmatrix}
1 & 0 & 0 & 0 \\
0 & 0 & \frac{1}{\kappa_2} & 1 \\
\frac{x_{itf}}{\kappa_1} & 1 & -\frac{x_{itf}}{\kappa_2} & -1 \\
1 & 0 & -1 & 0
\end{pmatrix}
\begin{pmatrix}
c_1 \\ c_2 \\ c_3 \\ c_4
\end{pmatrix}
=
\begin{pmatrix}
0 \\ F(1) \\ \frac{f_0x_{itf}^2}{2\kappa_1} - F(x_{itf}) \\ f_0x_{itf} - F_x(x_{itf})
\end{pmatrix},
\end{equation}
where $F_x(x)$ is the integral of $f(x)$,
\begin{equation}
F_x(x) = \frac{aw\sqrt{\pi}}{2\kappa_2}\,\mathrm{erf}\!\left(\frac{x-x_c}{w}\right).
\end{equation}

\subsection*{Problem 4: Slanted interface in a two-dimensional domain with multiple source terms}\label{subsec:2dpoisson}
Finally, we consider a two-dimensional rectangular domain with a slanted interface. This problem extends the study to two dimensions, where the slanted interface and mixed boundary conditions provide a more challenging test of robustness and geometric generalization. The domain is defined on $\bx = (x, y) \in \Omega=[0, 2]\times[0, 1]$. The interface is defined as $\Gamma = \{(x, y) | y = (x - x_b) / (x_t - x_b) \}$ with $x_b=0.8$ and $x_t=1.2$.
The diffusivity is defined as
\begin{equation}
\kappa
=
\begin{cases}
0.1 & y >\frac{x - x_b}{x_t - x_b} \\
1 & \text{otherwise}.
\end{cases}
\end{equation}
A homogeneous Dirichlet condition is assigned on the right-side subdomain boundary,
\begin{equation}
u(\bx)=0 \qquad\text{on }\partial\Omega_R,
\end{equation}
with
$\partial\Omega_R=\{(x,y)| y < (x - x_b) / (x_t - x_b)\} \cap \{(x, y)|y=0 \text{ or } 1\text{, or } x=2\}$.
The $\partial\Omega_R$ includes a subset of the outer boundary that is restricted to the right-side subdomain defined by $y < (x-x_b)/(x_t-x_b)$.
On the left-side subdomain, a
homogeneous Neumann condition is assigned,
\begin{equation}
\bn\cdot\kappa\nabla u(\bx)=0 \qquad\text{on }\partial\Omega_L,
\end{equation}
with
$\partial\Omega_L=\{(x,y)| y > (x - x_b) / (x_t - x_b)\} \cap \{(x, y)|y=0 \text{ or } 1\text{, or } x=0\}$. The $\partial\Omega_L$ includes a subset of the outer boundary that is restricted to the left-side subdomain defined by $y > (x-x_b)/(x_t-x_b)$.
The source term is defined as the sum of
three Gaussian sources,
\begin{equation}\label{eq:2dpoisson-f}
f(\bx)
=
\sum_{i=1}^3
A_i\exp\left( - \frac{\Vert \bx - \bx_i \Vert^2}{r_i^2} \right),
\end{equation}
with $\{\bx_i\} = \{(0.3, 0.6), (1, 0.2), (1.6, 0.7)\}$,
$\{r_i\} = \{0.08, 0.2, 0.1\}$
and
$\{A_i\} = \{10, 20, 15\}$.
\begin{figure}[tbph]
\hspace*{-1.8cm}
    \begin{tikzpicture}[
]
    \begin{groupplot}[
        group style={
            group name = my plots,
            group size= 2 by 1,
            xlabels at =edge bottom,
            horizontal sep=3cm,
            vertical sep=2.2cm,
        },
        height = 0.3\textwidth,
        width = 0.5\textwidth,
        name=chung,
    ]    
\pgfplotsset{set layers=standard}%

        \nextgroupplot[
            ylabel={$y$},
            xlabel={$x$},
            tick scale binop ={\times},
            xmin = 0, xmax = 2,
            ymin = 0, ymax = 1,
            colorbar,
            colormap name=viridis,
            point meta min=0,   
            point meta max=2.66,   
        ]
        
            \edef\imagepath{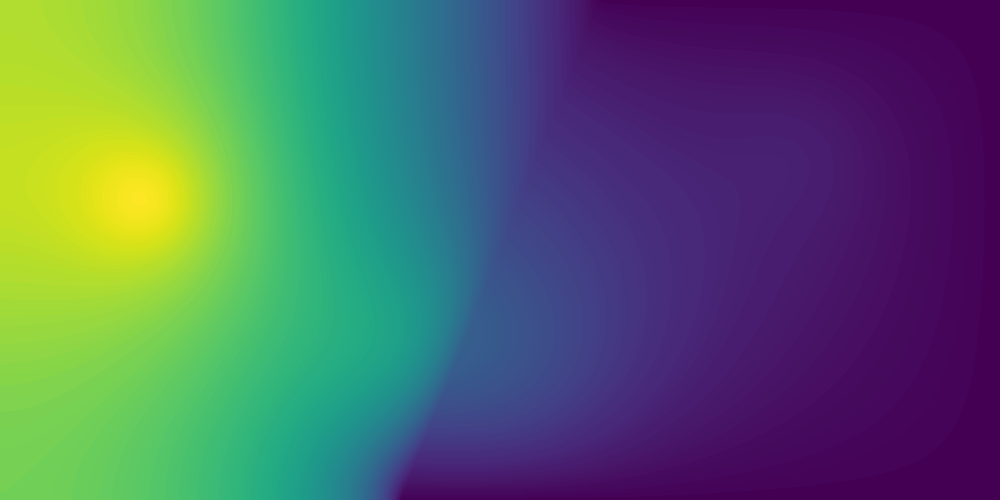}
            \addplot graphics[xmin=0, xmax=2, ymin=0, ymax=1]{\imagepath};

            \addplot[
                domain=0.8:1.2, samples=100, dashed, red,
            ]
            {(x-0.8)/0.4};

        \nextgroupplot[
            ylabel={$y$},
            xlabel={$x$},
            tick scale binop ={\times},
            xmin = 0, xmax = 2,
            ymin = 0, ymax = 1,
            colorbar,
            colormap name=viridis,
            point meta min=0,   
            point meta max=20
            ,   
        ]
        
            \edef\imagepath{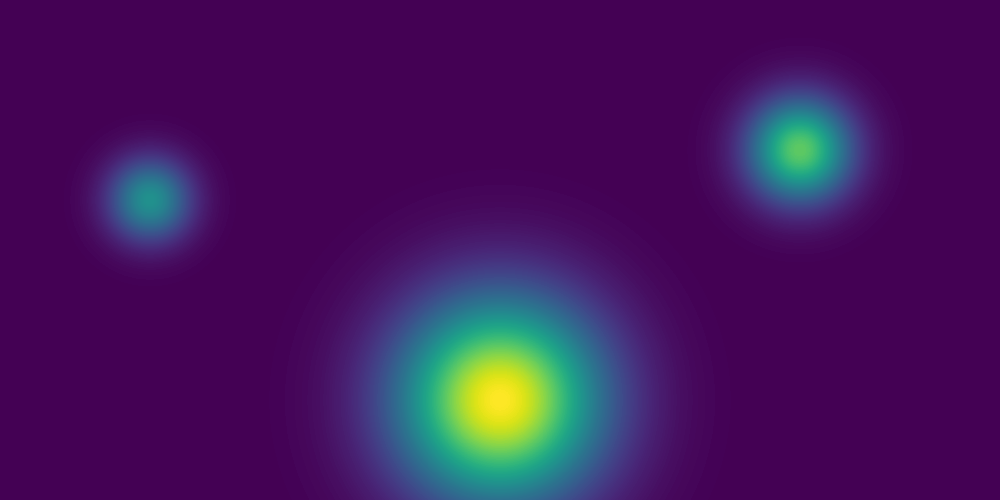}
            \addplot graphics[xmin=0, xmax=2, ymin=0, ymax=1]{\imagepath};

  \end{groupplot}
\node[below = 2cm of my plots c1r1.south west,
    anchor=west,
] {(a) $u(\bx)$};
\node[below = 2cm of my plots c2r1.south west,
    anchor=west,
] {(b) $f(\bx)$};
\end{tikzpicture}
%
    \caption{Model problem with a slanted interface in a two-dimensional domain: (a) solution; and (b) source term.}
    \label{fig:2dpoisson_truth}
\end{figure}
Figure~\ref{fig:2dpoisson_truth} shows the solution and source term of the model problem.
The solution is obtained via a finite-element solver~\cite{mfem}
with a fine mesh of $\sim132,000$ grid points.
\subsection{Comparison between the windowing and buffer approaches}\label{subsec:window-vs-buffer}

Having established the two hard-constrained formulations, we now compare the training performance of the windowing and buffer approaches on representative interface problems. This comparison highlights the distinct structural properties of the two ansatz constructions and clarifies their respective strengths in practical settings.
\par
For the windowing approach, we first conducted a preliminary study to identify an effective window configuration in terms of polynomial order and boundary/interface subdomain overlap. Based on these results, we adopt the window functions and overlapping configuration that demonstrated the most accurate and stable training behavior. Details of this study are provided in~\ref{app:window-prelim}. In particular, we employ $\tW^{(1)}_d$ and $\tW^{(1)}_n$ for the boundary and interface subdomains (see Section~\ref{sec:window_functions}), together with a uniform interior subdomain size $\Delta x=0.25$. The boundary and interface subdomains are assigned size $2\Delta x$, resulting in perfect overlap between neighboring regions.
\par
The PINN models for the windowing approach \eqref{eq:window-ansatz} and the buffer approach \eqref{eq:buffer}
are trained for the target problem with a spatially-varying source term in Section~\ref{subsec:1int-spatially-varying}.
For the windowing approach,
the interior subdomain centers are located at $x_1=0.25$ and $x_2=0.75$.
For the buffer approach, the configuration described in Section~\ref{subsec:buffer-1d-interface} is used.
For both approaches, a fully connected network with hidden layers $[12, 12]$ is
used for each interior subdomain.
For all cases, the training loss \eqref{eq:J-int} is evaluated
with $40$ collocation points uniformly distributed in the domain.
For both approaches, the training is performed via SOAP optimizer\footnote{SOAP
(ShampoO with Adam in the Preconditioner's eigenbasis) is a second-order optimizer designed to bridge the gap between the high performance of Shampoo and the computational efficiency of AdamW.}~\cite{vyas2024soap},
with learning rate $5\times10^{-3}$ for $3\times10^4$ iterations.
\par
\begin{figure}[tbhp]
    \input{figures_config3_rhs.tex}
    \caption{The training performance of the windowing approach on the problem in Section~\ref{subsec:1int-spatially-varying}:
    the left-hand side of the solution compared with the source term,
    using the interior window function (a) $\tW^{(1)}_{int}$, (b) $\tW^{(2)}_{int}$, or (c) $\tW^{(3)}_{int}$; and
    (d) the left-hand side contribution of each subdomain solution for the best case of $\tW^{(1)}_{int}$.}
    \label{fig:config3-rhs}
\end{figure}
The training results are shown in Figure~\ref{fig:config3-rhs}
and Figure~\ref{fig:config3-buffer}.
While the windowing approach managed to achieve the solution relative $L_2$ error of $5.71\times10^{-4}$ in the best case,
the buffer approach achieved a smaller error of $7.36\times10^{-5}$.
To explain the performance gap, we assess the impact of the shape of window functions.
Figure~\ref{fig:config3-rhs}~(a--c)
show $\LHS[u_{NN}]$ with different interior window functions $\tW^{(k)}_{int}$, comparing with the analytic source term $f(x)$.
While LHS matches well with the source term overall,
a major discrepancy appears at $x=0.75$.
In particular, 
we recognize that the shapes of $\LHS[u_{NN}]$ are strongly influenced by the second derivatives of $\tW_{int}$.
The high-order polynomial $\tW^{(3)}_{int}$, with its smooth second derivatives,
achieves relative $L_2$ error of $1.19\%$ for $\LHS[u]$, smaller than $\tW^{(1)}_{int}$ or $\tW^{(2)}_{int}$.
However, the resulting solution $u$ turns out to be less accurate than $\tW^{(1)}_{int}$,
implying that its $\LHS[u]$ is subject to constant or low-frequency bias.
We also emphasize that the LHS contribution from the boundary and interface subdomains is minimal;
as shown in Figure~\ref{fig:config3-rhs}~(d),
LHS of the boundary/interface subdomains are summed up to be a constant function,
and the interior subdomain solutions mainly approximate the actual source term.
\par
These results strongly indicate that
the training performance of the windowing approach
depends on how well the derivatives of a given window function can span a source term (or other LHS terms).
Furthermore, the smoothness of the derivatives does not significantly impact the solution accuracy as long as they do not introduce unwanted discontinuities.
These behaviors stem from the fact that LHS of the windowing approach
is mainly expressed as a basis expansion, for example,
\begin{equation}
\LHS[W(x)NN(x)]
= \refB{-\kappa\Dparttwo{W}{x}NN - 2\kappa\Dpartial{W}{x}\Dpartial{NN}{x} - \kappa W\Dparttwo{NN}{x}},
\end{equation}
with window derivatives as basis functions.
\todo{Kevin: am I speaking too strong about this sentence? how do we show windowing approach is also subject to universal approximation theorem?}
However, if the window derivatives cannot span sufficiently general functions,
the expressivity of the associated NN can be restricted by the shape of these derivatives,
which in turn requires a larger NN architecture or longer training iterations.
\par
\begin{figure}[H]
    \begin{tikzpicture}[
]
    \begin{groupplot}[
        group style={
            group name = my plots,
            group size= 2 by 1,
            xlabels at =edge bottom,
            horizontal sep=2cm,
            vertical sep=2.2cm,
        },
        height = 0.45\textwidth,
        width = 0.5\textwidth,
        enlarge x limits={false, abs value = 5mm},
        enlarge y limits={false, abs value = 5mm},
        name=chung,
    ]    
\pgfplotsset{set layers=standard}%

        \nextgroupplot[
            ymin=-0.1, ymax=1.1,
            xlabel={$x$},
            ylabel={$\LHS[u]$},
            tick scale binop ={\times},
            xtick={0, 0.25, 0.5, 0.75, 1.},
        ]

            \addplot+[
                mark=none,
                gray!50,
                line width=1.5,
            ] coordinates {(0.5, -3) (0.5, 0.75)};

            \addplot+ [
                line width=1.0,
                solid,
                smooth,
                mark=none,
                black,
            ]
            table [
                x index=0, y index=1,
            ]{data_config3_buffer_rhs.trained.total.txt};
            \addplot+ [
                line width=1.0,
                dashed,
                mark=none,
                red,
            ]
            table [
                x index=0, y index=2,
            ]{data_config3_buffer_rhs.trained.total.txt};

            \node[anchor=west, font=\tiny] at (rel axis cs: 0., 0.95) (a) {Rel. $L_2$ error};
            \node[anchor=west, font=\tiny] at (rel axis cs: 0., 0.88) {$u$: $7.36\times10^{-5}$};
            \node[anchor=west, font=\tiny] at (rel axis cs: 0., 0.8) {$\LHS[u]$: $8.13\times10^{-4}$};

        \nextgroupplot[
            ymin=-0.3, ymax=1.3,
            xlabel={$x$},
            ylabel={$\LHS[u]$},
            tick scale binop ={\times},
            xtick={0, 0.25, 0.5, 0.75, 1.},
            legend style={
                draw=none, fill=white,
                at={(rel axis cs: 0., 1.0)},
                anchor=north west,
                nodes={scale=1.0},
                legend cell align={left},
                legend columns=1,
                font=\tiny,
                /tikz/every even column/.append style={column sep=0.5cm},
            },
        ]

            \addlegendimage{
                line width=1.0,
                solid,
                mark=none,
                blue,
            };
            \addlegendentry{$\LHS[NN_m]$};
            \addlegendimage{
                line width=1.0,
                solid,
                mark=none,
                orange,
            };
            \addlegendentry{$\LHS[g_m]$, buffer};
            \addlegendimage{
                line width=1.0,
                dashed,
                mark=none,
                black,
            };
            \addlegendentry{$\LHS[u_{NN}]$, total};

            \addplot+[
                mark=none,
                gray!50,
                line width=1.5,
            ] coordinates {(0.5, -3) (0.5, 3)};

            \addplot+ [
                line width=1.0,
                solid,
                mark=none,
                blue,
            ]
            table [
                x index=0, y expr=-\thisrowno{1},
            ]{data_config3_buffer_rhs.trained.component.txt};
            \addplot+ [
                line width=1.0,
                solid,
                mark=none,
                orange,
            ]
            table [
                x index=0, y expr=-\thisrowno{2},
            ]{data_config3_buffer_rhs.trained.component.txt};
            \addplot+ [
                line width=1.0,
                dashed,
                mark=none,
                black,
            ]
            table [
                x index=0, y expr=-\thisrowno{3},
            ]{data_config3_buffer_rhs.trained.component.txt};

  \end{groupplot}
\node[below = 1.5cm of my plots c1r1.south west,
    anchor=west,
] {(a) total $\LHS[u_{NN}]$};
\node[below = 1.5cm of my plots c2r1.south west,
    anchor=west,
] {(b) LHS contributions};
\end{tikzpicture}
%
    \caption{The training performance of the buffer approach on the problem in Section~\ref{subsec:1int-spatially-varying}:
    (a) the left-hand side of the solution compared with the source term; and
    (b) the left-hand side contributions of the neural networks and the buffer function.}
    \label{fig:config3-buffer}
\end{figure}
On the other hand, the NN terms in the buffer approach are not restricted by any window function,
which means that they can approximate the source term far better with much fewer training iterations.
Figure~\ref{fig:config3-buffer}~(a)
shows that the LHS of the buffer approach matches well with the source term.
In Figure~\ref{fig:config3-buffer}~(b), $\LHS[g_m]$ remains to be a constant function,
imposing minimal burden on NN terms in resolving the physics equation residual.
This underscores the advantage of the buffer approach for more general, complex problem settings.

\subsection{Comparison with soft-constrained methods in one-dimensional model problems}

\begin{table}[tbh]
    \centering
    \small
    \resizebox{\textwidth}{!}{%
    \begin{tabular}{|l|c|c|c|}
        \hline
        Problem 1 & Layers & AF & Relative $L_2$ error \\
        \hline
        $\phi$-PINNs & $[2, 12, 12, 1]$ & $\tanh$ & (1.89e-5, 1.18e-4, 1.54e0) \\
        I-PINNs & $[1, 12, 12, 1]$ & Sigmoid, $\tanh$ & (4.38e-5, 7.04e-4, 1.62e0) \\
        AdaI-PINNs & $[1, 12, 12, 1]$ & Adaptive $\tanh$ & (1.37e-5, 4.16e-4, \textcolor{red}{1.45e0}) \\
        M-PINN & $[1, 12, 12, 1]\times2$ & $\tanh$ & (8.37e-6, 7.68e-5, 3.58e0) \\
        \hline
        Window & $[1, 12, 12, 1]\times2$ & $\tanh$ & (\textcolor{blue}{2.10e-9}, 7.61e-2, 6.67e0) \\
        Buffer & $[1, 12, 12, 1]\times2$ & $\tanh$ & (1.26e-8, \textcolor{green}{2.61e-5}, 2.26e1) \\
        \hline
    \end{tabular}%
    }
    \caption{\refB{Comparison between different PINN approaches on Problem 1. Each cell reports the $(\min,\mathrm{median},\max)$ of the relative $L_2$ error over $400$ runs ($100$ seeds $\times$ $2$ initializers $\times$ $2$ scales); see \ref{app:init} for the per-initializer breakdown. The best minimum, median, and maximum across the six methods are highlighted in \textcolor{blue}{blue}, \textcolor{green}{green}, and \textcolor{red}{red}, respectively.}}
    \label{tab:init-prob1}
\end{table}
\begin{table}[tbh]
    \centering
    \small
    \resizebox{\textwidth}{!}{%
    \begin{tabular}{|l|c|c|c|}
        \hline
        Problem 2 & Layers & AF & Relative $L_2$ error \\
        \hline
        $\phi$-PINNs & $[2, 12, 12, 1]$ & $\tanh$ & (2.04e-5, 4.05e-4, 1.74e0) \\
        I-PINNs & $[1, 12, 12, 1]$ & Swish-$\tanh$-Sigmoid-Silu & (2.66e-4, 6.84e-2, \textcolor{red}{8.79e-1}) \\
        AdaI-PINNs & $[1, 12, 12, 1]$ & Adaptive $\tanh$ & (3.66e-5, 7.13e-3, 6.86e0) \\
        M-PINN & $[1, 12, 12, 1]\times 4$ & $\tanh$ & (7.99e-6, 1.01e-4, 6.63e0) \\
        \hline
        Window & $[1, 12, 12, 1]\times4$ & $\tanh$ & (\textcolor{blue}{6.06e-8}, 5.65e-4, 3.88e0) \\
        Buffer & $[1, 12, 12, 1]\times4$ & $\tanh$ & (4.03e-6, \textcolor{green}{8.71e-5}, 1.12e2) \\
        \hline
    \end{tabular}%
    }
    \caption{\refB{Comparison between different PINN approaches on Problem 2. Format and highlighting follow Table~\ref{tab:init-prob1}.}}
    \label{tab:init-prob2}
\end{table}
\begin{table}[tbh]
    \centering
    \small
    \resizebox{\textwidth}{!}{%
    \begin{tabular}{|l|c|c|c|}
        \hline
        Problem 3 & Layers & AF & Relative $L_2$ error \\
        \hline
        $\phi$-PINNs & $[1, 12, 12, 1]$ & $\tanh$ & (6.05e-4, 6.64e-2, \textcolor{red}{1.10e0}) \\
        I-PINNs & $[1, 12, 12, 1]$ & Sigmoid, $\tanh$ & (1.03e-3, 1.18e-1, 1.67e2) \\
        AdaI-PINNs & $[1, 12, 12, 1]$ & Adaptive $\tanh$ & (1.36e-2, 1.02e0, 7.64e1) \\
        M-PINN & $[1, 12, 12, 1]\times2$ & $\tanh$ & (7.48e-4, 6.73e-2, 8.91e2) \\
        \hline
        Window & $[1, 12, 12, 1]\times2$ & $\tanh$ & (9.36e-4, 6.42e-2, 1.07e2) \\
        Buffer & $[1, 12, 12, 1]\times2$ & $\tanh$ & (\textcolor{blue}{5.90e-4}, \textcolor{green}{6.21e-2}, 2.04e2) \\
        \hline
    \end{tabular}%
    }
    \caption{\refB{Comparison among different PINN approaches on Problem 3. Format and highlighting follow Table~\ref{tab:init-prob1}.}}
    \label{tab:init-prob3}
\end{table}
Tables~\ref{tab:init-prob1}-\ref{tab:init-prob3}
summarize the comparison among different PINN approaches including soft-constrained methods
on the interface problems presented in Section~\ref{subsec:1interface}. The soft-constrained baselines comprise M-PINN~\cite{zhang2022multi}, I-PINNs~\cite{sarma2024interface}, AdaI-PINNs~\cite{roy2024adaptive}, and $\phi$-PINNs\footnote{$\phi$-PINNs can be regarded as a special case of DCSNN~\cite{hu2022discontinuity}, in which the latent variable is fixed prior to training to one of $\phi=\{-1,0,1\}$.}.
For a consistent comparison, the hidden layer width is fixed at $12$ across all approaches.
The PDE loss \eqref{eq:J-int} is evaluated on
$20$ uniformly distributed collocation points for Problem 1
and $40$ for Problems 2 and 3.
Training is carried out with the SOAP optimizer at a learning rate of $5\times10^{-3}$,
for $10^4$ iterations on Problems 1 and 2 and $3\times10^4$ iterations on Problem 3.
\refB{Since the training results depend strongly on the initial weights of the neural networks and the optimal initialization varies across PINN approaches, we conducted a statistical ablation study across the three one-dimensional problems.
Tables~\ref{tab:init-prob1}--\ref{tab:init-prob3} report the $(\min,\mathrm{median},\max)$ of the relative $L_2$ error aggregated over $400$ runs per method: $100$ independent random seeds under each combination of two weight initializers (\texttt{glorot\_uniform}, \texttt{random\_normal}) and two scales ($1$, $0.1$). Each seed jointly randomizes the network weights and the piecewise diffusivities of the problems. The details of this ablation study is provided in \ref{app:init}.}
\par
\refB{Overall, the buffer approach delivers the smallest median error on all three one-dimensional problems, whereas the windowing approach attains occasional exceptional minimums but does not consistently improve on the soft-constrained baselines in typical performance.
The windowing approach attains the smallest minimum errors on Problems~1 and~2 ($\mathcal{O}(10^{-9})$ and $\mathcal{O}(10^{-8})$, respectively), which are also observed in~\ref{app:window-prelim} for perfectly overlapped subdomains.
However, such low errors are not expected in general cases,
and the median errors of the windowing approach are comparable to or worse than those of the soft-constrained baselines.
The buffer approach, by contrast, achieves the smallest median error over all cases,
regardless of the number of interfaces or the type of source term.
Its minimum errors are also lower than those of the soft-constrained approaches, and are the lowest overall on Problem~3, indicating more consistent behavior across random initializations and diffusivity samples.
Compared to M-PINN, which uses the same network size as the hard-constrained methods, the improvement of the buffer approach is only marginal.
A systematic investigation of buffer function forms that could widen this margin over soft-constrained counterparts is beyond the scope of this study and is left for future work.}

\subsection{Demonstration on a two-dimensional problem}\label{subsec:demo-2dpoisson}

To evaluate the scalability and robustness of the proposed hard-constrained formulations,
we extend our investigation to a two-dimensional elliptic interface problem defined in Section~\ref{subsec:2dpoisson}.
The presence of mixed Dirichlet and Neumann boundary conditions,
together with a nontrivial interface geometry,
makes this case an effective benchmark for assessing the generalization of the windowing and buffer approaches to higher-dimensional settings.

\subsubsection{\refB{Challenges of the windowing approach in higher dimensions}}\label{subsubsec:window2d}
\refB{Unlike the one-dimensional cases, in which the windowing approach showed comparable performance, extending it to higher dimensions is actually challenging.
While the windowing approach benefits from the perfectly overlapping configuration ($\beta=2$, \ref{app:window-prelim}), such a configuration is generally not available in higher dimensions.
In essence, the perfectly overlapping configuration requires an angle-preserving transformation between subdomains and a reference hyper-rectangle, so that the interface normal flux condition is preserved under the transformation.
However, it is not straightforward to find such a mapping for curved interfaces and multiple junctions.
Under any partial-overlap configuration, the window supports cannot reach the opposite sides of their subdomains cleanly---they either extend beyond the subdomain or terminate in its interior---which stiffens the training optimization.}
\par
We first \refB{illustrate this difficulty} with hard constraints applied only at the interface\refB{, using window supports that extend beyond each interior subdomain.}
Training is performed for $3\times10^4$ epochs with the SOAP optimizer
at a learning rate of $10^{-3}$.
The PDE loss \eqref{eq:J-int} is evaluated on $80\times40$ collocation points
uniformly distributed in the interior of the domain.
All Dirichlet and Neumann boundary conditions are enforced softly
through the losses \eqref{eq:J-dbc_1} and \eqref{eq:J-dbc_2}.
The boundary losses are evaluated at $80$ uniformly distributed points along the horizontal
boundaries and $40$ along the vertical boundaries, excluding the corner points.
All boundaries are assigned a unit weight, except for the bottom-left Neumann boundary,
whose weight is relaxed to $10^{-4}$.
The relative $L_2$ error is computed on the $\sim 132{,}000$ finite-element grid points.
The window functions are defined as one-dimensional functions along the normal direction of the interface,
each with size $\Delta x=1.2$, large enough to cover the respective subdomain.
The solution ansatz employs four neural networks:
two 2D networks for the interior subdomains
and two 1D networks for the interface value and slope functions.
\par
\begin{figure}[tbh]
\hspace*{-1.8cm}
    \begin{tikzpicture}[
]
    \begin{groupplot}[
        group style={
            group name = my plots,
            group size= 2 by 2,
            xlabels at =edge bottom,
            horizontal sep=3cm,
            vertical sep=2.2cm,
        },
        height = 0.3\textwidth,
        width = 0.5\textwidth,
        name=chung,
    ]    
\pgfplotsset{set layers=standard}%

        \nextgroupplot[
            ylabel={$y$},
            xlabel={$x$},
            tick scale binop ={\times},
            xmin = 0, xmax = 2,
            ymin = 0, ymax = 1,
            colorbar,
            colormap name=viridis,
            point meta min=-0.01,   
            point meta max=2.48,   
        ]
        
            \edef\imagepath{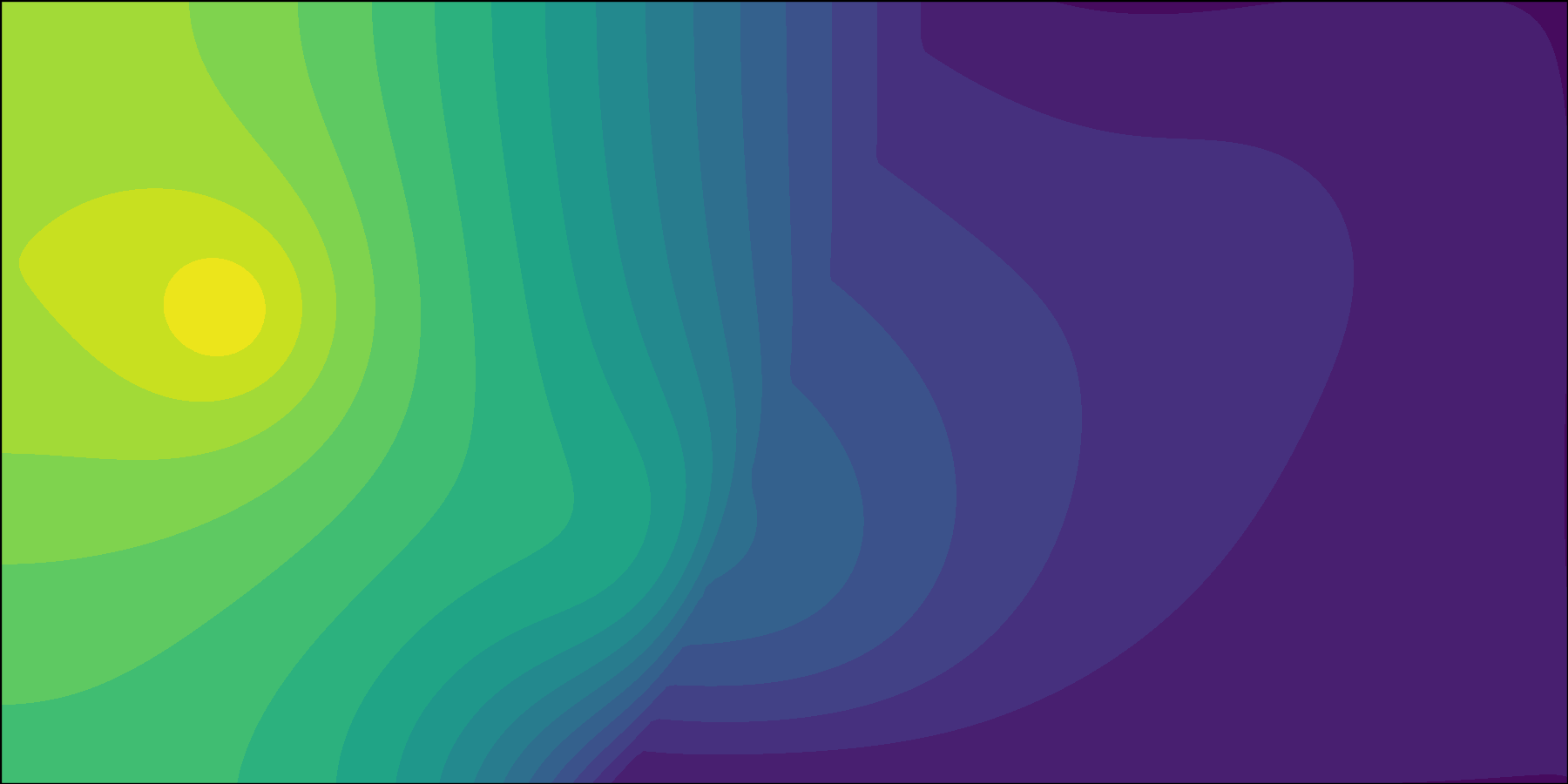}
            \addplot graphics[xmin=0, xmax=2, ymin=0, ymax=1]{\imagepath};

            \addplot[
                domain=0.8:1.2, samples=100, dashed, red,
            ]
            {(x-0.8)/0.4};

        \nextgroupplot[
            ylabel={$y$},
            xlabel={$x$},
            tick scale binop ={\times},
            xmin = 0, xmax = 2,
            ymin = 0, ymax = 1,
            colorbar,
            colormap name=viridis,
            point meta min=-0.22,   
            point meta max=20.0,   
        ]
        
            \edef\imagepath{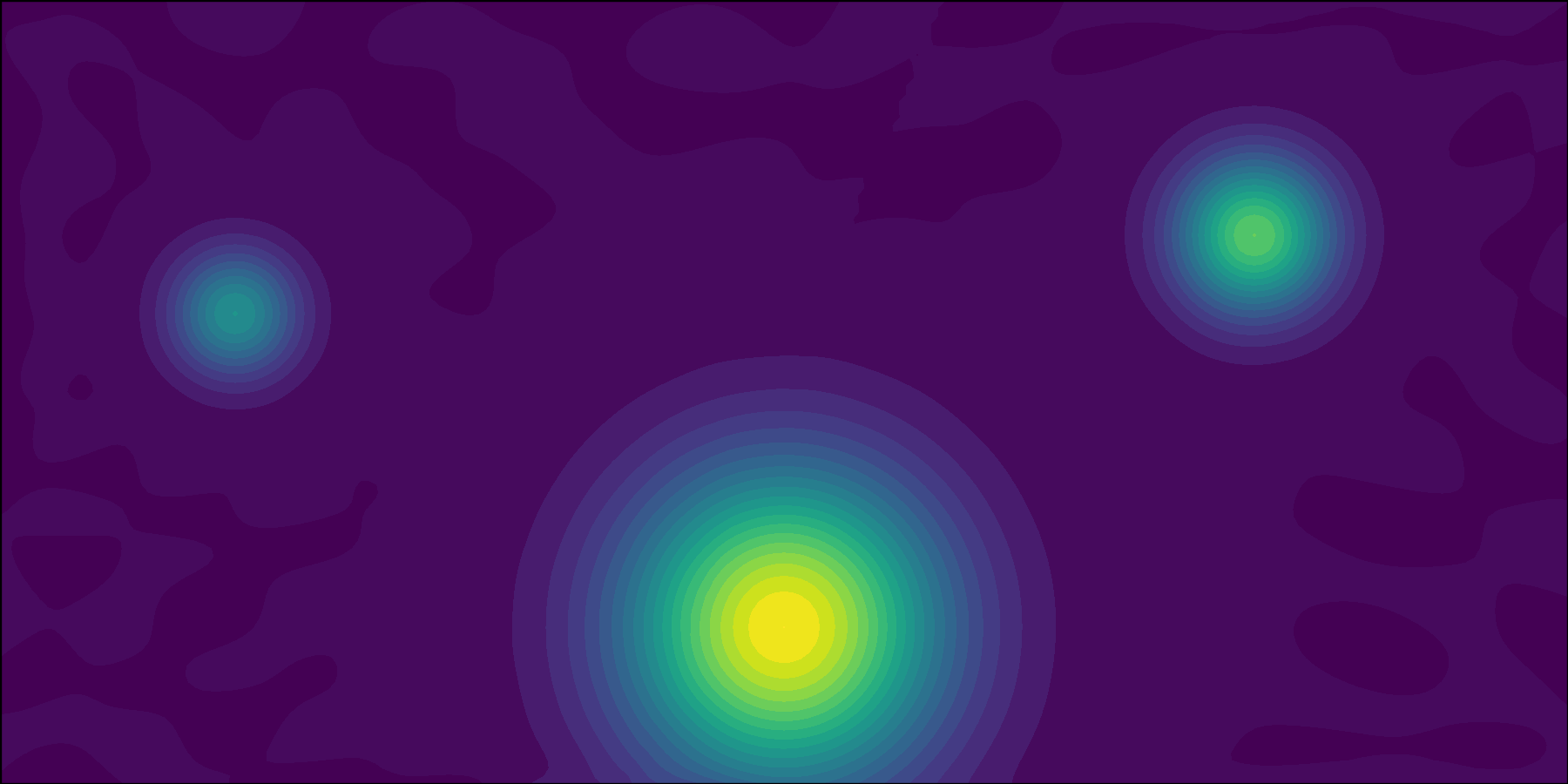}
            \addplot graphics[xmin=0, xmax=2, ymin=0, ymax=1]{\imagepath};

            \addplot[
                domain=0.8:1.2, samples=100, dashed, red,
            ]
            {(x-0.8)/0.4};


        \nextgroupplot[
            ylabel={$y$},
            xlabel={$x$},
            tick scale binop ={\times},
            xmin = 0., xmax = 2.0,
            ymin = 0., ymax = 1.0,
            colorbar,
            colormap name=viridis,
            point meta min=0.0,   
            point meta max=0.48,   
        ]
        
            \edef\imagepath{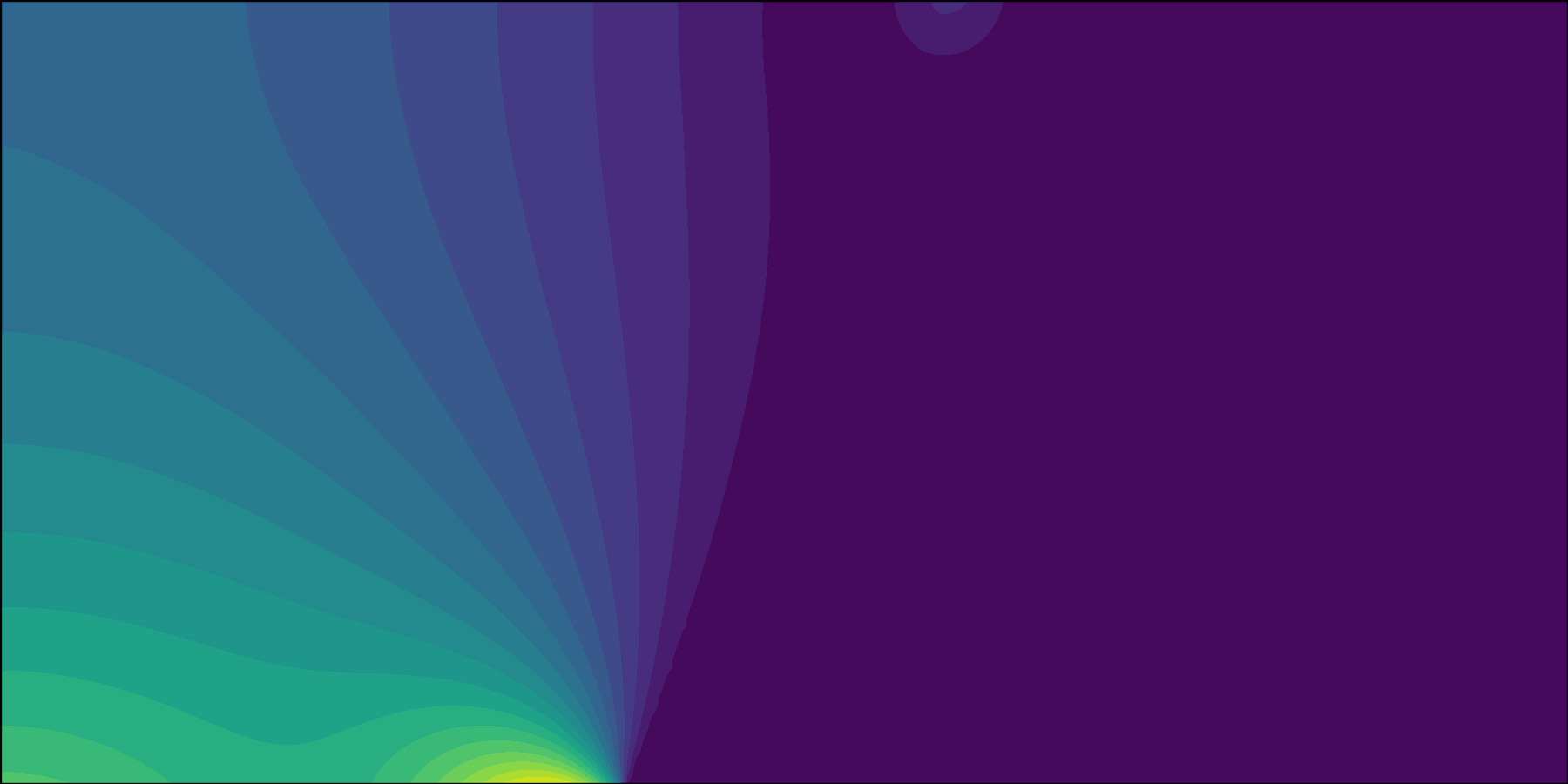}
            \addplot graphics[xmin=0., xmax=2., ymin=0, ymax=1.]{\imagepath};

            \addplot[
                domain=0.8:1.2, samples=100, dashed, red,
            ]
            {(x-0.8)/0.4};

        \nextgroupplot[
            xmin=0, xmax=0.8,
            ymin=0,
            xlabel={$x$},
            ylabel={$\bn\cdot\nabla u$},
            tick scale binop ={\times},
        ]
            \addplot+ [
                line width=1.0,
                solid,
                smooth,
                mark=none,
                blue,
            ]
            table [
                x index=0, y index=1,
            ]{data_2dpoisson_window.bl_neumann.txt};

  \end{groupplot}
\node[below = 1.5cm of my plots c1r1.south west,
    anchor=west,
] {(a) Trained solution \eqref{eq:window-ansatz}};
\node[below = 1.5cm of my plots c2r1.south west,
    anchor=west,
] {(b) $\LHS[u]$};
\node[below = 1.5cm of my plots c1r2.south west,
    anchor=west,
] {(c) Point-wise error $|u-u_{ref}|$};
\node[below = 1.5cm of my plots c2r2.south west,
    anchor=west,
] {(d) Bottom-left Neumann boundary};
\end{tikzpicture}
%
    \caption{Trained solution of \refB{windowing} approach \eqref{eq:window-ansatz} for Problem 4:
    (a) trained solution;
    (b) $\LHS[u]$;
    (c) point-wise absolute solution error; and
    (d) the solution slope on the bottom-left boundary.
    }
    \label{fig:2dpoisson_window}
\end{figure}
Figure~\ref{fig:2dpoisson_window} shows the trained solution of the windowing approach.
Overall, the solution qualitatively agrees with the finite-element ground truth,
achieving a relative $L_2$ error of $10.1\%$.
Figure~\ref{fig:2dpoisson_window}~(c) shows that the solution error
propagates from the bottom-left boundary,
where the Neumann condition is significantly violated, particularly for $x\in[0.5, 0.8]$.
This boundary uses the relaxed loss weight of $10^{-4}$;
however, in this study, increasing the weight further degraded the PDE loss
and ultimately yielded a larger overall error.
This region corresponds to where the interface solution ($u^{(I)}_{ij}$ in Eq. \ref{eq:window-ansatz}) overlaps with the Neumann boundary.
Although the interface ansatz exactly satisfies the interface condition,
its neural networks are defined only along the tangential coordinate of the interface,
which is not aligned with the Neumann boundary.
Consequently, the interface ansatz struggles to satisfy the Neumann boundary condition,
while the interior ansatz decays near the interface by design and cannot compensate for this deficiency.
\par
A simple workaround for this corner issue is to reduce the window size so that the interface ansatz lies entirely within the domain,
and augment the solution ansatz with the corner contribution \eqref{eq:window-corner}.
However, this results in imperfect overlap between the solution ansatz components,
which causes the optimization to become very stiff with poor convergence.
This behavior is consistent with the one-dimensional results in \ref{subsec:comparison-overlap}.
The issue is further aggravated when all boundary conditions are hard-constrained,
as the solution ansatz struggles to resolve the physics equation residual.
\refB{
We illustrate this failure with a fully hard-constrained variant in \ref{app:2dpoisson-full-window}.
In Figure~\ref{fig:2dpoisson_window_full_lhsrhs} of \ref{app:2dpoisson-full-window}, the trained residual then exhibits pronounced artificial oscillations near the endpoints of the window supports, a signature of the imperfect overlap.}
\par

The difficulty associated with corners, kinks, and interface–boundary intersections is well known in both classical numerical methods and more recent PINN formulations. A commonly cited explanation is the presence of singularities or reduced regularity in the true solution, for instance, in problems with reentrant corners, where derivatives may become unbounded or exhibit rapid variation near the corner \cite{dauge2006elliptic}. This can lead to poor accuracy and convergence in both finite element and neural network-based solvers.
In addition, even in the absence of singularities, geometric intersections can cause issues due to the multiple constraints being imposed along different directions, which can lead to stiffness in the resulting optimization problem \cite{sharma2023stiff}. 

A related issue arises from the interaction between the chosen ansatz and the differential operators associated with the form of the problem. For example, Sukumar and Srivastava~\cite{sukumar2022exact} developed a hard-constrained ansatz based on approximate distance functions to enforce boundary conditions exactly. Due to the strong form collocation nature of PINNs, the Laplacian of the distance function appears in the residual, which blows up for collocation points near the corner. To mitigate this issue, they avoid sampling collocation points within a prescribed neighborhood of the corner. Similar ideas appear in other PINN formulations, where spatial weighting or adaptive sampling reduces the influence of regions near singular or geometrically complex features \cite{lee2026least, zeng2024adaptive}. Notably, we take a similar approach in \ref{app:2dpoisson-full-window}, where the residual is multiplied by $r^2$ (where $r$ is the distance to the corner) in order to regularize the problem. 
\refB{While this weighting can regularize the corner singularity, the imperfect overlap of the windows still introduces spurious modes near the endpoints of their supports, degrading the solution accuracy in the present demonstration.}

\subsubsection{\refB{Two-dimensional demonstration of the buffer approach}}\label{subsubsec:buffer2d}
In the buffer formulation, the radial basis buffer function described in (\ref{eq:g-rbf}) was used,
with its degrees of freedom determined by the boundary/interface constraints (\ref{eq:g-bc-high-dim}).
\refB{As discussed in Section~\ref{subsec:buffer-multi-dim}, the constraints are enforced exactly only at the boundary and interface sample points, so the buffer method is \emph{discretely} hard-constrained in this two-dimensional setting.}
Boundary sample points on each edge were placed
following the guideline described in Section~\ref{subsec:buffer-multi-dim}:
four points for the Dirichlet boundaries and eight points elsewhere.
This choice provided a good balance between local enforcement and smooth interpolation across the sampled points.
\refB{The resulting coefficient systems are $40\times40$ and $28\times28$, corresponding to
$24$ Neumann and $16$ interface conditions in the left subdomain, and $12$ Dirichlet and
$16$ interface conditions in the right.}
\refA{The ratio parameters in the buffer interface conditions (\ref{eq:g-itf}) are set to $\gamma_0=\gamma_1=1$, corresponding to an even split between the two subdomains.}
\refB{Their condition numbers are $1.85\times10^{2}$ for the left subdomain and
$9.92\times10^{1}$ for the right, within the moderate range targeted by the guideline of
Section~\ref{subsec:buffer-multi-dim}.}
\par
\begin{figure}[tbh]
\hspace*{-1.8cm}
    \begin{tikzpicture}[
]
    \begin{groupplot}[
        group style={
            group name = my plots,
            group size= 2 by 2,
            xlabels at =edge bottom,
            horizontal sep=3cm,
            vertical sep=2.2cm,
        },
        height = 0.3\textwidth,
        width = 0.5\textwidth,
        name=chung,
    ]    
\pgfplotsset{set layers=standard}%

        \nextgroupplot[
            ylabel={$y$},
            xlabel={$x$},
            tick scale binop ={\times},
            xmin = 0, xmax = 2,
            ymin = 0, ymax = 1,
            colorbar,
            colormap name=viridis,
            point meta min=-3.159e-2,   
            point meta max=2.645,   
        ]
        
            \edef\imagepath{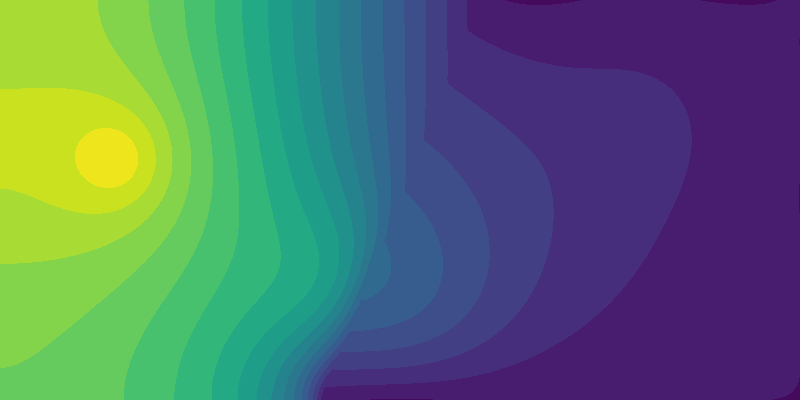}
            \addplot graphics[xmin=0, xmax=2, ymin=0, ymax=1]{\imagepath};

            \addplot[
                domain=0.8:1.2, samples=100, dashed, red,
            ]
            {(x-0.8)/0.4};

        \nextgroupplot[
            ylabel={$y$},
            xlabel={$x$},
            tick scale binop ={\times},
            xmin = 0, xmax = 2,
            ymin = 0, ymax = 1,
            colorbar,
            colormap name=viridis,
            point meta min=0,   
            point meta max=2.66,
            ,   
        ]
        
            \edef\imagepath{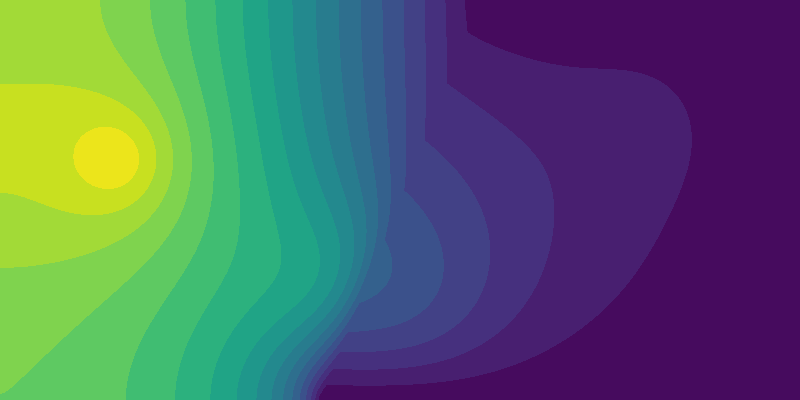}
            \addplot graphics[xmin=0, xmax=2, ymin=0, ymax=1]{\imagepath};

            \addplot[
                domain=0.8:1.2, samples=100, dashed, red,
            ]
            {(x-0.8)/0.4};

        \nextgroupplot[
            ylabel={$y$},
            xlabel={$x$},
            tick scale binop ={\times},
            xmin = 0, xmax = 2,
            ymin = 0, ymax = 1,
            colorbar,
            colormap name=viridis,
            point meta min=-1.013e-1,   
            point meta max=2.645,   
        ]
        
            \edef\imagepath{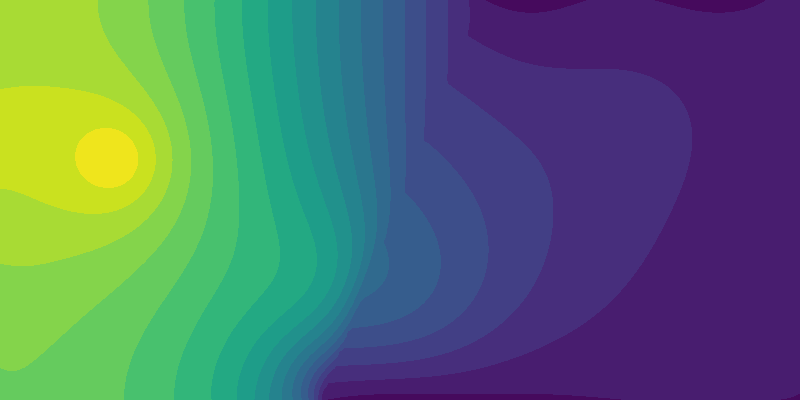}
            \addplot graphics[xmin=0, xmax=2, ymin=0, ymax=1]{\imagepath};

            \addplot[
                domain=0.8:1.2, samples=100, dashed, red,
            ]
            {(x-0.8)/0.4};

        \nextgroupplot[
            ylabel={$y$},
            xlabel={$x$},
            tick scale binop ={\times},
            xmin = 0, xmax = 2,
            ymin = 0, ymax = 1,
            colorbar,
            colormap name=viridis,
            point meta min=-5.906e-2,   
            point meta max=7.780e-2,
            ,   
        ]
        
            \edef\imagepath{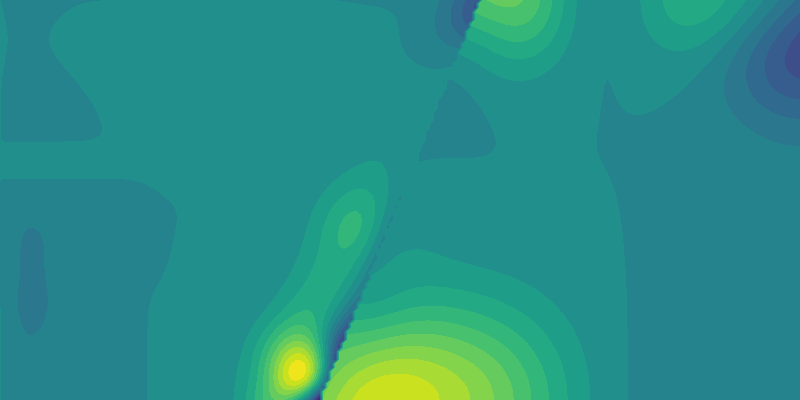}
            \addplot graphics[xmin=0, xmax=2, ymin=0, ymax=1]{\imagepath};

            \addplot[
                domain=0.8:1.2, samples=100, dashed, red,
            ]
            {(x-0.8)/0.4};

  \end{groupplot}
\node[below = 1.5cm of my plots c1r1.south west,
    anchor=west,
] {(a) Trained solution (\ref{eq:buffer})};
\node[below = 1.5cm of my plots c2r1.south west,
    anchor=west,
] {(b) Ground truth};
\node[below = 1.5cm of my plots c1r2.south west,
    anchor=west,
] {\refB{(c)} PINN term $NN_m(\bx)$};
\node[below = 1.5cm of my plots c2r2.south west,
    anchor=west,
] {\refB{(d)} Buffer term $g_m(\bx)$ (\ref{eq:g-rbf})};
\end{tikzpicture}
%
    \caption{Trained solution of buffer approach (Eq.~\ref{eq:buffer}) for Problem 4: (a) trained solution; (b) the ground truth; (c) PINN term $NN_m(\bx)$; and (d) buffer term $g_m(\bx)$ (Eq.~\ref{eq:g-rbf}).}
    \label{fig:2dpoisson_buffer}
\end{figure}
Figure~\ref{fig:2dpoisson_buffer} compares the trained solution of the buffer approach (\ref{eq:buffer}) with the finite-element ground truth.
The learned solution accurately captures the \refB{discontinuous slope} across the slanted interface and reproduces the asymmetric field distribution induced by the localized Gaussian sources.
The decomposition plots---showing the neural network term $NN_m(\bx)$ and the buffer term $g_m(\bx)$---illustrate that the buffer correction remains small (on the order of $10^{-2}$) relative to the dominant neural component, confirming that the buffer function primarily serves as a lightweight mechanism for enforcing the hard constraints without distorting the learned field.
\par
\begin{figure}[p]
    \begin{tikzpicture}[
font=\small,
]
    \begin{groupplot}[
        group style={
            group name = my plots,
            group size= 2 by 4,
            xlabels at =edge bottom,
            horizontal sep=2cm,
            vertical sep=2.2cm,
        },
        height = 0.25\textwidth,
        width = 0.45\textwidth,
        enlarge x limits={false, abs value = 5mm},
        enlarge y limits={true, abs value = 2mm},
        name=chung,
    ]    
\pgfplotsset{set layers=standard}%

        \nextgroupplot[
            xlabel={$x$},
            ylabel={$\partial_y u$},
            tick scale binop ={\times},
            ymin=-0.75, ymax=0.75,
        ]
            \addplot[black, thin] coordinates {(0,0) (0.8,0)};

            \addplot+ [
                line width=1.0,
                solid,
                smooth,
                mark=none,
                blue,
            ]
            table [
                x index=0, y index=5,
            ]{data_2dpoisson_pred.bnd_line.txt};

            \addplot+ [
                line width=1.0,
                draw=none,
                smooth,
                mark=*,
                mark options={fill=white,},
                red,
            ]
            table [
                x index=0, y index=1,
            ]{data_2dpoisson_pred.bnd_pt.BL.txt};

        \nextgroupplot[
            xlabel={$x$},
            ylabel={$u$},
            tick scale binop ={\times},
        ]

            \addplot[black, thin] coordinates {(0.8,0) (2,0)};

            \addplot+ [
                line width=1.0,
                solid,
                smooth,
                mark=none,
                blue,
            ]
            table [
                x index=1, y index=6,
            ]{data_2dpoisson_pred.bnd_line.txt};

            \addplot+ [
                line width=1.0,
                draw=none,
                smooth,
                mark=*,
                mark options={fill=white,},
                red,
            ]
            table [
                x index=0, y index=1,
            ]{data_2dpoisson_pred.bnd_pt.BR.txt};

        \nextgroupplot[
            xlabel={$x$},
            ylabel={$\partial_y u$},
            tick scale binop ={\times},
        ]
            \addplot[black, thin] coordinates {(0,0) (1.2,0)};

            \addplot+ [
                line width=1.0,
                solid,
                smooth,
                mark=none,
                blue,
            ]
            table [
                x index=2, y index=7,
            ]{data_2dpoisson_pred.bnd_line.txt};

            \addplot+ [
                line width=1.0,
                draw=none,
                smooth,
                mark=*,
                mark options={fill=white,},
                red,
            ]
            table [
                x index=0, y index=1,
            ]{data_2dpoisson_pred.bnd_pt.TL.txt};

        \nextgroupplot[
            xlabel={$x$},
            ylabel={$u$},
            tick scale binop ={\times},
        ]
            \addplot[black, thin] coordinates {(1.2,0) (2,0)};

            \addplot+ [
                line width=1.0,
                solid,
                smooth,
                mark=none,
                blue,
            ]
            table [
                x index=3, y index=8,
            ]{data_2dpoisson_pred.bnd_line.txt};

            \addplot+ [
                line width=1.0,
                draw=none,
                smooth,
                mark=*,
                mark options={fill=white,},
                red,
            ]
            table [
                x index=0, y index=1,
            ]{data_2dpoisson_pred.bnd_pt.TR.txt};

        \nextgroupplot[
            xlabel={$y$},
            ylabel={$\partial_x u$},
            tick scale binop ={\times},
        ]
            \addplot[black, thin] coordinates {(0,0) (1,0)};

            \addplot+ [
                line width=1.0,
                solid,
                smooth,
                mark=none,
                blue,
            ]
            table [
                x index=4, y index=9,
            ]{data_2dpoisson_pred.bnd_line.txt};

            \addplot+ [
                line width=1.0,
                draw=none,
                smooth,
                mark=*,
                mark options={fill=white,},
                red,
            ]
            table [
                x index=0, y index=1,
            ]{data_2dpoisson_pred.bnd_pt.L.txt};

        \nextgroupplot[
            xlabel={$y$},
            ylabel={$u$},
            tick scale binop ={\times},
        ]
            \addplot[black, thin] coordinates {(0,0) (1,0)};

            \addplot+ [
                line width=1.0,
                solid,
                smooth,
                mark=none,
                blue,
            ]
            table [
                x index=4, y index=10,
            ]{data_2dpoisson_pred.bnd_line.txt};

            \addplot+ [
                line width=1.0,
                draw=none,
                smooth,
                mark=*,
                mark options={fill=white,},
                red,
            ]
            table [
                x index=0, y index=1,
            ]{data_2dpoisson_pred.bnd_pt.R.txt};

        \nextgroupplot[
            xlabel={$y$},
            ylabel={$\jump{u}$},
            tick scale binop ={\times},
        ]
            \addplot[black, thin] coordinates {(0,0) (1,0)};

            \addplot+ [
                line width=1.0,
                solid,
                smooth,
                mark=none,
                blue,
            ]
            table [
                x index=4,
                y expr=\thisrowno{11} - \thisrowno{12},
            ]{data_2dpoisson_pred.bnd_line.txt};

            \addplot+ [
                line width=1.0,
                draw=none,
                smooth,
                mark=*,
                mark options={fill=white,},
                red,
            ]
            table [
                x index=0, y index=1,
            ]{data_2dpoisson_pred.bnd_pt.ITF-sol.txt};


        \nextgroupplot[
            xlabel={$y$},
            ylabel={$\avg{\bn\cdot\kappa\nabla u}$},
            tick scale binop ={\times},
        ]
            \addplot[black, thin] coordinates {(0,0) (1,0)};

            \addplot+ [
                line width=1.0,
                solid,
                smooth,
                mark=none,
                blue,
            ]
            table [
                x index=4, y index=13,
            ]{data_2dpoisson_pred.bnd_line.txt};

            \addplot+ [
                line width=1.0,
                draw=none,
                smooth,
                mark=*,
                mark options={fill=white,},
                red,
            ]
            table [
                x index=0, y index=1,
            ]{data_2dpoisson_pred.bnd_pt.ITF-flux.txt};

  \end{groupplot}
\node[below = 1.2cm of my plots c1r1.south west,
    anchor=west, xshift=-5mm,
] {(a) Bottom left, Neumann};
\node[below = 1.2cm of my plots c2r1.south west,
    anchor=west, xshift=-5mm,
] {(b) Bottom right, Dirichlet};
\node[below = 1.2cm of my plots c1r2.south west,
    anchor=west, xshift=-5mm,
] {(c) Top left, Neumann};
\node[below = 1.2cm of my plots c2r2.south west,
    anchor=west, xshift=-5mm,
] {(d) Top right, Dirichlet};
\node[below = 1.2cm of my plots c1r3.south west,
    anchor=west, xshift=-5mm,
] {(e) Left, Neumann};
\node[below = 1.2cm of my plots c2r3.south west,
    anchor=west, xshift=-5mm,
] {(f) Right, Dirichlet};
\node[below = 1.2cm of my plots c1r4.south west,
    anchor=west, xshift=-5mm,
] {(g) Interface, solution};
\node[below = 1.2cm of my plots c2r4.south west,
    anchor=west, xshift=-5mm,
] {(h) Interface, flux};
\end{tikzpicture}
%
    \caption{Boundary enforcement of the buffer approach:
    (a) the Neumann condition on the bottom-left boundary;
    (b) the Dirichlet condition on the bottom-right boundary;
    (c) the Neumann condition on the top-left boundary;
    (d) the Dirichlet condition on the top-right boundary;
    (e) the Neumann condition on the left boundary;
    (f) the Dirichlet condition on the right boundary;
    (g) the solution condition on the interface; and
    (h) the flux condition on the interface.
    Blue lines indicate the trained solution,
    and red markers indicate the boundary sample points
    where the buffer term hard-constrains the boundary conditions.}
    \label{fig:2dpoisson_buffer_bnd}
\end{figure}
Boundary and interface conditions are analyzed in Figure~\ref{fig:2dpoisson_buffer_bnd}.
The results confirm precise satisfaction of the boundary conditions across all domain edges.
The buffer corrections strictly enforce boundary conditions at the sampled points (red markers), while remaining smooth between samples.
Figure~\ref{fig:2dpoisson_buffer_bnd}~(g–h) verifies that both the solution and the flux continuity are satisfied along the interface, with deviations \refB{$\lesssim10^{-2}$},
demonstrating effective enforcement of interfacial physics.
\refB{The corresponding root-mean-square errors are reported in Table~\ref{tab:2dpoisson-bnd-rmse},
evaluated at $100$ uniformly spaced points on each boundary and interface side
that do not coincide with the buffer sample points.}
\par
\begin{table}[tbh]
    \centering
    \small
    \setlength{\tabcolsep}{4pt}
    \begin{tabular}{|l|c|c|c|c|c|c|c|c|}
        \hline
        Side
        & BL & BR & TL & TR & L & R
        & $\jump{u}$ & $\avg{\bn\cdot\kappa\nabla u}$ \\
        \hline
        RMSE
        & 2.18e-2 & 7.40e-3 & 6.40e-4 & 9.22e-3 & 8.95e-4 & 7.71e-4
        & 1.94e-2 & 8.23e-3 \\
        \hline
    \end{tabular}
    \caption{\refB{Root-mean-square error of the boundary and interface conditions for the buffer approach on Problem 4. B, T, L, and R denote the bottom, top, left, and right edges.}}
    \label{tab:2dpoisson-bnd-rmse}
\end{table}
\par
\begin{table}[tbh] 
    \begin{tabular}{|l|c|c|c|} 
        \hline 
        Problem 4 & Layers & AF & Relative $L_2$ error \\ 
        \hline 
        $\phi$-PINNs & $[3, 25, 25, 25, 1]$ & $\tanh$ & 8.60e-3 \\
        I-PINNs & $[2, 25, 25, 25, 1]$ & swish, $\tanh$ & 7.81e-3 \\
        AdaI-PINNs & $[2, 25, 25, 25, 1]$ & Adaptive $\tanh$ & 
        8.33e-1 \\
        M-PINN & $[2, 25, 25, 25, 1]\times2$ & $\tanh$ & 7.34e-3 \\
        \hline
        Window$^*$ & $[2, 25, 25, 25, 1]\times2$ & $\tanh$ & 1.01e-1 \\
        Buffer & $[2, 25, 25, 25, 1]\times2$ & $\tanh$ & 3.48e-3 \\
        \hline 
    \end{tabular}
    \centering 
    \caption{Comparison between different PINN approaches on problem 4.
    Note that, for windowing approach, only the interface is hard-constrained while all the other boundary conditions are soft-constrained.
    } 
    \label{tab:init-prob4} 
\end{table}
Table~\ref{tab:init-prob4} summarizes the relative $L_2$ errors across all tested formulations.
For a consistent comparison, all models share the training setup described above
for the windowing approach.
For the soft-constrained methods, the interface losses are evaluated at
$40$ uniformly distributed collocation points along the interface,
excluding the corner points.
Despite the overlap difficulties discussed above, the windowing approach with interface-only hard constraints achieves a relatively poor error at $10.1\%$.
Among the soft-constrained baselines, AdaI-PINNs shows degraded performance in this two-dimensional case,
with relative $L_2$ errors exceeding $80\%$ due to its difficulty in resolving the flux discontinuity and Neumann conditions.
The remaining soft-constrained methods perform well, all achieving errors of $\gtrsim 7\times10^{-3}$.
By contrast, the buffer approach maintains a low error of \refB{$3.48\times10^{-3}$}.
Although the improvement over the soft-constrained methods is only modest,
this shows that the buffer approach provides a more straightforward and stable path for practical implementation.

\subsubsection{\refB{Sensitivity to the hyperparameters in the buffer approach}}\label{subsubsec:buffer2d-sensitivity}
\refB{
The buffer approach introduces its own hyperparameters for discrete hard-constraining: the number of sample points at which the conditions \eqref{eq:g-bc-high-dim} are enforced exactly.
In addition, the nominal radius of the radial basis functions \eqref{eq:g-rbf} attached to each sample point also affects the training performance and the solution accuracy.
We assessed their influence on Problem~4 through a sweep over $72{,}000$ configurations.
While a local sensitivity study around the best configuration gives no clear guidance,
the global trend across the full sweep shows that the training loss correlates strongly with the condition number of the buffer coefficient system,
supporting the conditioning-based guideline of Section~\ref{subsec:buffer-multi-dim}.
The correlation with the solution accuracy itself is not monotone: the best-conditioned configurations use the fewest sample points and the narrowest basis functions, and therefore enforce the constraints only loosely away from the samples,
so the lowest error is attained at an intermediate condition number rather than at the best-conditioned end of the sweep.
Choosing these hyperparameters is thus a balance between the conditioning of the
coefficient system and the accuracy with which the buffer functions represent the boundary
and interface conditions.
}
\refA{The interface split ratios $\gamma_0$ and $\gamma_1$ (\ref{eq:g-bc-high-dim}c--d) were also
included in the same sweep, and little sensitivity was observed to either.
This is expected, since they alter only the right-hand side of the buffer system
and thus leave its conditioning unchanged.
Full details of the sensitivity study are reported in~\ref{app:buffer2d-sensitivity}.}

\subsubsection{\refA{Buffer approach on a curved interface}}\label{subsubsec:buffer2d-petal}
\refA{
We further demonstrate the buffer approach on a more complex, curved interface.
Let $\Omega = [-1, 1]^{2}$ be partitioned by a two-lobe interface
\begin{equation}
\Gamma_{\mathrm{int}} = \{(x, y) \in \Omega \, : \, x^{2} + y^{2} = r(\theta)^{2}\},
\end{equation}
whose radius is parametrized in the polar angle $\theta = \operatorname{atan2}(y, x)$ as
\begin{equation}
r(\theta) = 0.4 + 0.2 \sin(2\theta).
\end{equation}
The interior lobe $\Omega_{1}$ and the exterior region $\Omega_{2}$ carry piecewise-constant diffusivities $\kappa_{1} = 0.1$ and $\kappa_{2} = 1$, respectively, and are driven by
\begin{equation}
    f(\bx) =
    \begin{cases}
        4 & \bx \in \Omega_{1}, \\
        1/(x^{2} + y^{2}) & \bx \in \Omega_{2}.
    \end{cases}
\end{equation}
Homogeneous Dirichlet conditions are imposed on the outer boundary,
\begin{subequations}
\begin{equation}
    u = 0 \qquad\text{on } \partial\Omega,
\end{equation}
and standard continuity of the solution and the normal flux is enforced across the curved interface,
\begin{align}
    \jump{u} &= 0 &\text{on } \Gamma_{\mathrm{int}}, \\
    \jump{\bn \cdot \kappa \nabla u} &= 0 &\text{on } \Gamma_{\mathrm{int}}.
\end{align}
\end{subequations}
The reference solution is produced by a standalone finite-element solve on the same geometry.}
\par
\refA{
Each subdomain is represented by a fully connected network with hidden layers
$[10, 10, 10]$ and $\tanh$ activation, initialized with the Glorot uniform scheme
scaled by $0.1$.
Training is performed for $5\times10^4$ iterations with the Adam optimizer~\cite{kingma2014adam},
with the learning rate cosine-decayed from $10^{-2}$ to $10^{-4}$.
The PDE loss \eqref{eq:J-int} is evaluated on a $40\times40$ uniform grid over the
bounding box, yielding $240$ collocation points in the interior lobe $\Omega_1$ and
$1360$ in the exterior region $\Omega_2$.
The buffer functions are assigned at $8$ Gauss--Legendre sample points on
each of the four Dirichlet boundary edges and $24$ sample points along the two
interface lobes.
The radial basis widths in (\ref{eq:g-rbf}) follow the guideline of
Section~\ref{subsec:buffer-multi-dim}, with the interface width scaled by $0.6$,
and the interface ratio parameters are set to $\gamma_0=\gamma_1=10$.
This configuration attains a relative $L_2$ error of $1.58\times10^{-3}$,
compared with $3.76\times10^{-3}$ for M-PINN, trained with the same architecture and training budget.
}
\par
\begin{figure}[p]
    \begin{tikzpicture}[
]
    \begin{groupplot}[
        group style={
            group name = my plots,
            group size= 2 by 2,
            xlabels at =edge bottom,
            horizontal sep=3.2cm,
            vertical sep=2.3cm,
        },
        height = 0.4\textwidth,
        width  = 0.4\textwidth,
        colorbar style={xshift=-3mm},
        name=chung,
    ]
\pgfplotsset{set layers=standard}%

        \nextgroupplot[
            ylabel={$y$},
            xlabel={$x$},
            tick scale binop ={\times},
            xmin = -1, xmax = 1,
            ymin = -1, ymax = 1,
            colorbar,
            colormap name=viridis,
            point meta min=-1.384e-16,   
            point meta max=2.264,        
        ]
            \edef\imagepath{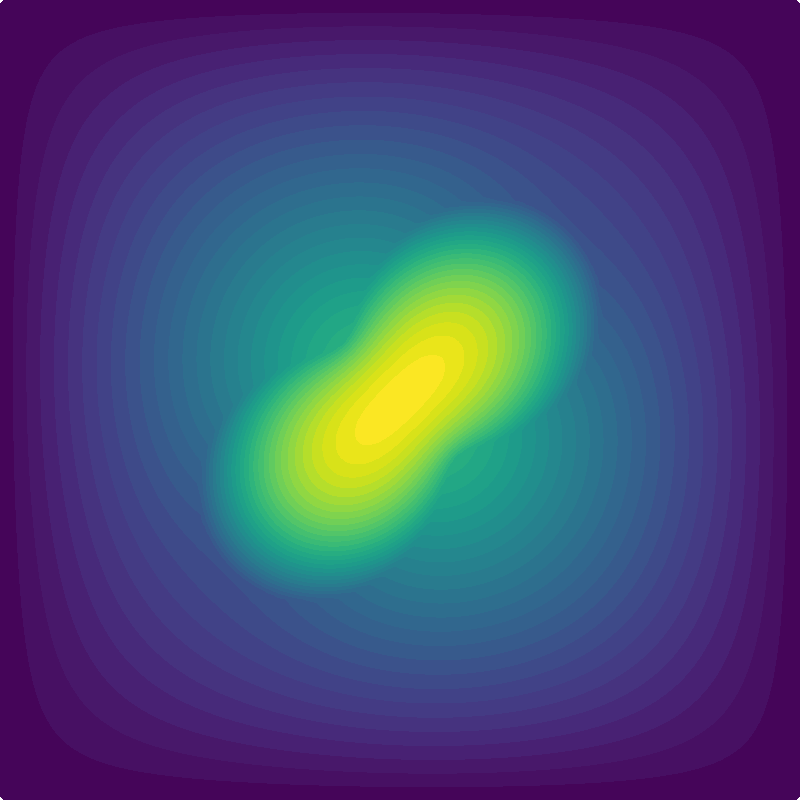}
            \addplot graphics[xmin=-1, xmax=1, ymin=-1, ymax=1]{\imagepath};

            \addplot+ [
                line width=0.6pt, dashed, black, mark=none, smooth,
            ]
            table [x index=0, y index=1, skip first n=1]
                {data_petal_interface.txt};

        \nextgroupplot[
            ylabel={$y$},
            xlabel={$x$},
            tick scale binop ={\times},
            xmin = -1, xmax = 1,
            ymin = -1, ymax = 1,
            colorbar,
            colormap name=viridis,
            point meta min=-9.363e-4,
            point meta max=2.260,
        ]
            \edef\imagepath{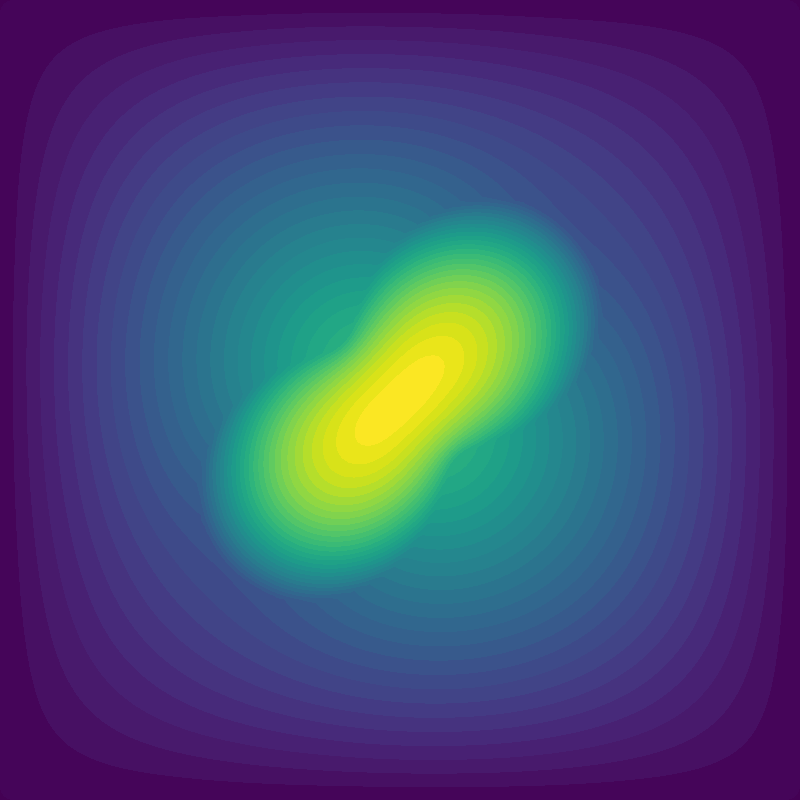}
            \addplot graphics[xmin=-1, xmax=1, ymin=-1, ymax=1]{\imagepath};

            \addplot+ [
                line width=0.6pt, dashed, black, mark=none, smooth,
            ]
            table [x index=0, y index=1, skip first n=1]
                {data_petal_interface.txt};

        \nextgroupplot[
            ylabel={$y$},
            xlabel={$x$},
            tick scale binop ={\times},
            xmin = -1, xmax = 1,
            ymin = -1, ymax = 1,
            colorbar,
            colormap name=viridis,
            point meta min=-1.122e-2,
            point meta max=1.671e-2,
        ]
            \edef\imagepath{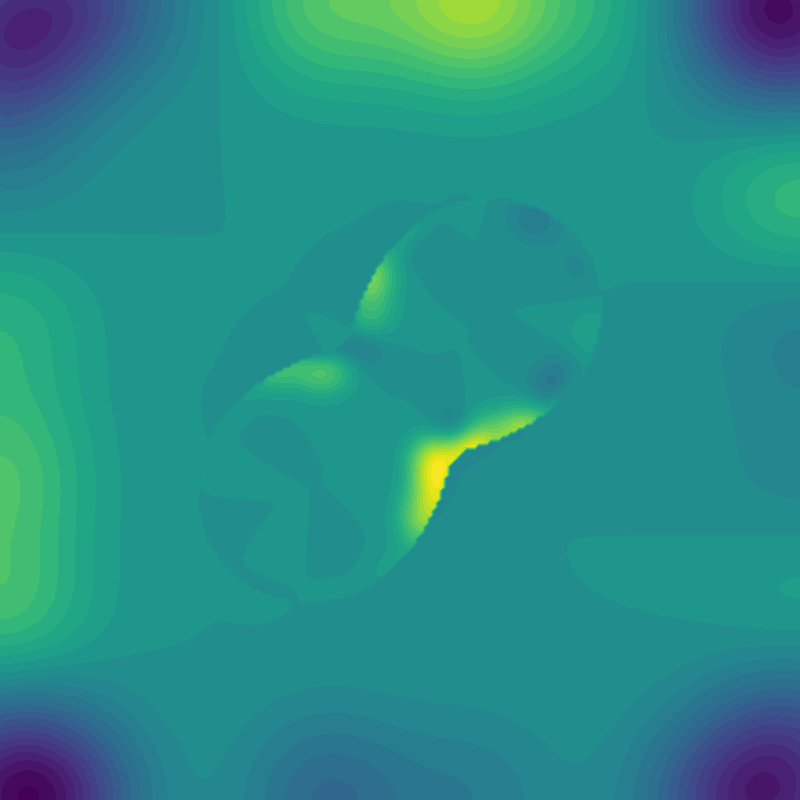}
            \addplot graphics[xmin=-1, xmax=1, ymin=-1, ymax=1]{\imagepath};

            \addplot+ [
                line width=0.6pt, dashed, black, mark=none, smooth,
            ]
            table [x index=0, y index=1, skip first n=1]
                {data_petal_interface.txt};

        \nextgroupplot[
            ylabel={$y$},
            xlabel={$x$},
            tick scale binop ={\times},
            xmin = -1, xmax = 1,
            ymin = -1, ymax = 1,
            colorbar,
            colormap name=viridis,
            point meta min=9.433e-9,
            point meta max=1.755e-2,
        ]
            \edef\imagepath{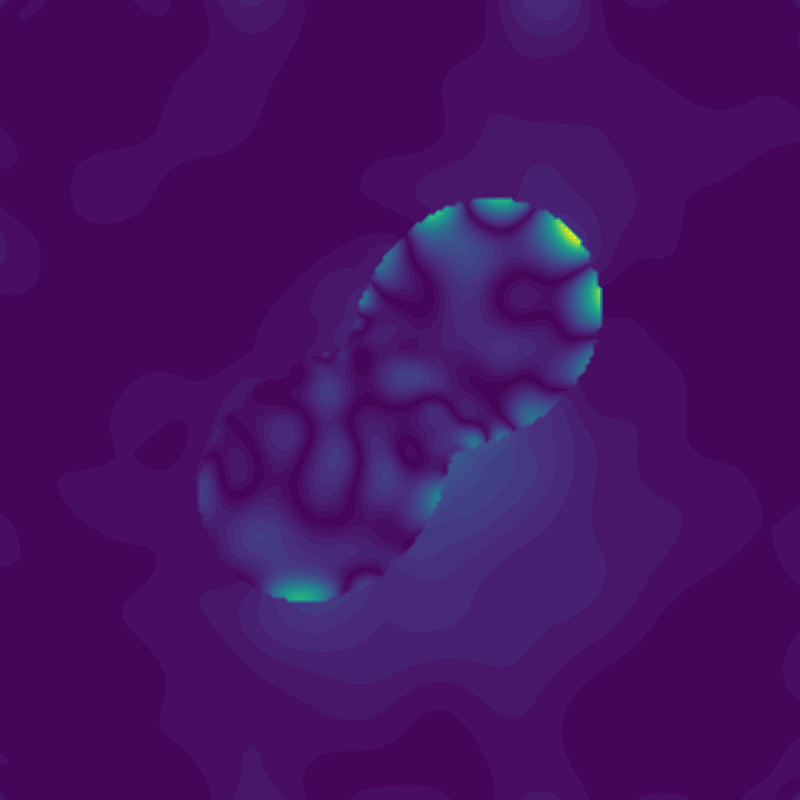}
            \addplot graphics[xmin=-1, xmax=1, ymin=-1, ymax=1]{\imagepath};

            \addplot+ [
                line width=0.6pt, dashed, black, mark=none, smooth,
            ]
            table [x index=0, y index=1, skip first n=1]
                {data_petal_interface.txt};

  \end{groupplot}
\node[below = 1.4cm of my plots c1r1.south west,
    anchor=west, xshift=-10pt,
] {(a) Ground truth};
\node[below = 1.4cm of my plots c2r1.south west,
    anchor=west, xshift=-10pt,
] {(b) Trained solution (\ref{eq:buffer})};
\node[below = 1.4cm of my plots c1r2.south west,
    anchor=west, xshift=-10pt,
] {(c) Buffer term $g_m(\bx)$ (\ref{eq:g-rbf})};
\node[below = 1.4cm of my plots c2r2.south west,
    anchor=west, xshift=-10pt,
] {(d) Point-wise error $|U-u_{\mathrm{exact}}|$};
\end{tikzpicture}
%
\vspace*{-.5cm}
    \caption{
        \refA{Buffer approach on the petal-shaped interface problem: (a) finite-element ground truth; (b) trained solution (Eq.~\ref{eq:buffer}); (c) buffer term $g_{m}(\bx)$ (Eq.~\ref{eq:g-rbf}); and (d) point-wise solution error $|U-u_{\mathrm{exact}}|$. The curved interface $\Gamma_{\mathrm{int}}$ is overlaid as a dashed black line.}
    }
    \label{fig:petal_buffer}
\end{figure}
\refA{Figure~\ref{fig:petal_buffer} shows that the buffer approach reproduces the two-lobe finite-element solution accurately across the curved interface, with the point-wise error concentrated near $\Gamma_{\mathrm{int}}$ and remaining below $\sim2\times10^{-2}$. The buffer term itself stays on the order of $10^{-2}$ and localizes at the boundary and interface sample points, so the dominant contribution to $U$ comes from the neural component while the buffer serves only as a lightweight correction.}
\par
\begin{figure}[p]
    \begin{tikzpicture}[
font=\small,
]
    \begin{groupplot}[
        group style={
            group name = my plots,
            group size= 2 by 3,
            xlabels at =edge bottom,
            horizontal sep=2cm,
            vertical sep=2.2cm,
        },
        height = 0.25\textwidth,
        width  = 0.45\textwidth,
        enlarge x limits={false, abs value = 5mm},
        enlarge y limits={true, abs value = 2mm},
        name=chung,
    ]
\pgfplotsset{set layers=standard}%

        \nextgroupplot[
            xlabel={$y$},
            ylabel={$u_{2}(x=-1,y)$},
            tick scale binop ={\times},
        ]
            \addplot[black, thin] coordinates {(-1,0) (1,0)};

            \addplot+ [
                line width=1.0,
                solid,
                smooth,
                mark=none,
                blue,
            ]
            table [x index=0, y index=1, skip first n=1]
                {data_petal_pred.bnd_line.L.txt};

            \addplot+ [
                line width=1.0,
                draw=none,
                smooth,
                only marks,
                mark=*,
                mark size=1.2pt,
                mark options={fill=red, draw=red},
                red,
            ]
            table [x index=0, y index=1, skip first n=1]
                {data_petal_pred.bnd_pt.L.txt};

        \nextgroupplot[
            xlabel={$y$},
            ylabel={$u_{2}(x=1,y)$},
            tick scale binop ={\times},
        ]
            \addplot[black, thin] coordinates {(-1,0) (1,0)};

            \addplot+ [
                line width=1.0,
                solid,
                smooth,
                mark=none,
                blue,
            ]
            table [x index=0, y index=1, skip first n=1]
                {data_petal_pred.bnd_line.R.txt};

            \addplot+ [
                line width=1.0,
                draw=none,
                smooth,
                only marks,
                mark=*,
                mark size=1.2pt,
                mark options={fill=red, draw=red},
                red,
            ]
            table [x index=0, y index=1, skip first n=1]
                {data_petal_pred.bnd_pt.R.txt};

        \nextgroupplot[
            xlabel={$x$},
            ylabel={$u_{2}(x,y=-1)$},
            tick scale binop ={\times},
        ]
            \addplot[black, thin] coordinates {(-1,0) (1,0)};

            \addplot+ [
                line width=1.0,
                solid,
                smooth,
                mark=none,
                blue,
            ]
            table [x index=0, y index=1, skip first n=1]
                {data_petal_pred.bnd_line.B.txt};

            \addplot+ [
                line width=1.0,
                draw=none,
                smooth,
                only marks,
                mark=*,
                mark size=1.2pt,
                mark options={fill=red, draw=red},
                red,
            ]
            table [x index=0, y index=1, skip first n=1]
                {data_petal_pred.bnd_pt.B.txt};

        \nextgroupplot[
            xlabel={$x$},
            ylabel={$u_{2}(x,y=1)$},
            tick scale binop ={\times},
        ]
            \addplot[black, thin] coordinates {(-1,0) (1,0)};

            \addplot+ [
                line width=1.0,
                solid,
                smooth,
                mark=none,
                blue,
            ]
            table [x index=0, y index=1, skip first n=1]
                {data_petal_pred.bnd_line.T.txt};

            \addplot+ [
                line width=1.0,
                draw=none,
                smooth,
                only marks,
                mark=*,
                mark size=1.2pt,
                mark options={fill=red, draw=red},
                red,
            ]
            table [x index=0, y index=1, skip first n=1]
                {data_petal_pred.bnd_pt.T.txt};

        \nextgroupplot[
            xlabel={$\theta$},
            ylabel={$\jump{u}$},
            tick scale binop ={\times},
            xmin=-3.14159, xmax=3.14159,
            xtick={-3.14159, 0, 3.14159},
            xticklabels={$-\pi$, $0$, $\pi$},
        ]
            \addplot[black, thin] coordinates {(-3.14159,0) (3.14159,0)};

            \addplot+ [
                line width=1.0,
                solid,
                smooth,
                mark=none,
                blue,
            ]
            table [x index=0, y index=1, skip first n=1]
                {data_petal_pred.bnd_line.ITF-sol.txt};

            \addplot+ [
                line width=1.0,
                draw=none,
                smooth,
                only marks,
                mark=*,
                mark size=1.2pt,
                mark options={fill=red, draw=red},
                red,
            ]
            table [x index=0, y index=1, skip first n=1]
                {data_petal_pred.bnd_pt.ITF-sol.txt};

        \nextgroupplot[
            xlabel={$\theta$},
            ylabel={$\jump{\bn\cdot\kappa\nabla u}$},
            tick scale binop ={\times},
            xmin=-3.14159, xmax=3.14159,
            xtick={-3.14159, 0, 3.14159},
            xticklabels={$-\pi$, $0$, $\pi$},
        ]
            \addplot[black, thin] coordinates {(-3.14159,0) (3.14159,0)};

            \addplot+ [
                line width=1.0,
                solid,
                smooth,
                mark=none,
                blue,
            ]
            table [x index=0, y index=1, skip first n=1]
                {data_petal_pred.bnd_line.ITF-flux.txt};

            \addplot+ [
                line width=1.0,
                draw=none,
                smooth,
                only marks,
                mark=*,
                mark size=1.2pt,
                mark options={fill=red, draw=red},
                red,
            ]
            table [x index=0, y index=1, skip first n=1]
                {data_petal_pred.bnd_pt.ITF-flux.txt};

  \end{groupplot}
\node[below = 1.3cm of my plots c1r1.south west,
    anchor=west, xshift=-5mm,
] {(a) Left edge, Dirichlet};
\node[below = 1.3cm of my plots c2r1.south west,
    anchor=west, xshift=-5mm,
] {(b) Right edge, Dirichlet};
\node[below = 1.3cm of my plots c1r2.south west,
    anchor=west, xshift=-5mm,
] {(c) Bottom edge, Dirichlet};
\node[below = 1.3cm of my plots c2r2.south west,
    anchor=west, xshift=-5mm,
] {(d) Top edge, Dirichlet};
\node[below = 1.3cm of my plots c1r3.south west,
    anchor=west, xshift=-5mm,
] {(e) Interface, solution};
\node[below = 1.3cm of my plots c2r3.south west,
    anchor=west, xshift=-5mm,
] {(f) Interface, flux};
\end{tikzpicture}
%
    \caption{\refA{Boundary and interface enforcement of the buffer approach on the petal problem:
    (a--d) the Dirichlet condition on the left, right, bottom, and top edges of $\partial\Omega$;
    (e) the solution jump $\jump{u}$ across the interface parametrized by $\theta$; and
    (f) the flux jump $\jump{\bn\cdot\kappa\nabla u}$ across the interface.
    Blue lines indicate the trained solution, and red markers indicate the boundary/interface sample points where the buffer term hard-constrains the constraints.}}
    \label{fig:petal_buffer_bnd}
\end{figure}
\refA{Figure~\ref{fig:petal_buffer_bnd} confirms that the Dirichlet condition on $\partial\Omega$ and the solution/flux continuity across the curved interface are enforced with the same qualitative behaviour observed on the slanted-interface problem: the buffer collocation points sit at machine precision, and the between-sample deviations stay small ($\sim10^{-3}$ on the boundaries and $\sim10^{-2}$ on the interface), despite the geometry being no longer aligned with any coordinate axis.}
\refA{The corresponding root-mean-square errors are reported in
Table~\ref{tab:petal-bnd-rmse}, evaluated at $401$ uniformly spaced points on each outer
edge and $801$ along the interface.}
\par
\begin{table}[tbh]
    \centering
    \small
    \setlength{\tabcolsep}{6pt}
    \begin{tabular}{|l|c|c|c|c|c|c|}
        \hline
        Side & L & R & B & T
        & $\jump{u}$ & $\jump{\bn\cdot\kappa\nabla u}$ \\
        \hline
        RMSE & 3.62e-4 & 3.73e-4 & 5.13e-4 & 6.92e-4 & 4.94e-3 & 1.29e-2 \\
        \hline
    \end{tabular}
    \caption{\refA{Root-mean-square error of the boundary and interface conditions for the
    buffer approach on the petal-shaped interface problem.
    L, R, B, and T denote the left, right, bottom, and top edges of $\partial\Omega$}.}
    \label{tab:petal-bnd-rmse}
\end{table}
\par
The two-dimensional experiment highlights the effectiveness of the buffer approach.
By embedding interface and boundary conditions into a low-rank buffer representation,
the method successfully avoids penalty-weight tuning while maintaining fast convergence.
These results suggest that the buffer-based formulation provides a promising foundation for extending hard-constrained PINNs to multidimensional and multiphysics systems involving complex interfaces.

\subsubsection{\refB{Computational overhead of the buffer approach}}\label{subsubsec:buffer2d-timing}
\refB{
The buffer approach determines the coefficients $\bc(\theta_m)$ from a linear system at
every training iteration, which the soft-constrained methods do not require. We therefore
quantify the resulting computational overhead on the petal problem, both against an
M-PINN model of identical network size and as a function of the number of boundary and
interface sample points.
}
\par
\refB{
All timings were measured on a single node of the Livermore Computing Dane cluster~\cite{llnl_dane},
equipped with two
Intel Xeon Platinum 8480+ processors at $2.0$~GHz and $256$~GB of memory, using the same
JAX implementation~\cite{jax2018github} employed for all results in this work. Each
component was compiled once and then profiled over $100$
repeated calls with the network parameters redrawn at each call.
Markers report the median and the bars the full range over the profiled calls.
The sweep is centered on the best-performing
configuration of Section~\ref{subsubsec:buffer2d-petal}, varying one of the two sample
counts at a time---$N_{bnd}$ per boundary edge and $N_{itf}$ per interface lobe---while
holding the network size, collocation grid, and all other settings fixed.
}
\begin{figure}[tbh]
    \begin{tikzpicture}[
font=\small,
]
    \begin{axis}[
        height = 0.30\textwidth,
        width  = 1.0\textwidth,
        ymode=log,
        ylabel={Time per call [ms]},
        ymin=1.5e-2, ymax=3e1,
        ytick={1e-2, 1e-1, 1e0, 1e1},
        xmin=-0.7, xmax=9.7,
        xtick={0,1,2,3,4,5,6,7,8,9},
        xticklabels={
            {$\bc$ solve},
            {$\cJ_{physics}, \Omega_1$: $NN_1$},
            {$\cJ_{physics}, \Omega_1$: $g_1$},
            {$\cJ_{physics}, \Omega_2$: $NN_2$},
            {$\cJ_{physics}, \Omega_2$: $g_2$},
            {$\cJ_{dbc}$},
            {$\cJ_{int}$},
            {$\cJ_{fint}$},
            {$\cJ$},
            {$\nabla_\theta\cJ$},
        },
        xticklabel style={rotate=40, anchor=east, font=\scriptsize},
        tick scale binop ={\times},
        every axis plot/.append style={
            only marks,
            mark size=1.6pt,
            line width=0.8,
            error bars/y dir=both,
            error bars/y explicit,
            error bars/error bar style={line width=0.7, solid},
            error bars/error mark options={line width=0.7, mark size=2.5pt},
        },
        legend style={
            draw=none, fill=none,
            at={(rel axis cs: 0.02, 0.97)},
            anchor=north west,
            legend cell align={left},
            legend columns=2,
            /tikz/every even column/.append style={column sep=0.4cm},
        },
    ]
\pgfplotsset{set layers=standard}%

        \addplot [
            mark=*,
            blue,
            mark options={fill=blue, draw=blue},
        ]
        table [
            x expr=\thisrowno{0}-0.13, y index=1,
            y error minus index=2, y error plus index=3,
            skip first n=1,
        ]{data_petal_timing_cmp_buffer.txt};
        \addlegendentry{Buffer}

        \addplot [
            mark=square*,
            red,
            mark options={fill=red, draw=red},
        ]
        table [
            x expr=\thisrowno{0}+0.13, y index=1,
            y error minus index=2, y error plus index=3,
            skip first n=1,
        ]{data_petal_timing_cmp_mpinn.txt};
        \addlegendentry{M-PINN}

    \end{axis}
\end{tikzpicture}
%
    \caption{\refB{Per-call runtime of the buffer and M-PINN components on the
    petal problem, for $N_{bnd}=N_{itf}=10$.
    Markers are medians and bars span the minimum and maximum over profiled calls.}}
    \label{fig:petal_timing_comparison}
\end{figure}
\refB{
Figure~\ref{fig:petal_timing_comparison} compares the buffer approach against an M-PINN
model with identical network sizes. Both methods evaluate the physics residual
$\cJ_{physics}$ through a neural network in each subdomain, at essentially the same cost
($0.17$ and $0.65$~ms for $\Omega_1$ and $\Omega_2$ with the buffer, against $0.14$ and
$0.56$~ms for M-PINN). Relative to this shared work, the buffer approach carries two
additional costs: the coefficient solve, and the evaluation of the buffer term $g_m$
within the residual. The coefficient solve is inexpensive, at $0.14$~ms, because the
matrix in (\ref{eq:g-bc-high-dim}) is independent of $\theta$ for a fixed geometry, sample
set, and radius, and is therefore LU-factorized once before training and reused at every
iteration; only the right-hand side is rebuilt. Evaluating $g_m$ costs about the same as
the network term in $\Omega_1$ ($0.17$~ms) and somewhat more in $\Omega_2$ ($1.61$~ms),
where more collocation points are located. M-PINN, in turn, must evaluate the
boundary and interface penalty terms $\cJ_{dbc}$, $\cJ_{int}$ and $\cJ_{fint}$, which the
buffer approach does not, although these are individually cheap ($\sim0.04$~ms each).
The net effect is that a full loss evaluation $\cJ$ takes $3.6\times$ longer than its
M-PINN counterpart, $2.56$~ms against $0.71$~ms.
The overhead of the gradient $\nabla_\theta\cJ$, however, is considerably smaller:
$4.32$~ms against $3.14$~ms, a factor of only $1.4\times$. This is presumably because the
buffer contributes only a small number of coefficients to backpropagate through, far fewer
than the network parameters that dominate the reverse pass in both methods, so the
additional forward work is largely amortized in the gradient computation.
Note that the individually profiled subcomponents do not add up to the total cost, since
JAX's just-in-time compilation fuses the terms and eliminates redundant work across them. The
component timings should therefore be read as indicative of relative cost rather than as
an exact decomposition of $\cJ$.
}
\par
\begin{figure}[tbh]
    \begin{tikzpicture}[
font=\small,
]
    \begin{groupplot}[
        group style={
            group name = my plots,
            group size= 2 by 2,
            xlabels at =edge bottom,
            ylabels at =edge left,
            horizontal sep=1.2cm,
            vertical sep=2.6cm,
        },
        height = 0.30\textwidth,
        width  = 0.45\textwidth,
        xmode=log,
        ymode=log,
        xmin=5, xmax=58,
        ymin=6e-2, ymax=1.5e1,
        xtick={6,10,12,24,48},
        xticklabels={6,10,12,24,48},
        minor xtick={},
        xticklabel style={font=\scriptsize},
        ytick={1e-1, 1e0, 1e1},
        ylabel={Time per call [ms]},
        tick scale binop ={\times},
        every axis plot/.append style={
            mark=*,
            mark size=1.4pt,
            line width=0.9,
            error bars/y dir=both,
            error bars/y explicit,
            error bars/error bar style={line width=0.7, solid},
            error bars/error mark options={line width=0.7, mark size=2.5pt},
        },
        legend style={
            draw=none, fill=none,
            at={(rel axis cs: 0.0, 1.04)},
            anchor=south west,
            legend cell align={left},
            nodes={scale=0.9},
            /tikz/every even column/.append style={column sep=0.35cm},
        },
        name=chung,
    ]
\pgfplotsset{set layers=standard}%

        \nextgroupplot[
            xlabel={$N_{bnd}$},
            legend columns=5,
        ]

            \addplot [blue, mark options={fill=blue, draw=blue}]
            table [x index=0, y index=1,
                   y error minus index=2, y error plus index=3,
                   skip first n=1,
            ]{data_petal_timing_vs_n_bnd_total.txt};
            \addlegendentry{$\bc$ solve}

            \addplot [green!60!black, mark options={fill=green!60!black, draw=green!60!black}]
            table [x index=0, y index=4,
                   y error minus index=5, y error plus index=6,
                   skip first n=1,
            ]{data_petal_timing_vs_n_bnd_total.txt};
            \addlegendentry{$\cJ_{physics}, \Omega_1$}

            \addplot [orange, mark options={fill=orange, draw=orange}]
            table [x index=0, y index=7,
                   y error minus index=8, y error plus index=9,
                   skip first n=1,
            ]{data_petal_timing_vs_n_bnd_total.txt};
            \addlegendentry{$\cJ_{physics}, \Omega_2$}

            \addplot [violet, mark options={fill=violet, draw=violet}]
            table [x index=0, y index=10,
                   y error minus index=11, y error plus index=12,
                   skip first n=1,
            ]{data_petal_timing_vs_n_bnd_total.txt};
            \addlegendentry{$\cJ$}

            \addplot [red, mark options={fill=red, draw=red}]
            table [x index=0, y index=13,
                   y error minus index=14, y error plus index=15,
                   skip first n=1,
            ]{data_petal_timing_vs_n_bnd_total.txt};
            \addlegendentry{$\nabla_\theta\cJ$}

        \nextgroupplot[
            xlabel={$N_{itf}$},
        ]

            \addplot [blue, mark options={fill=blue, draw=blue}]
            table [x index=0, y index=1,
                   y error minus index=2, y error plus index=3,
                   skip first n=1,
            ]{data_petal_timing_vs_n_itf_total.txt};

            \addplot [green!60!black, mark options={fill=green!60!black, draw=green!60!black}]
            table [x index=0, y index=4,
                   y error minus index=5, y error plus index=6,
                   skip first n=1,
            ]{data_petal_timing_vs_n_itf_total.txt};

            \addplot [orange, mark options={fill=orange, draw=orange}]
            table [x index=0, y index=7,
                   y error minus index=8, y error plus index=9,
                   skip first n=1,
            ]{data_petal_timing_vs_n_itf_total.txt};

            \addplot [violet, mark options={fill=violet, draw=violet}]
            table [x index=0, y index=10,
                   y error minus index=11, y error plus index=12,
                   skip first n=1,
            ]{data_petal_timing_vs_n_itf_total.txt};

            \addplot [red, mark options={fill=red, draw=red}]
            table [x index=0, y index=13,
                   y error minus index=14, y error plus index=15,
                   skip first n=1,
            ]{data_petal_timing_vs_n_itf_total.txt};

        \nextgroupplot[
            xlabel={$N_{bnd}$},
            legend columns=2,
        ]

            \addplot [green!60!black, mark options={fill=green!60!black, draw=green!60!black}]
            table [x index=0, y index=1,
                   y error minus index=2, y error plus index=3,
                   skip first n=1,
            ]{data_petal_timing_vs_n_bnd_split.txt};
            \addlegendentry{$\cJ_{physics}, \Omega_1$: $NN_1$}

            \addplot [green!60!black, dashed, mark options={solid, fill=white, draw=green!60!black}]
            table [x index=0, y index=4,
                   y error minus index=5, y error plus index=6,
                   skip first n=1,
            ]{data_petal_timing_vs_n_bnd_split.txt};
            \addlegendentry{$\cJ_{physics}, \Omega_1$: $g_1$}

            \addplot [orange, mark options={fill=orange, draw=orange}]
            table [x index=0, y index=7,
                   y error minus index=8, y error plus index=9,
                   skip first n=1,
            ]{data_petal_timing_vs_n_bnd_split.txt};
            \addlegendentry{$\cJ_{physics}, \Omega_2$: $NN_2$}

            \addplot [orange, dashed, mark options={solid, fill=white, draw=orange}]
            table [x index=0, y index=10,
                   y error minus index=11, y error plus index=12,
                   skip first n=1,
            ]{data_petal_timing_vs_n_bnd_split.txt};
            \addlegendentry{$\cJ_{physics}, \Omega_2$: $g_2$}

        \nextgroupplot[
            xlabel={$N_{itf}$},
        ]

            \addplot [green!60!black, mark options={fill=green!60!black, draw=green!60!black}]
            table [x index=0, y index=1,
                   y error minus index=2, y error plus index=3,
                   skip first n=1,
            ]{data_petal_timing_vs_n_itf_split.txt};

            \addplot [green!60!black, dashed, mark options={solid, fill=white, draw=green!60!black}]
            table [x index=0, y index=4,
                   y error minus index=5, y error plus index=6,
                   skip first n=1,
            ]{data_petal_timing_vs_n_itf_split.txt};

            \addplot [orange, mark options={fill=orange, draw=orange}]
            table [x index=0, y index=7,
                   y error minus index=8, y error plus index=9,
                   skip first n=1,
            ]{data_petal_timing_vs_n_itf_split.txt};

            \addplot [orange, dashed, mark options={solid, fill=white, draw=orange}]
            table [x index=0, y index=10,
                   y error minus index=11, y error plus index=12,
                   skip first n=1,
            ]{data_petal_timing_vs_n_itf_split.txt};

  \end{groupplot}
\node[below = 1.2cm of my plots c1r1.south west,
    anchor=west, xshift=-5mm,
] {(a) Total, vs.\ $N_{bnd}$};
\node[below = 1.2cm of my plots c2r1.south west,
    anchor=west, xshift=-5mm,
] {(b) Total, vs.\ $N_{itf}$};
\node[below = 1.2cm of my plots c1r2.south west,
    anchor=west, xshift=-5mm,
] {(c) Split, vs.\ $N_{bnd}$};
\node[below = 1.2cm of my plots c2r2.south west,
    anchor=west, xshift=-5mm,
] {(d) Split, vs.\ $N_{itf}$};
\end{tikzpicture}
%
    \caption{\refB{Per-call runtime of the buffer components against the numbers of
    boundary and interface sample points, $N_{bnd}$ and $N_{itf}$.
    (a, b) total cost per component; (c, d) the neural and buffer parts of the
    subdomain residuals.
    Markers are medians and bars span the minimum and maximum over profiled calls.}}
    \label{fig:petal_timing_vs_counts}
\end{figure}
\refB{
Figure~\ref{fig:petal_timing_vs_counts} shows how these costs scale with the number of
sample points, which also sets the number of buffer degrees of freedom. The growth is
markedly sublinear in every component. The buffer part of the subdomain residuals scales
as $\cO(N_{bnd}^{0.05})$ and $\cO(N_{bnd}^{0.40})$ for $\Omega_1$ and $\Omega_2$, and as
$\cO(N_{itf}^{0.32})$ and $\cO(N_{itf}^{0.10})$ with respect to the interface points. The
full loss and its gradient scale as $\cO(N_{bnd}^{0.31})$ and $\cO(N_{bnd}^{0.27})$, and as
$\cO(N_{itf}^{0.08})$ and $\cO(N_{itf}^{-0.01})$, so an eightfold increase in the number of
sample points, from $6$ to $48$, raises the cost of a gradient evaluation by well under a
factor of two.
The coefficient solve likewise remains inexpensive across the whole range, staying below
$0.3$~ms with nearly flat scaling.
While the triangular solve is asymptotically $\cO(N^2)$, at these system sizes it amounts
to only a few thousand floating-point operations and is dominated by kernel-launch
overhead, so the measured cost stays nearly flat.
It should be noted that a larger $N_{bnd}$ or $N_{itf}$ generally requires smaller radii for each radial
basis function in order to keep the coefficient system well conditioned.
This in turn calls for more collocation points to resolve the losses near the boundary and interface
and thus lengthens the loss evaluation, though this coupled effect has not been investigated here.
Nonetheless, the slow scaling of the overhead itself indicates that the buffer
approach remains practical as constraints are refined, and suggests that it should extend
to geometries with many interfaces or to three dimensions, where the number of sample
points required grows substantially.
}
\section{Conclusion}\label{sec:conclusion}
We presented two hard-constrained PINN methodologies for interface problems: a \textbf{windowing approach} that embeds interface continuity and flux balance exactly in the trial space, and a \textbf{buffer approach} that enforces the jump physics through auxiliary buffer functions in a structurally consistent, \refB{discretely} hard-constrained manner. By decoupling interface enforcement from the PDE residual minimization, both designs eliminate the loss-balancing pathology, leading to more stable training and sharper resolution near discontinuities.

We compared these approaches against soft-constrained counterparts---interface PINNs, adaptive interface PINNs, $\phi$-PINNs and \refB{M-PINN}. \refB{A statistical study shows that the buffer approach attains the smallest median error on all three one-dimensional benchmarks.
The windowing approach achieves the smallest minimum errors in occasional runs but does not consistently improve on the soft-constrained baselines in typical performance. The two-dimensional benchmarks show that the buffer approach} can maintain \refB{better} accuracy across different training configurations, while avoiding the need for loss-weight tuning.

A comparison between the windowing approach and the buffer approach reveals that,
although product form with window functions may strictly enforce boundary and interface conditions,
it can also limit the expressive power of the neural networks and make the training overly constrained.
The training performance of the windowing approach is inherently tied to the shape of the window function derivatives, i.e. how much they span the source term or other equation contributions.
In higher dimensions, the windowing approach faces additional challenges at domain corners where boundary and interface windows overlap.
The buffer approach, by contrast, leaves the neural networks unrestricted by any window function,
enabling faster convergence and more robust performance across problem settings.
On the two-dimensional benchmark, the buffer approach achieves a lower error than the windowing approach while avoiding the corner difficulties entirely.
\par
While the buffer approach achieved lower errors, the improvement was
marginal, as the soft-constrained methods also performed comparably well.
Although not explored in this study,
we expect the advantage of the proposed approach to become more pronounced under more extreme conditions.
For example, much larger contrasts in the diffusivities,
or stronger inhomogeneities in the boundary and interface conditions,
would demand more careful loss-weight tuning for the soft-constrained approaches.


\refA{The present study is confined to scalar second-order elliptic interface problems, for which each boundary or interface condition constrains a single scalar quantity. Systems in which the solution components are coupled through the constraints---elasticity being the canonical example, where a zero-traction condition couples all components along a boundary or interface---require the buffer degrees of freedom to be determined from a coupled, vector-valued constraint system rather than the componentwise one used here. Such higher-order and coupled formulations are beyond the scope of this work and are left as a direction for future study.}

In future work, the same constraint-enforcement principles can be extended to more complex settings, including two-dimensional unsteady and three-dimensional elliptic problems, multi-physics couplings, and time-dependent interface dynamics. \refB{A systematic investigation of alternative buffer function forms is also warranted, since the present buffer construction only marginally outperforms the strongest soft-constrained baselines on some benchmarks; identifying buffer parametrizations that yield a larger and more consistent advantage over soft-constrained counterparts is left for future work.} Also, while the buffer formulation was implemented here in a domain-decomposed form, its underlying correction mechanism is not tied to this particular architecture. A natural extension is to incorporate analogous buffer corrections into other interface-aware parametrizations, such as I-PINN or $\phi$-PINN architectures. Additional promising directions include adaptive window or buffer construction, \textit{a posteriori} error indicators, and multi-fidelity or data-assimilative training. Another natural direction is to investigate whether the buffer formulation can be combined with operator-learning architectures for parametric interface problems. Together, these directions point toward broader and more robust constraint-embedded learning frameworks in which additive constructions—ranging from analytic lifting to local corrective strategies—offer a flexible alternative to traditional soft enforcement or multiplicative masking approaches.

\section*{Acknowledgement}

The research described in this paper was conducted under the auspices of Lawrence Livermore National Laboratory (LLNL, Contract No. DE-AC52-07NA27344), Pacific Northwest National Laboratory (PNNL, Contract
No. DE-AC05–76RL01830), and Sandia National Laboratories (SNL, Contract
No.DE-NA0003525) as part of the DOE Interlaboratory (IL) Laboratory Directed Research and Development (LDRD) Pilot Program. 
The LLNL LDRD project number is 25-ERD-052 and the release number is
LLNL-JRNL-2010925.

\appendix

\section{Analytic solution for Problem 2}\label{app:prob2-sol}
The problem 2 in Section~\ref{subsec:3interface} has an analytic solution of the form,
\begin{equation}
u(x) =
\begin{cases}
- \frac{x^2}{2k_1} + \frac{c_1 x}{Kk_1} & x \in [0, 0.25] \\
- \frac{x^2}{2k_2} + \frac{c_2 x}{Kk_2} + \frac{c_3}{Kk_2} & x \in [0.25, 0.5] \\
- \frac{x^2}{2k_3} + \frac{c_4 x}{Kk_3} + \frac{c_5}{Kk_3} & x \in [0.5, 0.75] \\
- \frac{x^2}{2k_4} + \frac{c_6 x}{Kk_4} + \frac{c_7}{Kk_4} & x \in [0.75, 1],
\end{cases}
\end{equation}
with
$K = k_{1} k_{2} k_{3} + k_{1} k_{2} k_{4} + k_{1} k_{3} k_{4}
+ k_{2} k_{3} k_{4}$.
The coefficients are given as
\begin{subequations}
\begin{equation}
\begin{split}
c_1 = c_2 &= c_4 = c_6 =
\frac{7}{8} k_{1} k_{2} k_{3}
+ \frac{5}{8} k_{1} k_{2} k_{4}
+ \frac{3}{8} k_{1} k_{3} k_{4}
+ \frac{1}{8} k_{2} k_{3} k_{4}
\end{split}
\end{equation}
\begin{equation}
\begin{split}
c_3 =
&- \frac{3}{16} k_{1} k_{2} k_{3}
- \frac{1}{8} k_{1} k_{2} k_{4}
- \frac{1}{16} k_{1} k_{3} k_{4}
+ \frac{3}{16} k_{2}^{2} k_{3}
+ \frac{1}{8} k_{2}^{2} k_{4}
+ \frac{1}{16} k_{2} k_{3} k_{4}
\end{split}
\end{equation}
\begin{equation}
\begin{split}
c_5 =
&- \frac{5}{16} k_{1} k_{2} k_{3}
- \frac{3}{16} k_{1} k_{2} k_{4}
+ \frac{1}{8} k_{1} k_{3}^{2}
+ \frac{3}{16} k_{2} k_{3}^{2}
+ \frac{3}{16} k_{2} k_{3} k_{4}
\end{split}
\end{equation}
\begin{equation}
\begin{split}
c_7 =
&- \frac{3}{8} k_{1} k_{2} k_{3}
- \frac{1}{8} k_{1} k_{2} k_{4}
+ \frac{1}{8} k_{1} k_{3} k_{4}
+ \frac{3}{8} k_{2} k_{3} k_{4}.
\end{split}
\end{equation}
\end{subequations}

\section{A preliminary study on optimal window configuration}\label{app:window-prelim}

\subsection{Effect of window functions in windowing approach}\label{subsec:comparison-window-order}
First, we investigate the impact of the choice of window function shape on the solution accuracy. While there is no unique window function that satisfies the condition \eqref{eq:wint-bc}
or \eqref{eq:wn-bc},
not all windows are optimal for training convergence.
\par
\begin{figure}[tbph]
    \begin{tikzpicture}[
]
    \begin{groupplot}[
        group style={
            group name = my plots,
            group size= 2 by 1,
            xlabels at =edge bottom,
            horizontal sep=2cm,
            vertical sep=2.2cm,
        },
        height = 0.45\textwidth,
        width = 0.5\textwidth,
        name=chung,
    ]    
\pgfplotsset{set layers=standard}%

        \nextgroupplot[
            xlabel={$x$},
            ylabel={$W(x)$},
            tick scale binop ={\times},
            xmin=0., xmax=1.,
            ymin=0., ymax=1.1,
            xtick={0, 0.25, 0.5, 0.75, 1.},
            legend style={
                draw=none, fill=none,
                at={(rel axis cs: -0.2, 1.0)},
                anchor=south west,
                nodes={scale=1.0},
                legend cell align={left},
                legend columns=3,
                /tikz/every even column/.append style={column sep=0.5cm},
            },
        ]

            \addlegendimage{
                line width=1.0,
                solid,
                mark=none,
                green,};
            \addlegendentry{$W_m(\bx)$, interior};
            \addlegendimage{
                line width=1.0,
                dashed,
                mark=none,
                blue,};
            \addlegendentry{$W_{b,d}(\bx)$, boundary};
            \addlegendimage{
                line width=1.0,
                dashed,
                mark=none,
                red,};
            \addlegendentry{$W_{b,n}(\bx)$, boundary};
            \addlegendimage{
                line width=1.0,
                dashed,
                mark=none,
                brown,};
            \addlegendentry{$W_{ij,d}(\bx)$, interface};
            \addlegendimage{
                line width=1.0,
                dashed,
                mark=none,
                purple,};
            \addlegendentry{$W_{ij,n}(\bx)$, interface};

            \addplot+[
                mark=none,
                gray,
                line width=1.5,
            ] coordinates {(0.5, 0) (0.5, 2)};

            \addplot+ [
                line width=1.0,
                solid,
                smooth,
                mark=none,
                green,
            ]
            table [
                x index=0, y index=1,
            ]{data_config1_window_function_int0L.txt};
            \addplot+ [
                line width=1.0,
                solid,
                smooth,
                mark=none,
                green,
            ]
            table [
                x index=0, y index=1,
            ]{data_config1_window_function_int0M.txt};
            \addplot+ [
                line width=1.0,
                solid,
                smooth,
                mark=none,
                green,
            ]
            table [
                x index=0, y index=1,
            ]{data_config1_window_function_int0R.txt};
            \addplot+ [
                line width=1.0,
                solid,
                smooth,
                mark=none,
                green,
            ]
            table [
                x index=0, y index=1,
            ]{data_config1_window_function_int1L.txt};
            \addplot+ [
                line width=1.0,
                solid,
                smooth,
                mark=none,
                green,
            ]
            table [
                x index=0, y index=1,
            ]{data_config1_window_function_int1M.txt};
            \addplot+ [
                line width=1.0,
                solid,
                smooth,
                mark=none,
                green,
            ]
            table [
                x index=0, y index=1,
            ]{data_config1_window_function_int1R.txt};

            \addplot+ [
                line width=1.0,
                dashed,
                smooth,
                mark=none,
                red,
            ]
            table [
                x index=0, y index=1,
            ]{data_config1_window_function_bndL.txt};
            \addplot+ [
                line width=1.0,
                dashed,
                smooth,
                mark=none,
                red,
            ]
            table [
                x index=0, y index=1,
            ]{data_config1_window_function_bndR.txt};
            \addplot+ [
                line width=1.0,
                dashed,
                smooth,
                mark=none,
                purple,
            ]
            table [
                x index=0, y index=1,
            ]{data_config1_window_function_itf.txt};

            \addplot+ [
                line width=1.0,
                dashed,
                smooth,
                mark=none,
                blue,
            ]
            table [
                x index=0, y index=2,
            ]{data_config1_window_function_bndL.txt};
            \addplot+ [
                line width=1.0,
                dashed,
                smooth,
                mark=none,
                blue,
            ]
            table [
                x index=0, y index=2,
            ]{data_config1_window_function_bndR.txt};
            \addplot+ [
                line width=1.0,
                dashed,
                smooth,
                mark=none,
                brown,
            ]
            table [
                x index=0, y index=2,
            ]{data_config1_window_function_itf.txt};

        \nextgroupplot[
            tick scale binop ={\times},
            xmin = 0, xmax = 3,
            ymin = 0, ymax = 3,
            xtick={0.5, 1.5, 2.5},
            xticklabels={{$\tW^{(1)}_{d}$, $\tW^{(1)}_{n}$}, {$\tW^{(2)}_{d}$, $\tW^{(2)}_{n}$}, {$\tW^{(3)}_{d}$, $\tW^{(3)}_{n}$}},
            x tick label style={xshift=20pt, rotate=45,anchor=east},
            ytick={0.5, 1.5, 2.5},
            yticklabels={$\tW^{(3)}_{int}$, $\tW^{(2)}_{int}$, $\tW^{(1)}_{int}$},
        ]
        
            \edef\imagepath{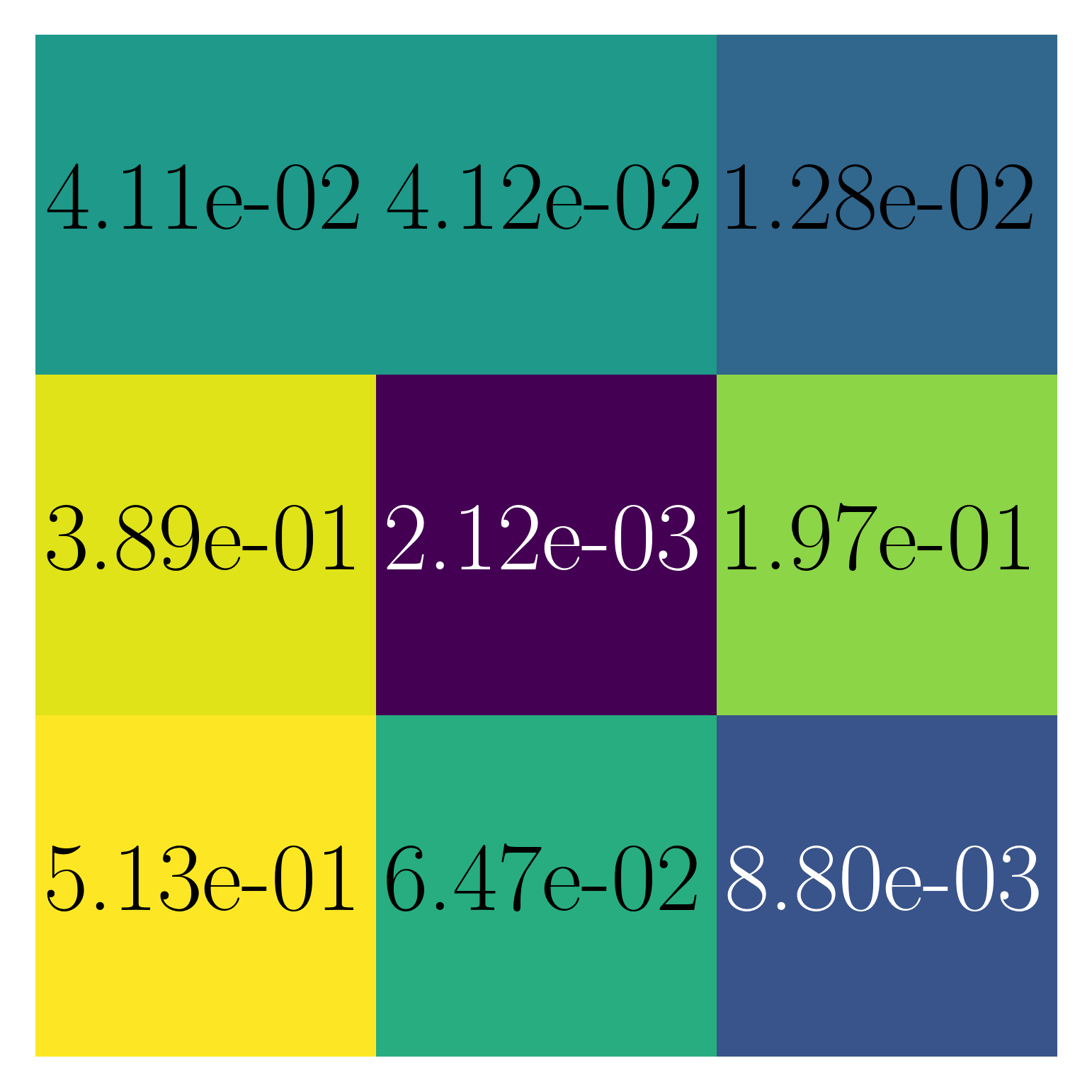}
            \addplot graphics[xmin=-0.1, xmax=3.1, ymin=-0.1, ymax=3.1]{\imagepath};

  \end{groupplot}
\node[below = 2cm of my plots c1r1.south west,
    anchor=west,
] {(a) Window configuration};
\node[below = 2cm of my plots c2r1.south west,
    anchor=west,
] {(b) Relative $L_2$ error};
\end{tikzpicture}
%
    \caption{Training performance of the windowing approach with various window functions:
    (a) window configuration for the target problem in Section~\ref{subsec:1interface}; and
    (b) relative $L_2$ error of the trained solutions depending on the combination of different window functions.}
    \label{fig:config1}
\end{figure}
Problem 1 in Section~\ref{subsec:1interface} is used to compare various window functions.
Figure~\ref{fig:config1}(a) shows the base window configuration used for the problem.
$M=2$ interior subdomains are used to represent the solutions for distinct diffusivities,
with their centers at $\bx_1=0.25$ and $\bx_2=0.75$, respectively.
Two boundary subdomains are located at both ends of the domain $\bx=0$ and $\bx=1$,
and one interface subdomain at the interface $\bx=0.5$.
We used a uniform subdomain size $\Delta\bx_m=0.25$ for all subdomains.
\par
The solution ansatz \eqref{eq:window-ansatz} is trained
with target window functions listed in Table~\ref{tab:wint}-\ref{tab:wn}.
Since $\tW_d$ and $\tW_n$ are used as pairs at the boundary and interface,
their polynomial orders are set to the same in all cases.
A fully connected network with hidden layers $[25, 25]$ is
used for each subdomain interior.
For all cases, the training loss \eqref{eq:J-int} is evaluated
with $60$ collocation points uniformly distributed in the domain.
The training is performed via SOAP optimizer~\cite{vyas2024soap},
with learning rate $5\times10^{-3}$ for $10^4$ iterations.
\par
Figure~\ref{fig:config1}~(b) shows the relative error $L_2$ of the trained solution
 relative to the exact solution.
Although the shapes of the functions in Figure~\ref{fig:windows} vary only moderately, the resulting precision changes drastically, spanning nearly three orders of magnitude.
With the choice of $\tW_d^{(1)}$ and $\tW_n^{(1)}$ i.e. cubic polynomials,
the solution error remains higher than $1\%$ regardless of $\tW_{int}$.
The highest accuracy is achieved for the case of $\tW_d^{(2)}$, $\tW_n^{(2)}$,
and $\tW_{int}^{(2)}$, whose second-order derivatives vanish at the subdomain centers.
Interestingly, the choice of smoother window functions e.g., fifth order polynomials, does not necessarily lead to improved solution accuracy.
This result shows that, while many window functions can hard-constrain the boundary and interface conditions, the specific choice of window functions can impact the training performance and solution accuracy. Furthermore, a smoother window function is not necessarily the best-performing.
\par
\begin{figure}[tbph]
    \begin{tikzpicture}[
]
    \begin{groupplot}[
        group style={
            group name = my plots,
            group size= 3 by 1,
            xlabels at =edge bottom,
            horizontal sep=2cm,
            vertical sep=2.2cm,
        },
        height = 0.35\textwidth,
        width = 0.35\textwidth,
        enlarge x limits={false, abs value = 5mm},
        enlarge y limits={false, abs value = 5mm},
        name=chung,
    ]    
\pgfplotsset{set layers=standard}%

        \nextgroupplot[
            xlabel={$\tau$},
            ylabel={$\Dparttwo{}{\tau}\tW^{(k)}_{int}$},
            tick scale binop ={\times},
            xmin=-1, xmax=1,
            xtick={-1, -0.5, 0, 0.5, 1},
            xticklabels={1, 0.5, 0, 0.5, 1},
            legend style={
                draw=none, fill=none,
                at={(rel axis cs: 0, 1.0)},
                anchor=south west,
                nodes={scale=1.0},
                legend cell align={left},
                legend columns=3,
                /tikz/every even column/.append style={column sep=0.5cm},
            },
        ]

            \addplot+ [
                line width=1.0,
                solid,
                mark=none,
                blue,
            ]
            table [
                x index=0, y index=4,
            ]{data_window_derivatives_w_derivatives.int.txt};

            \addplot+ [
                line width=1.0,
                solid,
                mark=none,
                brown,
            ]
            table [
                x index=0, y index=5,
            ]{data_window_derivatives_w_derivatives.int.txt};

            \addplot+ [
                line width=1.0,
                solid,
                mark=none,
                red,
            ]
            table [
                x index=0, y index=6,
            ]{data_window_derivatives_w_derivatives.int.txt};

            \legend{$k=1$, $k=2$, $k=3$};

        \nextgroupplot[
            xlabel={$\tau$},
            ylabel={$\Dparttwo{}{\tau}\tW^{(k)}_{d}$},
            tick scale binop ={\times},
        ]

            \addplot+ [
                line width=1.0,
                solid,
                mark=none,
                blue,
            ]
            table [
                x index=0, y index=4,
            ]{data_window_derivatives_w_derivatives.dir.txt};

            \addplot+ [
                line width=1.0,
                solid,
                mark=none,
                brown,
            ]
            table [
                x index=0, y index=5,
            ]{data_window_derivatives_w_derivatives.dir.txt};

            \addplot+ [
                line width=1.0,
                solid,
                mark=none,
                red,
            ]
            table [
                x index=0, y index=6,
            ]{data_window_derivatives_w_derivatives.dir.txt};

        \nextgroupplot[
            xmin=0., xmax=1.,
            xlabel={$\tau$},
            ylabel={$\Dparttwo{}{\tau}\tW^{(k)}_{n}$},
            tick scale binop ={\times},
        ]

            \addplot+ [
                line width=1.0,
                solid,
                mark=none,
                blue,
            ]
            table [
                x index=0, y index=4,
            ]{data_window_derivatives_w_derivatives.neu.txt};

            \addplot+ [
                line width=1.0,
                solid,
                mark=none,
                brown,
            ]
            table [
                x index=0, y index=5,
            ]{data_window_derivatives_w_derivatives.neu.txt};

            \addplot+ [
                line width=1.0,
                solid,
                mark=none,
                red,
            ]
            table [
                x index=0, y index=6,
            ]{data_window_derivatives_w_derivatives.neu.txt};

  \end{groupplot}
\node[below = 1.5cm of my plots c1r1.south west,
    anchor=west,
] {(a)};
\node[below = 1.5cm of my plots c2r1.south west,
    anchor=west,
] {(b)};
\node[below = 1.5cm of my plots c3r1.south west,
    anchor=west,
] {(c)};
\end{tikzpicture}
%
    \caption{Second derivatives of target window functions:
    (a) interior;
    (b) Dirichlet; and
    (c) Neumann.
    For the interior window, the coordinate was mirrored about $\tau=0$ to
    illustrate the behavior of the functions across the subdomain.}
    \label{fig:window-derivatives}
\end{figure}
The impact of the window function shape on the training performance
can be explained with the shape of their higher derivatives,
specifically up to second derivatives in this study.
This stems from the fact that PINN is trained
with the residual loss for the physics equation,
which is presented mainly as a differential equation.
For example, substituting an interior subdomain solution $u^{(S)}_m$ in \eqref{eq:uD-m}
into the left-hand side of the physics equation \eqref{eq:poisson} yields
\begin{equation}
\LHS[u^{(S)}_m]
=
\refB{-\kappa_m W_m \nabla^2NN_m
-2\kappa_m \nabla W_m\cdot\nabla NN_m
-\kappa_m NN_m \nabla^2 W_m},
\end{equation}
where $\nabla W_m$ and $\nabla^2W_m$ have contributions to resolve the physics equation.
While the function shape itself may not vary significantly with the polynomial order, as shown in Figure~\ref{fig:windows}, the shape of its higher derivatives do.
Figure~\ref{fig:window-derivatives}
shows the second derivatives of target window functions in Table~\ref{tab:wint}-\ref{tab:wn},
where the derivative shape changes drastically with the polynomial order.
Note that the second derivatives of $\tW^{(k)}_d$ and $\tW^{(k)}_n$
do not vanish for $k=1$ at $\tau=1$, where a boundary/interface subdomain ends and faces another boundary/interface subdomain.
Such non-vanishing second derivatives can introduce discontinuities in the residual for the physics equation, which is challenging for a NN to resolve.
\begin{figure}[tbph]
    \begin{tikzpicture}[
]
    \begin{groupplot}[
        group style={
            group name = my plots,
            group size= 3 by 1,
            xlabels at =edge bottom,
            horizontal sep=2cm,
            vertical sep=2.2cm,
        },
        height = 0.45\textwidth,
        width = 0.5\textwidth,
        enlarge x limits={false, abs value = 5mm},
        enlarge y limits={false, abs value = 5mm},
        xtick={0, 0.25, 0.5, 0.75, 1.},
        name=chung,
    ]    
\pgfplotsset{set layers=standard}%

        \nextgroupplot[
            ymin=-20, ymax=25,
            xlabel={$x$},
            ylabel={$\LHS[u]$},
            tick scale binop ={\times},
            legend style={
                draw=none, fill=none,
                at={(rel axis cs: -0.4, 1.05)},
                anchor=south west,
                nodes={scale=1.0},
                legend cell align={left},
                legend columns=4,
                /tikz/every even column/.append style={column sep=0.5cm},
            },
        ]

            \addlegendimage{
                line width=1.0,
                solid,
                mark=none,
                blue,
            };
            \addlegendentry{$\LHS[u^{(S)}_m]$};
            \addlegendimage{
                line width=1.0,
                solid,
                mark=none,
                green,
            };
            \addlegendentry{$\LHS[u^{(B)}_b]$};
            \addlegendimage{
                line width=1.0,
                solid,
                mark=none,
                red,
            };
            \addlegendentry{$\LHS[u^{(I)}_{ij}]$};
            \addlegendimage{
                line width=1.0,
                dashed,
                mark=none,
                black,
            };
            \addlegendentry{$\LHS[u_{NN}]$};

            \addplot+[
                mark=none,
                gray!50,
                line width=1.5,
            ] coordinates {(0.5, -30) (0.5, 30)};
            \addplot+[
                mark=none,
                gray!50,
                line width=1.5,
            ] coordinates {(0.25, -30) (0.25, 30)};
            \addplot+[
                mark=none,
                gray!50,
                line width=1.5,
            ] coordinates {(0.75, -30) (0.75, 30)};

            \addplot+ [
                line width=1.0,
                solid,
                mark=none,
                green,
            ]
            table [
                x index=0, y expr=-\thisrowno{1},
            ]{data_config1_window1_rhs_BCL.txt};

            \addplot+ [
                line width=1.0,
                solid,
                mark=none,
                green,
            ]
            table [
                x index=0, y expr=-\thisrowno{1},
            ]{data_config1_window1_rhs_BCR.txt};

            \addplot+ [
                line width=1.0,
                solid,
                mark=none,
            ]
            table [
                x index=0, y expr=-\thisrowno{1},
            ]{data_config1_window1_rhs_ITF.txt};

            \addplot+ [
                line width=1.0,
                solid,
                smooth,
                mark=none,
                blue,
            ]
            table [
                x index=0, y expr=-\thisrowno{1},
            ]{data_config1_window1_rhs_pinn0.txt};

            \addplot+ [
                line width=1.0,
                solid,
                smooth,
                mark=none,
                blue,
            ]
            table [
                x index=0, y expr=-\thisrowno{1},
            ]{data_config1_window1_rhs_pinn1.txt};

            \addplot+ [
                line width=1.0,
                dashed,
                mark=none,
                black,
            ]
            table [
                x index=0, y expr=-\thisrowno{1},
            ]{data_config1_window1_rhs_total.txt};

        \nextgroupplot[
            ymin=-40, ymax=40,
            xlabel={$x$},
            ylabel={$\LHS[u]$},
            tick scale binop ={\times},
        ]

            \addplot+[
                mark=none,
                gray!50,
                line width=1.5,
            ] coordinates {(0.5, -50) (0.5, 50)};
            \addplot+[
                mark=none,
                gray!50,
                line width=1.5,
            ] coordinates {(0.25, -50) (0.25, 50)};
            \addplot+[
                mark=none,
                gray!50,
                line width=1.5,
            ] coordinates {(0.75, -50) (0.75, 50)};

            \addplot+ [
                line width=1.0,
                solid,
                mark=none,
                green,
            ]
            table [
                x index=0, y expr=-\thisrowno{1},
            ]{data_config1_window5_rhs_BCL.txt};

            \addplot+ [
                line width=1.0,
                solid,
                mark=none,
                green,
            ]
            table [
                x index=0, y expr=-\thisrowno{1},
            ]{data_config1_window5_rhs_BCR.txt};

            \addplot+ [
                line width=1.0,
                solid,
                mark=none,
            ]
            table [
                x index=0, y expr=-\thisrowno{1},
            ]{data_config1_window5_rhs_ITF.txt};

            \addplot+ [
                line width=1.0,
                solid,
                smooth,
                mark=none,
                blue,
            ]
            table [
                x index=0, y expr=-\thisrowno{1},
            ]{data_config1_window5_rhs_pinn0.txt};

            \addplot+ [
                line width=1.0,
                solid,
                smooth,
                mark=none,
                blue,
            ]
            table [
                x index=0, y expr=-\thisrowno{1},
            ]{data_config1_window5_rhs_pinn1.txt};

            \addplot+ [
                line width=1.0,
                dashed,
                mark=none,
                black,
            ]
            table [
                x index=0, y expr=-\thisrowno{1},
            ]{data_config1_window5_rhs_total.txt};

  \end{groupplot}
\node[below = 1.5cm of my plots c1r1.south west,
    anchor=west,
] {(a) $\tW^{(1)}_{int}$, $\tW^{(1)}_{d}$, $\tW^{(1)}_{n}$};
\node[below = 1.5cm of my plots c2r1.south west,
    anchor=west,
] {(b) $\tW^{(2)}_{int}$, $\tW^{(2)}_{d}$, $\tW^{(2)}_{n}$};
\end{tikzpicture}
%
    \caption{Each subdomain solution's contribution to the left-hand side of the physics equation \eqref{eq:poisson}:
    (a) the worst case $\tW^{(1)}_{int}$, $\tW^{(1)}_{d}$, $\tW^{(1)}_{n}$; and
    (b) the best case $\tW^{(2)}_{int}$, $\tW^{(2)}_{d}$, $\tW^{(2)}_{n}$.}
    \label{fig:config1-rhs}
\end{figure}
Figure~\ref{fig:config1-rhs}~(a) shows each subdomain solution's contribution to the left-hand side of the physics equation \eqref{eq:poisson}, i.e. $\LHS[u]\equiv-\nabla\cdot\kappa\nabla u$,
for the case of $\tW^{(1)}_{int}$, $\tW^{(1)}_{d}$, $\tW^{(1)}_{n}$.
The sum of all contributions should approximately amount to $f(\bx)=1$ for the model problem in Section~\ref{subsec:1interface}.
Clearly, $\LHS[u^{(B)}_b]$ and $\LHS[u^{(I)}_{ij}]$ do not vanish at the end of their subdomains,
introducing an undesired discontinuity at $x=0.25, 0.75$.
$\LHS[u^{(S)}_m]$ does compensate those from the boundary/interface to a certain extent;
however it suffers from the well-known spectral bias when resolving the unwanted discontinuities~\cite{rahaman2019spectral}.
This discontinuity cannot be resolved with other interior window functions,
yielding high solution errors for all cases of $\tW^{(1)}_{d}$, $\tW^{(1)}_{n}$.
\par
The discontinuity disappears with higher order window functions,
as their second derivatives vanish at the subdomain boundary $\tau=1$.
Figure~\ref{fig:config1-rhs}~(b)
shows $\LHS[u]$ of each subdomain solution
for the case of $\tW^{(2)}_{int}$, $\tW^{(2)}_{d}$, $\tW^{(2)}_{n}$.
As $\LHS[u^{(B)}_b]$ and $\LHS[u^{(I)}_{ij}]$ vanishes at $x=0.25, 0.75$,
the interior component $\LHS[u^{(S)}_m]$ only has to resolve a continuous function,
which leads to a significant improvement in solution accuracy.
Interestingly, enforcing smoothness of $\LHS[u]$ did not improve the solution accuracy.
As shown in Figure~\ref{fig:config1}~(b),
the use of higher order polynomials $\tW_{int}^{(3)}$, $\tW_{d}^{(3)}$, or $\tW_{n}^{(3)}$
rather degrades the solution accuracy,
though their $\LHS[u]$ becomes smooth at $x=0.25, 0.75$.
This is presumably because the shapes of second derivatives match the best at $k=2$ for all windows,
which makes it easier for interior NNs to resolve $\LHS[u]$ for the physics equation.
The impact of the shape of second derivatives has been further discussed in Section~\ref{subsec:window-vs-buffer}.

\subsection{Effect of boundary/interface subdomain size in windowing approach}\label{subsec:comparison-overlap}
From the previous section,
we also observe from Figure~\ref{fig:config1-rhs}
that $\LHS$ from the boundary or interface are significantly larger than the actual forcing term $f(\bx)$.
This can make the interior NNs focus more on approximating these $\LHS$ contributions rather than actual forcing,
deteriorating the overall solution accuracy.
A way to reduce the boundary/interface $\LHS$ contribution
is to increase the size of the boundary/interface subdomain, i.e., $\Delta \bx_b$ and $\Delta \bx_{ij}$.
Since derivatives of window functions will be multiplied by the factor of $\Delta \bx^{-1}$
($\Delta \bx^{-2}$ for second derivatives),
increasing $\Delta \bx_b$ and $\Delta \bx_{ij}$ will reduce the boundary/interface $\LHS$ contributions.
This motivates us to investigate the impact of the boundary/interface subdomain size on the training performance.
\par
\begin{figure}[H]
    \begin{tikzpicture}[
]
    \begin{groupplot}[
        group style={
            group name = my plots,
            group size= 2 by 1,
            xlabels at =edge bottom,
            horizontal sep=2cm,
            vertical sep=2.2cm,
        },
        height = 0.45\textwidth,
        width = 0.5\textwidth,
        name=chung,
    ]    
\pgfplotsset{set layers=standard}%

        \nextgroupplot[
            xlabel={$x$},
            ylabel={$W(x)$},
            tick scale binop ={\times},
            xmin=0., xmax=1.,
            ymin=0., ymax=1.1,
            xtick={0, 0.5, 1.},
            legend style={
                draw=none, fill=none,
                at={(rel axis cs: -0.2, 1.0)},
                anchor=south west,
                nodes={scale=1.0},
                legend cell align={left},
                legend columns=3,
                /tikz/every even column/.append style={column sep=0.5cm},
            },
        ]

            \addlegendimage{
                line width=1.0,
                solid,
                mark=none,
                green,};
            \addlegendentry{$W_m(\bx)$, interior};
            \addlegendimage{
                line width=1.0,
                dashed,
                mark=none,
                blue,};
            \addlegendentry{$W_{b,d}(\bx)$, boundary};
            \addlegendimage{
                line width=1.0,
                dashed,
                mark=none,
                red,};
            \addlegendentry{$W_{b,n}(\bx)$, boundary};
            \addlegendimage{
                line width=1.0,
                dashed,
                mark=none,
                brown,};
            \addlegendentry{$W_{ij,d}(\bx)$, interface};
            \addlegendimage{
                line width=1.0,
                dashed,
                mark=none,
                purple,};
            \addlegendentry{$W_{ij,n}(\bx)$, interface};

            \addplot+[
                mark=none,
                gray,
                line width=1.5,
            ] coordinates {(0.5, 0) (0.5, 2)};

            \addplot+ [
                line width=1.0,
                solid,
                mark=none,
                green,
            ]
            table [
                x index=0, y index=1,
            ]{data_config2_window_function_int0L.txt};
            \addplot+ [
                line width=1.0,
                solid,
                mark=none,
                green,
            ]
            table [
                x index=0, y index=1,
            ]{data_config2_window_function_int0M.txt};
            \addplot+ [
                line width=1.0,
                solid,
                mark=none,
                green,
            ]
            table [
                x index=0, y index=1,
            ]{data_config2_window_function_int0R.txt};
            \addplot+ [
                line width=1.0,
                solid,
                mark=none,
                green,
            ]
            table [
                x index=0, y index=1,
            ]{data_config2_window_function_int1L.txt};
            \addplot+ [
                line width=1.0,
                solid,
                mark=none,
                green,
            ]
            table [
                x index=0, y index=1,
            ]{data_config2_window_function_int1M.txt};
            \addplot+ [
                line width=1.0,
                solid,
                mark=none,
                green,
            ]
            table [
                x index=0, y index=1,
            ]{data_config2_window_function_int1R.txt};

            \addplot+ [
                line width=1.0,
                dashed,
                mark=none,
                red,
            ]
            table [
                x index=0, y index=1,
            ]{data_config2_window_function_bndL.txt};
            \addplot+ [
                line width=1.0,
                dashed,
                mark=none,
                red,
            ]
            table [
                x index=0, y index=1,
            ]{data_config2_window_function_bndR.txt};
            \addplot+ [
                line width=1.0,
                dashed,
                mark=none,
                purple,
            ]
            table [
                x index=0, y index=1,
            ]{data_config2_window_function_itf.txt};

            \addplot+ [
                line width=1.0,
                dashed,
                mark=none,
                blue,
            ]
            table [
                x index=0, y index=2,
            ]{data_config2_window_function_bndL.txt};
            \addplot+ [
                line width=1.0,
                dashed,
                mark=none,
                blue,
            ]
            table [
                x index=0, y index=2,
            ]{data_config2_window_function_bndR.txt};
            \addplot+ [
                line width=1.0,
                dashed,
                mark=none,
                brown,
            ]
            table [
                x index=0, y index=2,
            ]{data_config2_window_function_itf.txt};

        \nextgroupplot[
            tick scale binop ={\times},
            ymode=log,
            xlabel={$\beta$},
            ylabel={Relative $L_2$ error},
            legend pos=south west,
            legend style={draw=none},
        ]
        
            \addplot+ [
                mark=*,
                blue,
                mark options={fill=white,},
            ]
            table [
                x index=0, y index=1,
            ]{data_config2_bnd_factors_window0_l2_error.txt};

            \addplot+ [
                mark=*,
                mark options={fill=white,},
            ]
            table [
                x index=0, y index=1,
            ]{data_config2_bnd_factors_window1_l2_error.txt};

            \addplot+ [
                mark=*,
                mark options={fill=white,},
            ]
            table [
                x index=0, y index=1,
            ]{data_config2_bnd_factors_window2_l2_error.txt};

            \legend{{$\tW^{(1)}_{d}$, $\tW^{(1)}_{n}$}, {$\tW^{(2)}_{d}$, $\tW^{(2)}_{n}$}, {$\tW^{(3)}_{d}$, $\tW^{(3)}_{n}$}};

  \end{groupplot}
\draw[<->, solid,black, line width=1.0,] (my plots c1r1.south west) + (0pt, -20pt) -- node[anchor=south] {$\beta \Delta x$} ++(2cm, -20pt);
\draw[-, solid,black, line width=0.5,] (my plots c1r1.south west) + (0pt, -15pt) -- ++(0pt, -25pt);
\draw[-, solid,black, line width=0.5,] (my plots c1r1.south west) + (2cm, 0pt) -- ++(2cm, -25pt);
\node[below = 2cm of my plots c1r1.south west,
    anchor=west,
] {(a) Configuration};
\node[below = 2cm of my plots c2r1.south west,
    anchor=west,
] {(b) Relative $L_2$ error};
\end{tikzpicture}
%
    \caption{Impact of the boundary/interface subdomain size:
    (a) window configuration with various boundary/interface element size; and
    (b) Relative $L_2$ error of the trained solutions depending on the boundary/interface element size.
    For interior subdomains, $\tW_{int}^{(1)}$ is used for all cases.}
    \label{fig:config2}
\end{figure}
The solution ansatz \eqref{eq:window-ansatz} for the windowing approach is
trained with the target problem in Section~\ref{subsec:1interface},
with a varying boundary/interface subdomain.
Figure~\ref{fig:config2}~(a) shows the window configuration for this investigation.
The same subdomain locations from the previous section are used.
While a uniform subdomain size $\Delta x=0.25$ is used for the interior subdomains,
the boundary/interface subdomain sizes are multiplied by a factor $\beta$, i.e. $\beta\Delta x$,
with $\beta$ ranging from 1 to 2.
At $\beta=2$, a boundary/interface subdomain reaches the boundary or the interface on the opposite side.
\par
For this investigation,
various polynomial orders for $\tW^{(k)}_d$ and $\tW^{(k)}_n$ in Table~\ref{tab:wd}-\ref{tab:wn}
are considered,
while $\tW^{(1)}_{int}$ is used for interior subdomains in all cases.
A fully connected network with hidden layers $[25, 25]$ is
used for each interior subdomain.
For all cases, the training loss \eqref{eq:J-int} is evaluated
with $60$ collocation points uniformly distributed in the domain.
The training is performed via SOAP optimizer~\cite{vyas2024soap},
with learning rate $5\times10^{-3}$ for $10^4$ iterations.
\par
Figure~\ref{fig:config2}~(b) shows the relative $L_2$ errors 
of the trained solutions depending on the subdomain size factor $\beta$.
For $\tW^{(2)}$ and $\tW^{(3)}$,
the solution accuracy is gradually improved with increasing $\beta$.
At $\beta=2$, the relative error decreased by an order of magnitude, compared to $\beta=1$.
This is presumably because the LHS contribution from the boundary or interface decreases
with increasing subdomain size, as discussed earlier.
\par
Meanwhile,
the solution accuracy with $\tW^{(1)}$ does not seem to be affected by $\beta$ up to $\beta=1.9$.
This trend drastically changes at $\beta=2$,
achieving the lowest relative error of $2.55\times10^{-7}$.
This behavior is mainly driven by the discontinuities introduced by $\LHS[u^{(B)}_b]$ and $\LHS[u^{(I)}_{ij}]$,
which were discussed in \ref{subsec:comparison-window-order}.
For $\beta < 2$, these discontinuities still persist in the domain,
degrading the training convergence regardless of the subdomain size.
At $\beta=2$, however, these discontinuities lie exactly on the boundary or the interface,
where the physics equation residual \eqref{eq:J-int} is not evaluated.
\par
\begin{figure}[H]
    \begin{tikzpicture}[
]
    \begin{groupplot}[
        group style={
            group name = my plots,
            group size= 3 by 1,
            xlabels at =edge bottom,
            horizontal sep=2cm,
            vertical sep=2.2cm,
        },
        height = 0.45\textwidth,
        width = 0.5\textwidth,
        enlarge x limits={false, abs value = 5mm},
        enlarge y limits={false, abs value = 5mm},
        xtick={0, 0.25, 0.5, 0.75, 1.},
        name=chung,
    ]    
\pgfplotsset{set layers=standard}%

        \nextgroupplot[
            ymin=-6, ymax=6,
            xlabel={$x$},
            ylabel={$\LHS[u]$},
            tick scale binop ={\times},
            legend style={
                draw=none, fill=none,
                at={(rel axis cs: -0.4, 1.05)},
                anchor=south west,
                nodes={scale=1.0},
                legend cell align={left},
                legend columns=4,
                /tikz/every even column/.append style={column sep=0.5cm},
            },
        ]

            \addlegendimage{
                line width=1.0,
                solid,
                mark=none,
                blue,
            };
            \addlegendentry{$\LHS[u^{(S)}_m]$};
            \addlegendimage{
                line width=1.0,
                solid,
                mark=none,
                green,
            };
            \addlegendentry{$\LHS[u^{(B)}_b]$};
            \addlegendimage{
                line width=1.0,
                solid,
                mark=none,
                red,
            };
            \addlegendentry{$\LHS[u^{(I)}_{ij}]$};
            \addlegendimage{
                line width=1.0,
                dashed,
                mark=none,
                black,
            };
            \addlegendentry{$\LHS[u_{NN}]$};

            \addplot+[
                mark=none,
                gray!50,
                line width=1.5,
            ] coordinates {(0.5, -30) (0.5, 30)};
            \addplot+[
                mark=none,
                gray!50,
                line width=1.5,
            ] coordinates {(0.25, -30) (0.25, 30)};
            \addplot+[
                mark=none,
                gray!50,
                line width=1.5,
            ] coordinates {(0.75, -30) (0.75, 30)};

            \addplot+ [
                line width=1.0,
                solid,
                mark=none,
                green,
            ]
            table [
                x index=0, y expr=-\thisrowno{1},
            ]{data_config2_rhs_rhs_BCL.txt};

            \addplot+ [
                line width=1.0,
                solid,
                mark=none,
                green,
            ]
            table [
                x index=0, y expr=-\thisrowno{1},
            ]{data_config2_rhs_rhs_BCR.txt};

            \addplot+ [
                line width=1.0,
                solid,
                mark=none,
            ]
            table [
                x index=0, y expr=-\thisrowno{1},
            ]{data_config2_rhs_rhs_ITF.txt};

            \addplot+ [
                line width=1.0,
                solid,
                mark=none,
                blue,
            ]
            table [
                x index=0, y expr=-\thisrowno{1},
            ]{data_config2_rhs_rhs_pinn0.txt};

            \addplot+ [
                line width=1.0,
                solid,
                mark=none,
                blue,
            ]
            table [
                x index=0, y expr=-\thisrowno{1},
            ]{data_config2_rhs_rhs_pinn1.txt};

            \addplot+ [
                line width=1.0,
                dashed,
                mark=none,
                black,
            ]
            table [
                x index=0, y expr=-\thisrowno{1},
            ]{data_config2_rhs_rhs_total.txt};

        \nextgroupplot[
            tick scale binop ={\times},
            xmin = 0, xmax = 3,
            ymin = 0, ymax = 3,
            xtick={0.5, 1.5, 2.5},
            xticklabels={{$\tW^{(1)}_{d}$, $\tW^{(1)}_{n}$}, {$\tW^{(2)}_{d}$, $\tW^{(2)}_{n}$}, {$\tW^{(3)}_{d}$, $\tW^{(3)}_{n}$}},
            x tick label style={xshift=20pt, rotate=45,anchor=east},
            ytick={0.5, 1.5, 2.5},
            yticklabels={$\tW^{(3)}_{int}$, $\tW^{(2)}_{int}$, $\tW^{(1)}_{int}$},
        ]
        
            \edef\imagepath{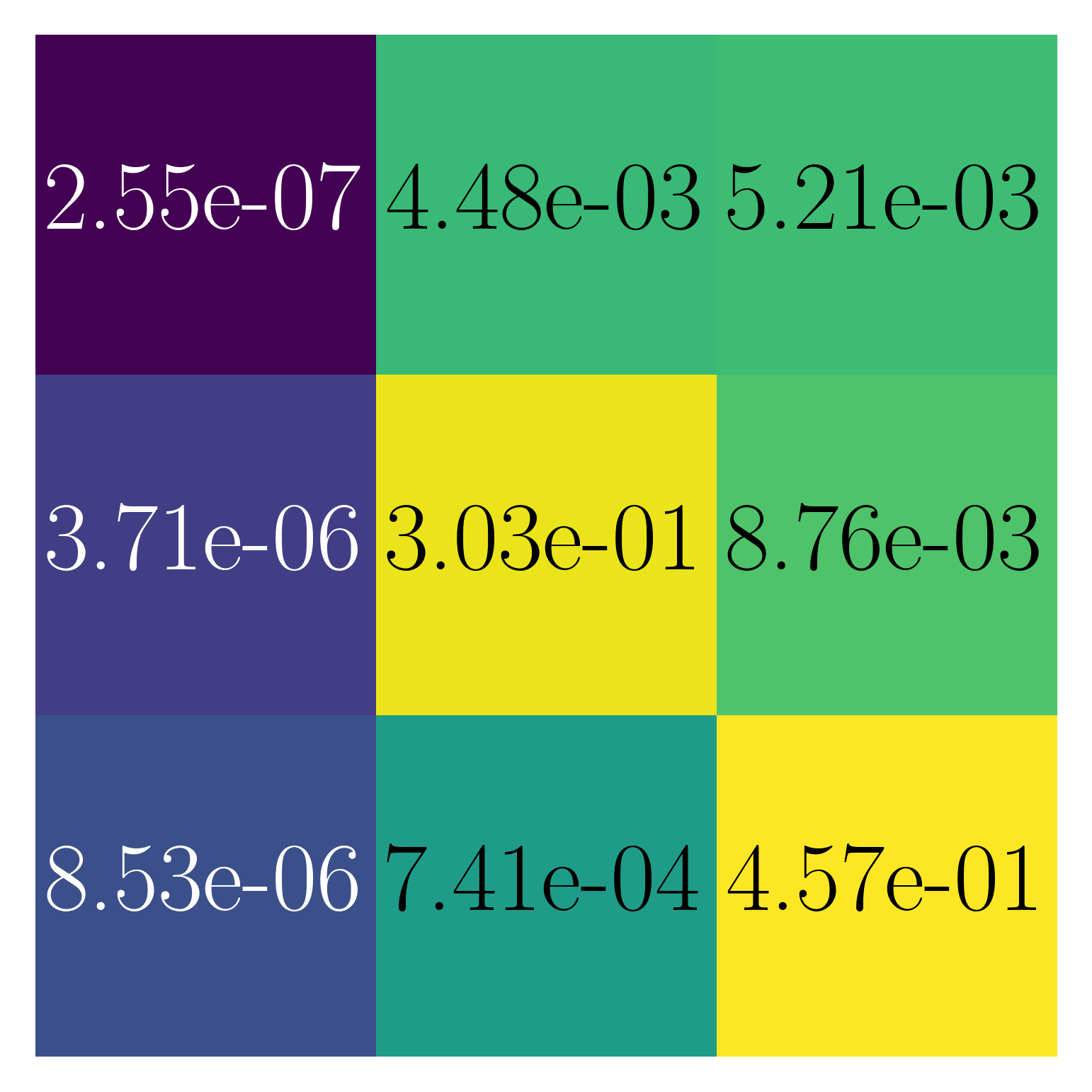}
            \addplot graphics[xmin=-0.1, xmax=3.1, ymin=-0.1, ymax=3.1]{\imagepath};

  \end{groupplot}
\node[below = 2cm of my plots c1r1.south west,
    anchor=west,
] {(a) $\tW^{(1)}_{int}$, $\tW^{(1)}_{d}$, $\tW^{(1)}_{n}$};
\node[below = 2cm of my plots c2r1.south west,
    anchor=west,
] {(b) Relative $L_2$ error};
\end{tikzpicture}
%
    \caption{Windowing approach with perfectly overlapping boundary/interface subdomains ($\beta=2$):
    (a) the left-hand side contribution of each subdomain for the best case of $\tW^{(1)}_{int}$, $\tW^{(1)}_{d}$, $\tW^{(1)}_{n}$; and
    (b) the relative $L_2$ error of the trained solution with various window function combinations.}
    \label{fig:config2-rhs}
\end{figure}
The shape of the window function also significantly contributes to this error reduction.
Figure~\ref{fig:config2-rhs}~(a) shows the $\LHS$ contribution of each subdomain.
Since $\tW^{(1)}_d$ and $\tW^{(1)}_n$ are cubic polynomials,
$\LHS[u^{(B)}_b]$ and $\LHS[u^{(I)}_{ij}]$ remain to be linear functions in interior subdomains,
which can be resolved by interior NNs.
In fact,
these $\LHS$ are summed up to be exactly the same as the forcing term $f(x)=1$,
so interior subdomain solutions $u^{(S)}_m$ merely approximates the trivial solution $u=0$.
This behavior occurs only at $(\tW^{(1)}_d, \tW^{(1)}_n)$ case.
Figure~\ref{fig:config2-rhs}~(b) shows the solution errors of different window function combinations
in the case of perfectly overlapping boundary and interface, i.e., $\beta=2$.
While all cases of $(\tW^{(1)}_d, \tW^{(1)}_n)$ achieve relative errors less than $10^{-4}$,
higher polynomials could not achieve the same accuracy even at $\beta=2$.
\par
This is, of course, a very particular case, and we do not expect to achieve similar solution accuracy for general problems.
In general,
$\LHS$ of boundary/interface subdomains would not exactly match the forcing term
and the interior subdomains would have to resolve the residual for the physics equation.
However, boundary/interface $\LHS$ with $(\tW^{(1)}_d, \tW^{(1)}_n)$
will remain constant or linear in general cases,
which is significantly beneficial for fast training convergence.
\par
\refB{While} perfectly overlapping boundary/interface subdomains ($\beta=2$) with the windows $(\tW^{(1)}_d, \tW^{(1)}_n)$ \refB{turned out to be} the optimal configuration \refB{here, this finding should be interpreted as an empirical guideline for the one-dimensional problems, rather than a generally optimal choice. The extension of this configuration to higher dimensions and more general geometries is not straightforward, and its consequences for the two-dimensional windowing experiment is discussed in Section~\ref{subsubsec:window2d}.}

\section{Ablation test on initialization}\label{app:init}
\refB{This section presents an ablation study on how the compared methods respond to random variations in both the neural-network initialization and the piecewise diffusivities. Six PINN methods are compared on the three one-dimensional problems. For each method, the network weights are drawn from two initializers (\texttt{glorot\_uniform} and \texttt{random\_normal}) at two scales ($1$ and $0.1$), with $100$ independent random seeds per initializer--scale combination. The piecewise diffusivities are likewise sampled independently for each seed: for Problems~1 and~3, $k_1\sim\mathrm{LogU}[10^{-2},10^{-1}]$ and $k_2\sim\mathrm{LogU}[10^{0},10^{1}]$; for Problem~2, $k_1,\ldots,k_4\sim\mathrm{LogU}[10^{-2},10^{1}]$ independently. Each cell reports the $(\min,\mathrm{median},\max)$ of the relative $L_2$ error over the $100$ seeds. Within each problem table, the overall best minimum, median, and maximum across all methods and all four initializer--scale columns are highlighted in \textcolor{blue}{blue}, \textcolor{green}{green}, and \textcolor{red}{red}, respectively.}

\begin{table}[h!]
    \centering
    \resizebox{\textwidth}{!}{%
    \begin{tabular}{|l|c|c|c|c|}
        \hline
        Problem 1 & \multicolumn{2}{c|}{\texttt{glorot\_uniform}} & \multicolumn{2}{c|}{\texttt{random\_normal}} \\
        \hline
        Method & scale $=0.1$ & scale $=1$ & scale $=0.1$ & scale $=1$ \\
        \hline
        $\phi$-PINNs & (2.68e-5, 1.20e-4, 9.99e-1) & (2.99e-5, 1.23e-4, 1.54e0) & (1.89e-5, 9.83e-5, 1.29e-3) & (2.12e-5, 1.10e-4, 1.03e0) \\
        I-PINNs      & (7.73e-5, 7.71e-4, 1.90e-1) & (1.36e-4, 1.14e-3, 1.03e-1) & (8.28e-5, 1.00e-3, 4.03e-1) & (4.38e-5, 3.77e-4, 1.62e0) \\
        AdaI-PINNs   & (3.95e-5, 4.69e-4, 1.00e0)  & (1.37e-5, 4.36e-4, 1.00e0)  & (4.56e-5, 3.60e-4, 1.45e0)  & (3.17e-5, 3.45e-4, 1.00e0) \\
        M-PINN       & (1.13e-5, 8.67e-5, 3.09e-2) & (1.23e-5, 7.35e-5, 3.58e0)  & (8.37e-6, 6.96e-5, 9.99e-1) & (1.41e-5, 7.70e-5, 2.36e0) \\
        \hline
        Window       & (2.12e-7, 8.27e-2, 2.64e0)  & (5.72e-7, 6.96e-2, 6.52e-1) & (\textcolor{blue}{2.10e-9}, 8.88e-2, 1.10e0) & (4.20e-6, 5.76e-2, 6.67e0) \\
        Buffer       & (1.26e-8, \textcolor{green}{9.93e-6}, 1.06e1) & (2.95e-7, 4.09e-5, \textcolor{red}{3.36e-4}) & (8.38e-8, 1.37e-5, 2.26e1) & (3.41e-6, 3.45e-5, 6.04e-4) \\
        \hline
    \end{tabular}%
    }
    \caption{\refB{Per-(initializer, scale) statistics for Problem~1. Each cell is $(\min,\mathrm{median},\max)$ of the relative $L_2$ error over $100$ random seeds.}}
    \label{tab:app-init-prob1}
\end{table}

\begin{table}[h!]
    \centering
    \resizebox{\textwidth}{!}{%
    \begin{tabular}{|l|c|c|c|c|}
        \hline
        Problem 2 & \multicolumn{2}{c|}{\texttt{glorot\_uniform}} & \multicolumn{2}{c|}{\texttt{random\_normal}} \\
        \hline
        Method & scale $=0.1$ & scale $=1$ & scale $=0.1$ & scale $=1$ \\
        \hline
        $\phi$-PINNs & (3.26e-5, 4.00e-4, 8.36e-1) & (4.05e-5, 3.29e-4, 4.72e-1) & (5.10e-5, 4.05e-4, 2.60e-1) & (2.04e-5, 5.05e-4, 1.74e0) \\
        I-PINNs      & (2.66e-4, 6.37e-2, 8.79e-1) & (1.40e-3, 6.92e-2, 8.12e-1) & (8.76e-4, 5.11e-2, 5.48e-1) & (1.56e-3, 8.83e-2, 8.44e-1) \\
        AdaI-PINNs   & (5.57e-5, 7.17e-3, 9.92e-1) & (5.39e-5, 1.01e-2, 6.86e0)  & (3.66e-5, 6.22e-3, 1.01e0)  & (6.96e-5, 8.10e-3, 1.11e0) \\
        M-PINN       & (7.99e-6, 1.07e-4, 6.63e0)  & (1.28e-5, 9.42e-5, 1.38e0)  & (1.90e-5, 1.24e-4, 7.17e-1) & (8.33e-6, 8.84e-5, 3.13e-1) \\
        \hline
        Window       & (\textcolor{blue}{6.06e-8}, 4.91e-4, 3.88e0) & (1.30e-6, 2.82e-4, 1.38e0) & (3.31e-7, 1.49e-4, 1.75e0) & (1.19e-5, 4.81e-3, 1.73e0) \\
        Buffer       & (4.03e-6, 7.45e-5, \textcolor{red}{3.22e-3}) & (5.97e-6, \textcolor{green}{6.50e-5}, 7.07e0) & (5.90e-6, 6.99e-5, 3.21e-1) & (1.33e-5, 1.50e-4, 1.12e2) \\
        \hline
    \end{tabular}%
    }
    \caption{\refB{Per-(initializer, scale) statistics for Problem~2. Same format as Table~\ref{tab:app-init-prob1}.}}
    \label{tab:app-init-prob2}
\end{table}

\begin{table}[h!]
    \centering
    \resizebox{\textwidth}{!}{%
    \begin{tabular}{|l|c|c|c|c|}
        \hline
        Problem 3 & \multicolumn{2}{c|}{\texttt{glorot\_uniform}} & \multicolumn{2}{c|}{\texttt{random\_normal}} \\
        \hline
        Method & scale $=0.1$ & scale $=1$ & scale $=0.1$ & scale $=1$ \\
        \hline
        $\phi$-PINNs & (9.18e-4, 6.63e-2, 9.98e-1) & (6.05e-4, 6.33e-2, 2.67e-1) & (1.18e-3, 6.81e-2, 5.26e-1) & (7.32e-4, 6.55e-2, 1.10e0) \\
        I-PINNs      & (1.87e-2, 2.49e-1, 1.67e2)  & (2.62e-3, 1.08e-1, 1.06e0)  & (7.91e-3, 3.35e-1, 1.12e0)  & (1.03e-3, 7.04e-2, 1.03e0) \\
        AdaI-PINNs   & (1.36e-2, 1.01e0, 1.59e0)   & (1.61e-2, 1.04e0, 1.68e0)   & (2.03e-2, 1.04e0, 1.44e0)   & (3.02e-2, 1.03e0, 7.64e1) \\
        M-PINN       & (9.59e-4, 6.82e-2, 1.05e0)  & (1.47e-3, 6.30e-2, 2.66e-1) & (2.13e-3, 6.86e-2, 1.02e0)  & (7.48e-4, 6.73e-2, 8.91e2) \\
        \hline
        Window       & (1.33e-3, 6.52e-2, 7.35e0)  & (9.51e-4, 6.48e-2, 1.07e2)  & (9.36e-4, 6.28e-2, 2.67e-1) & (1.01e-3, 6.47e-2, 3.58e0) \\
        Buffer       & (\textcolor{blue}{5.90e-4}, 6.27e-2, \textcolor{red}{2.66e-1}) & (9.98e-4, 6.56e-2, 2.70e-1) & (7.72e-4, \textcolor{green}{6.15e-2}, 2.04e2) & (8.20e-4, 6.18e-2, 2.71e-1) \\
        \hline
    \end{tabular}%
    }
    \caption{\refB{Per-(initializer, scale) statistics for Problem~3. Same format as Table~\ref{tab:app-init-prob1}.}}
    \label{tab:app-init-prob3}
\end{table}

\section{Full hard-constraining windowing approach
on two-dimensional Poisson equation}
\label{app:2dpoisson-full-window}

In this appendix, we present the numerical results of the windowing formulation (Section~\ref{subsec:windowing})
applied to the two-dimensional Poisson problem introduced in Section~\ref{subsec:2dpoisson},
with all boundary and interface conditions enforced as hard constraints.
Unlike the partially hard-constrained setting discussed in
Section~\ref{subsec:demo-2dpoisson}---where only the interface conditions were imposed exactly---here
both Dirichlet and Neumann boundary conditions,
as well as the interface continuity and flux balance conditions,
are embedded directly into the solution ansatz through the window
and corner constructions described in Section~\ref{subsec:window-multi-dim}.
\par
\begin{table}[ht]
\centering
\caption{Window element parameters. Interior half-widths marked with ${}^*$ are in reference coordinates $(\xi_1, \xi_2) \in [0,1]^2$. Interface edge sizes are expressed as fractions of the interface length $L_{\mathrm{itf}} = \sqrt{29}/5 \approx 1.08$. Interface corner angle widths are given as approximate (left, right) subdomain pairs.}
\label{tab:window_params}
\begin{tabular}{l c c c}
\hline
Window & Center & Normal size & Tangential size \\
\hline
Left interior   & $(0.50,\, 0.50)$ & $0.50^*$ & $0.50^*$ \\
Right interior  & $(1.50,\, 0.50)$ & $0.50^*$ & $0.50^*$ \\
\hline
Left edge        & $(0.00,\, 0.50)$ & $0.50$ & $0.50$ \\
Bottom-left edge & $(0.40,\, 0.00)$ & $0.50$ & $0.40$ \\
Top-left edge    & $(0.54,\, 1.00)$ & $0.30$ & $0.54$ \\
Right edge       & $(2.00,\, 0.50)$ & $0.50$ & $0.50$ \\
Bottom-right edge & $(1.46,\, 0.00)$ & $0.30$ & $0.54$ \\
Top-right edge   & $(1.60,\, 1.00)$ & $0.50$ & $0.40$ \\
Interface edge   & $(1.00,\, 0.50)$ & $0.25\, L_{{itf}}$ & $0.40\, L_{{itf}}$ \\
\hline
 & Center & Radius & Angle width ($\times \pi$) \\
\hline
Bottom-left corner  & $(0.00,\, 0.00)$ & $0.40$ & $1/2$ \\
Top-left corner     & $(0.00,\, 1.00)$ & $0.50$ & $1/2$ \\
Bottom-right corner & $(2.00,\, 0.00)$ & $0.50$ & $1/2$ \\
Top-right corner    & $(2.00,\, 1.00)$ & $0.40$ & $1/2$ \\
Bottom interface corner & $(0.80,\, 0.00)$ & $0.40$ & $\approx(0.62,\; 0.38)$ \\
Top interface corner    & $(1.20,\, 1.00)$ & $0.40$ & $\approx(0.38,\; 0.62)$ \\
\hline
\end{tabular}
\end{table}
\begin{figure}[tbh]
    \begin{tikzpicture}
    \begin{axis}[
        ylabel={$y$},
        xlabel={$x$},
        xmin = -0.02, xmax = 2.02,
        ymin = -0.02, ymax = 1.02,
        height = 0.4\textwidth,
        width = 0.65\textwidth,
        axis on top,
        legend style={
            at={(1.02,1)},
            anchor=north west,
            draw=black,
            font=\small,
            legend cell align={left},
        },
    ]

        \addlegendimage{
            only marks, mark=*, mark size=2.5pt, teal}
        \addlegendentry{Interior subdomains}

        \addlegendimage{
            only marks, mark=*, mark size=2.5pt, blue}
        \addlegendentry{Left subdomain corners}

        \addlegendimage{
            only marks, mark=+, mark size=3pt,
            blue, thick}
        \addlegendentry{Left subdomain edges}

        \addlegendimage{
            only marks, mark=*, mark size=2.5pt, red}
        \addlegendentry{Right subdomain corners}

        \addlegendimage{
            only marks, mark=+, mark size=3pt,
            red, thick}
        \addlegendentry{Right subdomain edges}

        \addlegendimage{
            only marks, mark=x, mark size=3pt,
            orange, thick}
        \addlegendentry{Interface edges}

        \addlegendimage{
            only marks, mark=*, mark size=2.5pt,
            orange}
        \addlegendentry{Interface corners}

        \edef\imagepath{%
            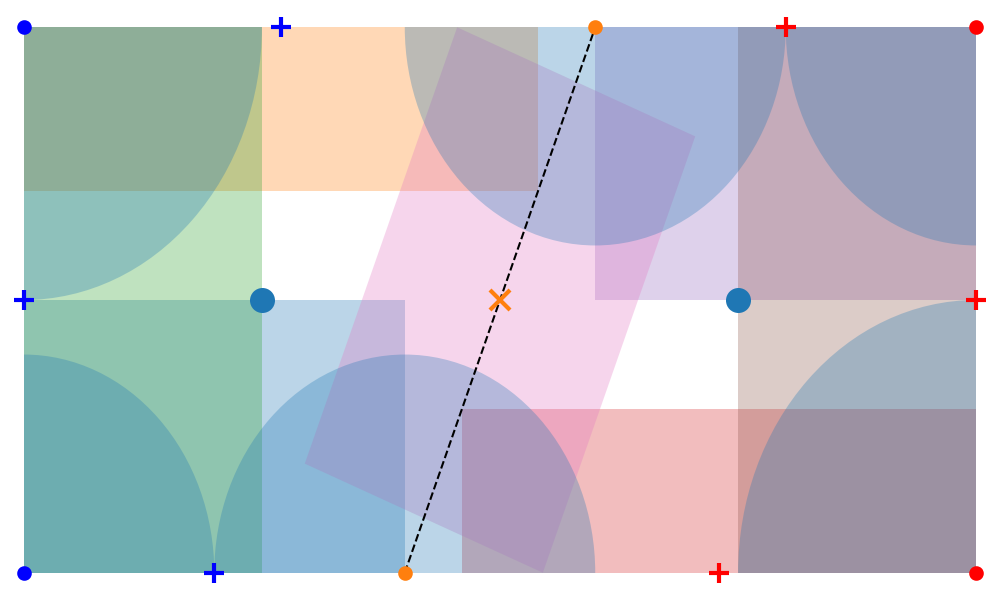}
        \addplot graphics[
            xmin=-0.05, xmax=2.05,
            ymin=-0.05, ymax=1.05
        ]{\imagepath};

    \end{axis}
\end{tikzpicture}
    \caption{Spatial layout of window elements for the fully hard-constrained
    windowing approach applied to the two-dimensional Poisson problem.}
    \label{fig:2dpoisson_window_full_extent}
\end{figure}
Table~\ref{tab:window_params} summarizes the geometric configuration
of all window elements used in the fully hard-constrained two-dimensional windowing implementation.
Figure~\ref{fig:2dpoisson_window_full_extent} visualizes this layout,
showing the spatial placement and overlap of the interior, edge, interface, and corner windows in the domain.
The interior windows are centered at the midpoints of the left and right subdomains and are defined
in reference coordinates with symmetric half-widths.
A guiding principle in selecting window sizes was that the support of each window
should not extend beyond neighboring nodes.
Edge windows are placed along each physical boundary and along the slanted interface,
with their normal extents chosen to ensure sufficient overlap with the adjacent interior windows.
The tangential sizes of the top-left edge, bottom-right edge, and interface edge windows
were chosen slightly smaller than the corresponding physical edge lengths.
This choice prevents window support from passing beyond domain boundaries or
interface endpoints while still covering the majority of each boundary or interface segment
to maintain effective constraint enforcement and overlap.
Corner windows are constructed in polar coordinates at each physical and interface junction,
with specified radii and angular spans to enforce compatibility between intersecting boundary
and interface conditions.
\par
For the window functions, we adopt $\tW_d^{(3)}$ and $\tW_n^{(3)}$,
the smoothest choices among the target window functions,
motivated by the observations in \ref{subsec:comparison-overlap}
that higher-order windows exhibit greater robustness under partial overlap configurations.
The solution ansatz employs a total of ten neural networks: two two-dimensional networks
for the interior subdomains, and seven one-dimensional networks defined along the tangential
coordinates of the domain boundaries and the interface.
The two-dimensional networks use the architecture $[2,32,32,32,1]$,
while the one-dimensional networks use $[1,25,25,25,1]$.
For the interior subdomains, each solution component is formulated on a unit-square reference domain 
$(\xi_1, \xi_2)\in[0,1]^2$,
with an appropriate spatial transformation applied to map the physical subdomain to this reference space.
Specifically, for the left subdomain we use the transformation
\begin{subequations}
\begin{equation}
\begin{split}
&\xi_1 = \frac{x_1}{\ell_L(y)},
\quad
\xi_2 = y,\\
&\ell_L(y) = x_{\mathrm{itf},b} + (x_{\mathrm{itf},t} - x_{\mathrm{itf},b}) \frac{y - y_{\min}}{y_{\max} - y_{\min}},
\end{split}
\end{equation}
with $x_{\mathrm{itf},b}=0.8$, $x_{\mathrm{itf},t}=1.2$, $y_{\min}=0$ and $y_{\max}=1$.
For the right subdomain, the mapping is defined as
\begin{equation}
\begin{split}
&\xi_1 = \frac{x - x_{\max}}{\ell_R(y)} + 1,
\quad
\xi_2 = y,\\
&\ell_R(y) = (x_{\max} - x_{\mathrm{itf},b}) + (x_{\mathrm{itf},b} - x_{\mathrm{itf},t}) \frac{y - y_{\min}}{y_{\max} - y_{\min}},
\end{split}
\end{equation}
with $x_{\max}=2$.
\end{subequations}
\par
For the evaluation of the loss \eqref{eq:J-int},
it should be noted that the corner window function in \eqref{eq:wc}
introduces a singularity in the Laplacian of the solution at the corner.
To avoid this issue, the loss is instead evaluated in polar coordinates,
with the governing equation multiplied by an $r^2$ factor,
\begin{equation}\label{eq:J-polar}
\tilde{\cJ}_{physics} =
\sum_{k}^K \left[
- r\Dpartial{}{r}\left(r\kappa(\bx_k)\Dpartial{u}{r}\right)
- \Dpartial{}{\alpha}\left(\kappa(\bx_k)\Dpartial{u}{\alpha}\right)
- r^2f(\bx_k)
\right]^2,
\end{equation}
where the polar coordinates $(r, \alpha)$ of collocation points
are defined relative to the nearest corner node.
\par
\begin{figure}[tbh]
\hspace*{-1.8cm}
    \begin{tikzpicture}[
]
    \begin{groupplot}[
        group style={
            group name = my plots,
            group size= 2 by 1,
            xlabels at =edge bottom,
            horizontal sep=3cm,
        },
        height = 0.3\textwidth,
        width = 0.5\textwidth,
    ]

        \nextgroupplot[
            ylabel={$y$},
            xlabel={$x$},
            tick scale binop ={\times},
            xmin = 0, xmax = 2,
            ymin = 0, ymax = 1,
            colorbar,
            colormap name=viridis,
        ]

            \edef\imagepath{%
                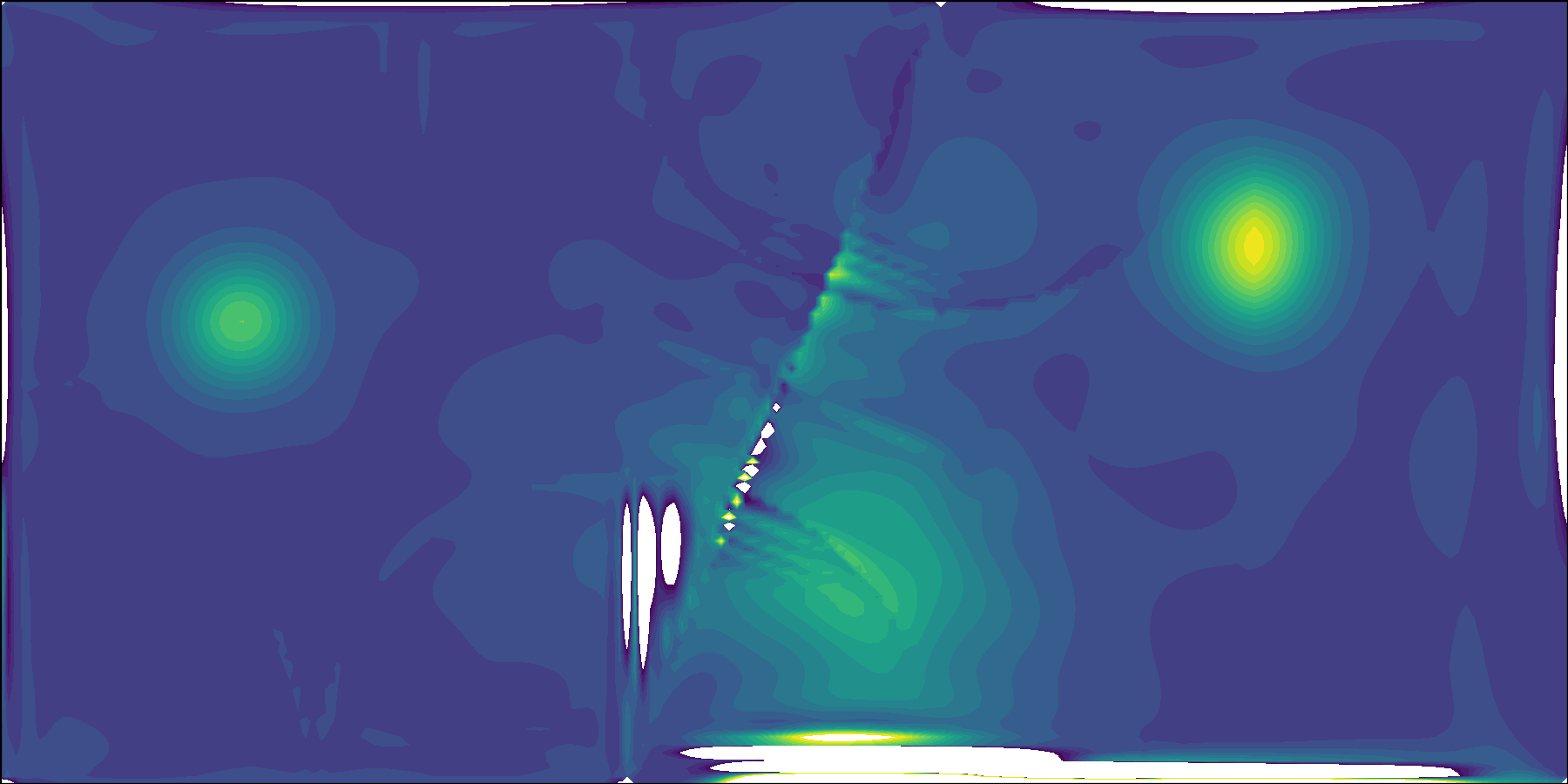}
            \addplot graphics[
                xmin=0, xmax=2,
                ymin=0, ymax=1
            ]{\imagepath};

            \addplot[
                domain=0.8:1.2, samples=100,
                dashed, black,
            ]
            {(x-0.8)/0.4};

        \nextgroupplot[
            ylabel={$y$},
            xlabel={$x$},
            tick scale binop ={\times},
            xmin = 0, xmax = 2,
            ymin = 0, ymax = 1,
            colorbar,
            colormap name=viridis,
        ]

            \edef\imagepath{%
                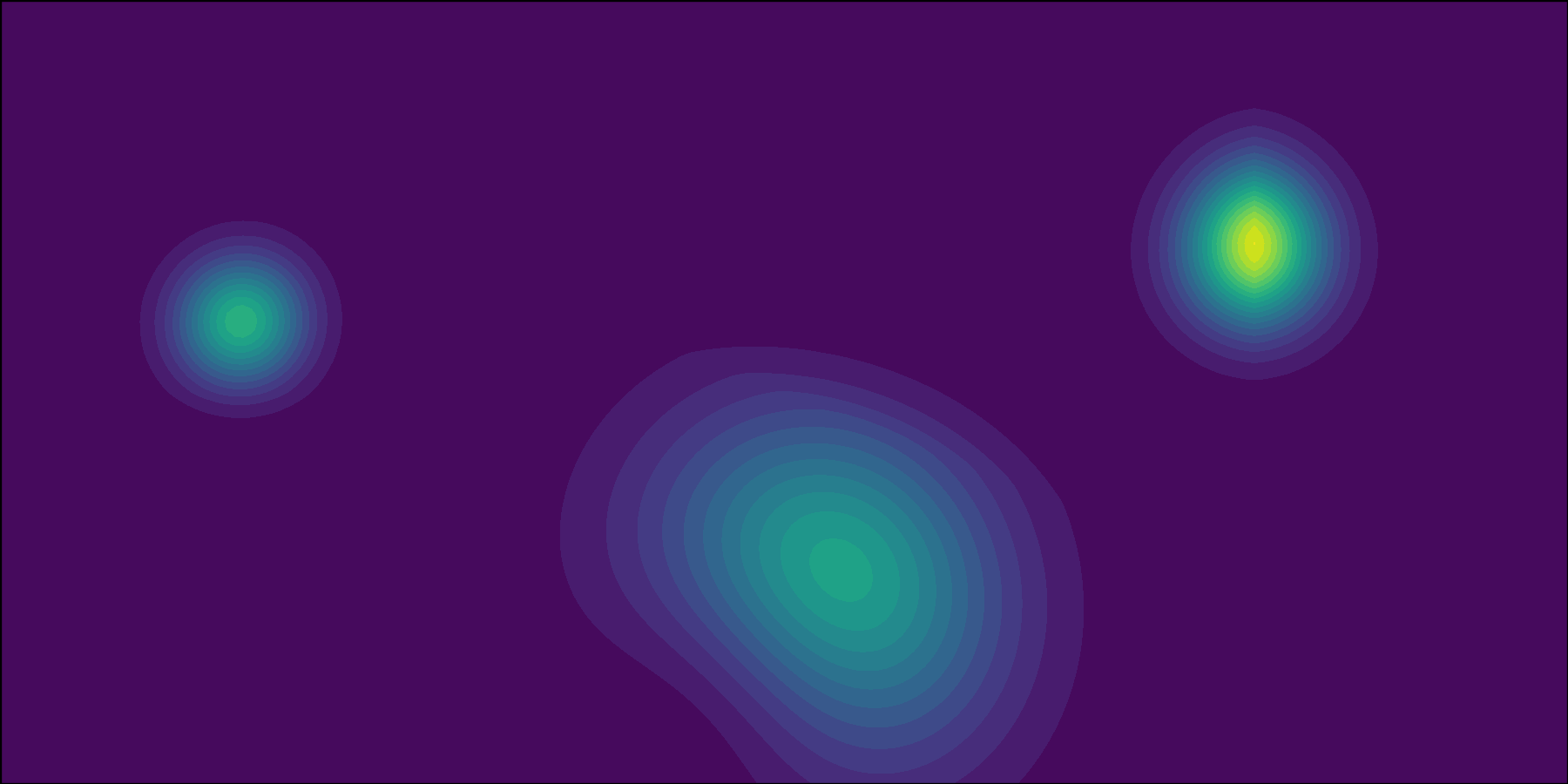}
            \addplot graphics[
                xmin=0, xmax=2,
                ymin=0, ymax=1
            ]{\imagepath};

            \addplot[
                domain=0.8:1.2, samples=100,
                dashed, black,
            ]
            {(x-0.8)/0.4};

    \end{groupplot}
\node[below = 1.5cm of my plots c1r1.south west,
    anchor=west,
] {(a) $r^2\LHS[u]$};
\node[below = 1.5cm of my plots c2r1.south west,
    anchor=west,
] {(b) $r^2f(\bx)$};
\end{tikzpicture}
    \caption{Comparison of the modified residual terms in \eqref{eq:J-polar} after training:
    (a) $r^2\LHS[u]$; and
    (b) $r^2f(\bx)$.}
    \label{fig:2dpoisson_window_full_lhsrhs}
\end{figure}
As in Section~\ref{subsec:demo-2dpoisson},
the model is trained using $80\times40$ collocation points uniformly distributed across the domain,
with $3\times10^4$ epochs under the SOAP optimizer and a learning rate of $10^{-3}$.
However, the solution ansatz exhibits difficulty converging to the correct solution
due to increased stiffness in the model.
Figure~\ref{fig:2dpoisson_window_full_lhsrhs}
presents the left-hand side and the forcing term in \eqref{eq:J-polar} after training.
Although the left-hand side provides a coarse approximation of the forcing term,
pronounced oscillations appear near the boundaries of the window supports,
suggesting that imperfect overlap among the windows contributes to the degraded convergence behavior.
Alternative window sizes and polynomial orders were also examined, yielding qualitatively similar results.

\section{\refB{Sensitivity to the hyperparameters in the buffer approach}}\label{app:buffer2d-sensitivity}
\refB{
This appendix examines how the buffer approach of Section~\ref{subsec:buffer-multi-dim}
responds to boundary-related hyperparameters, using Problem~4 of
Section~\ref{subsec:2dpoisson} as the target problem.
The first group of hyperparameters is the number of sample points at which the conditions
\eqref{eq:g-bc-high-dim} are enforced exactly:
$N_d$ per Dirichlet edge, $N_n$ per Neumann edge, and $N_i$ along the interface.
The second group is the radius of the radial basis functions \eqref{eq:g-rbf},
which we write as
\begin{equation}\label{eq:rho-factor}
r_0 = \rho_s\frac{L_s}{N_s+1},
\end{equation}
with $\rho_s$ a multiplier on the nominal radius, applied separately to the
Dirichlet, Neumann, and interface samples as $\rho_d$, $\rho_n$, and $\rho_i$.
\refA{
The last group is the pair of interface split ratios, $\gamma_0$ for the solution jump and
$\gamma_1$ for the flux, introduced in \eqref{eq:g-itf-ratio} and carried into the
higher-dimensional interface conditions \eqref{eq:g-bc-high-dim}.
}
The configuration reported in Section~\ref{subsubsec:buffer2d} corresponds to
$N_d=4$, $N_n=N_i=8$, $\rho_d=\rho_n=\rho_i=1$ and \refA{$\gamma_0=\gamma_1=1$.}
Throughout this appendix, $\bA_l$ and $\bA_r$ denote the systems that determine the
buffer coefficients in the left and the right subdomain,
the higher-dimensional counterparts of the one-dimensional system~\eqref{eq:gcoeffs}.
All runs reported below use single precision, as does the result in
Section~\ref{subsubsec:buffer2d}.
}
\par
\begin{table}[tbh]
\centering
\caption{\refB{Parameter values of the global sweep.
The reference configuration is marked in \textcolor{red}{red}.}}
\label{tab:2dpoisson-sens-sweep}
\begin{tabular}{l l l}
\hline
Parameter & Description & Values \\
\hline
$N_d$      & sample points per Dirichlet edge & $2$, \textcolor{red}{$4$}, $6$, $8$ \\
$N_n$      & sample points per Neumann edge   & $4$, $6$, \textcolor{red}{$8$}, $10$ \\
$N_i$      & sample points on interface       & $4$, $6$, \textcolor{red}{$8$}, $10$ \\
$\rho_d$   & Dirichlet radius multiplier      & $0.6$, $0.8$, \textcolor{red}{$1.0$}, $1.2$, $1.4$ \\
$\rho_n$   & Neumann radius multiplier        & $0.6$, $0.8$, $1.0$, \textcolor{red}{$1.2$}, $1.4$ \\
$\rho_i$   & interface radius multiplier      & \textcolor{red}{$0.6$}, $0.8$, $1.0$, $1.2$, $1.4$ \\
$\gamma_0$ & interface solution split \eqref{eq:g-bc-high-dim} & $0.1$, \textcolor{red}{$1.0$}, $10.0$ \\
$\gamma_1$ & interface flux split \eqref{eq:g-bc-high-dim}     & $0.1$, \textcolor{red}{$1.0$}, $10.0$ \\
\hline
\end{tabular}
\end{table}
\refB{
Table~\ref{tab:2dpoisson-sens-sweep} lists the parameter values considered in this
sensitivity study.
The three values taken by the interface split ratios $\gamma_0$ and $\gamma_1$
follow from the diffusivity ratio of Problem~4, $\kappa_l/\kappa_r=0.1$:
the jump and the flux are divided between the two subdomains in proportion to that
ratio, in proportion to its reciprocal, and evenly.
Every remaining setting is held at the values of Section~\ref{subsubsec:buffer2d}:
$80\times40$ collocation points, networks of $[2,25,25,25,1]$,
and $3\times10^4$ epochs of the SOAP optimizer at a learning rate of $10^{-3}$
from a fixed initialization seed.
In order to characterize the method beyond a single configuration, we swept the full
Cartesian product of these values, giving $72{,}000$ configurations in total.
}
\par
\refB{
The most accurate case found in this sensitivity study is marked as the reference
in Table~\ref{tab:2dpoisson-sens-sweep}, namely
$N_d=4$, $N_n=N_i=8$, $\gamma_0=\gamma_1=1$, $\rho_d=1$, $\rho_n=1.2$, and $\rho_i=0.6$.
It differs from the case reported in Section~\ref{subsubsec:buffer2d}
only in the two radius multipliers $\rho_n$ and $\rho_i$,
and achieves a lower error of $2.13\times10^{-3}$,
with $\mathrm{cond}(\bA_l)=2.47\times10^{2}$ and $\mathrm{cond}(\bA_r)=3.00\times10^{1}$.
We nonetheless do not report it there, since the other methods were not subjected to a
comparable parameter exploration, and quoting the best of $72{,}000$ buffer configurations
against a single run of each baseline would not be an even-handed comparison.
}
\par
\begin{figure}[tbph]
    \begin{tikzpicture}[
font=\small,
]
    \begin{groupplot}[
        group style={
            group name = my plots,
            group size= 2 by 4,
            horizontal sep=2.2cm,
            vertical sep=1.9cm,
        },
        height = 0.22\textwidth,
        width = 0.42\textwidth,
        ymode=log,
        tick scale binop ={\times},
        tick label style={font=\scriptsize},
        every axis plot/.append style={
            line width=1.0,
        },
        name=chung,
    ]
\pgfplotsset{set layers=standard}%

        \nextgroupplot[
            xlabel={$N_d$},
            ylabel={Relative $L_2$ error},
            xmin=1.5, xmax=8.5, xtick={2,4,6,8},
            ymin=1e-3, ymax=1e0,
        ]

            \addplot [gray, dotted, mark=none, line width=0.8]
            coordinates {(4, 1e-8) (4, 1e8)};

            \addplot [blue, mark=*, mark size=1.8pt,
                      mark options={fill=white,}]
            table [x index=0, y index=1,
            ]{data_2dpoisson_sensitivity_slice_n_dirichlet.txt};

        \nextgroupplot[
            xlabel={$N_d$},
            ylabel={$\mathrm{cond}(\bA_r)$},
            xmin=1.5, xmax=8.5, xtick={2,4,6,8},
            ymin=1e1, ymax=1e3,
        ]

            \addplot [gray, dotted, mark=none, line width=0.8]
            coordinates {(4, 1e-8) (4, 1e8)};

            \addplot [blue, mark=*, mark size=1.8pt,
                      mark options={fill=white,}]
            table [x index=0, y index=3,
            ]{data_2dpoisson_sensitivity_slice_n_dirichlet.txt};

        \nextgroupplot[
            xlabel={$\rho_d$},
            ylabel={Relative $L_2$ error},
            xmin=0.5, xmax=1.5, xtick={0.6,1.0,1.4},
            ymin=1e-3, ymax=1e0,
        ]

            \addplot [gray, dotted, mark=none, line width=0.8]
            coordinates {(1.0, 1e-8) (1.0, 1e8)};

            \addplot [blue, mark=*, mark size=1.8pt,
                      mark options={fill=white,}]
            table [x index=0, y index=1,
            ]{data_2dpoisson_sensitivity_slice_rb_dir_factor.txt};

        \nextgroupplot[
            xlabel={$\rho_d$},
            ylabel={$\mathrm{cond}(\bA_r)$},
            xmin=0.5, xmax=1.5, xtick={0.6,1.0,1.4},
            ymin=1e1, ymax=1e3,
        ]

            \addplot [gray, dotted, mark=none, line width=0.8]
            coordinates {(1.0, 1e-8) (1.0, 1e8)};

            \addplot [blue, mark=*, mark size=1.8pt,
                      mark options={fill=white,}]
            table [x index=0, y index=3,
            ]{data_2dpoisson_sensitivity_slice_rb_dir_factor.txt};

        \nextgroupplot[
            xlabel={$N_n$},
            ylabel={Relative $L_2$ error},
            xmin=3, xmax=11, xtick={4,6,8,10},
            ymin=1e-3, ymax=1e0,
        ]

            \addplot [gray, dotted, mark=none, line width=0.8]
            coordinates {(8, 1e-8) (8, 1e8)};

            \addplot [blue, mark=*, mark size=1.8pt,
                      mark options={fill=white,}]
            table [x index=0, y index=1,
            ]{data_2dpoisson_sensitivity_slice_n_neumann.txt};

        \nextgroupplot[
            xlabel={$N_n$},
            ylabel={$\mathrm{cond}(\bA_l)$},
            xmin=3, xmax=11, xtick={4,6,8,10},
            ymin=1e1, ymax=1e3,
        ]

            \addplot [gray, dotted, mark=none, line width=0.8]
            coordinates {(8, 1e-8) (8, 1e8)};

            \addplot [blue, mark=*, mark size=1.8pt,
                      mark options={fill=white,}]
            table [x index=0, y index=2,
            ]{data_2dpoisson_sensitivity_slice_n_neumann.txt};

        \nextgroupplot[
            xlabel={$\rho_n$},
            ylabel={Relative $L_2$ error},
            xmin=0.5, xmax=1.5, xtick={0.6,1.0,1.4},
            ymin=1e-3, ymax=1e0,
        ]

            \addplot [gray, dotted, mark=none, line width=0.8]
            coordinates {(1.2, 1e-8) (1.2, 1e8)};

            \addplot [blue, mark=*, mark size=1.8pt,
                      mark options={fill=white,}]
            table [x index=0, y index=1,
            ]{data_2dpoisson_sensitivity_slice_rb_neu_factor.txt};

        \nextgroupplot[
            xlabel={$\rho_n$},
            ylabel={$\mathrm{cond}(\bA_l)$},
            xmin=0.5, xmax=1.5, xtick={0.6,1.0,1.4},
            ymin=1e1, ymax=1e3,
        ]

            \addplot [gray, dotted, mark=none, line width=0.8]
            coordinates {(1.2, 1e-8) (1.2, 1e8)};

            \addplot [blue, mark=*, mark size=1.8pt,
                      mark options={fill=white,}]
            table [x index=0, y index=2,
            ]{data_2dpoisson_sensitivity_slice_rb_neu_factor.txt};

  \end{groupplot}
\node[below = 1.2cm of my plots c1r1.south west, anchor=west] {(a)};
\node[below = 1.2cm of my plots c2r1.south west, anchor=west] {(b)};
\node[below = 1.2cm of my plots c1r2.south west, anchor=west] {(c)};
\node[below = 1.2cm of my plots c2r2.south west, anchor=west] {(d)};
\node[below = 1.2cm of my plots c1r3.south west, anchor=west] {(e)};
\node[below = 1.2cm of my plots c2r3.south west, anchor=west] {(f)};
\node[below = 1.2cm of my plots c1r4.south west, anchor=west] {(g)};
\node[below = 1.2cm of my plots c2r4.south west, anchor=west] {(h)};
\end{tikzpicture}
%
    \caption{\refB{A local sensitivity study for the boundary sampling parameters:
    (a) relative $L_2$ error and (b) $\mathrm{cond}(\bA_r)$
    against the number of sample points per Dirichlet edge $N_d$;
    (c) relative $L_2$ error and (d) $\mathrm{cond}(\bA_r)$
    against the Dirichlet radius multiplier $\rho_d$;
    (e) relative $L_2$ error and (f) $\mathrm{cond}(\bA_l)$
    against the number of sample points per Neumann edge $N_n$; and
    (g) relative $L_2$ error and (h) $\mathrm{cond}(\bA_l)$
    against the Neumann radius multiplier $\rho_n$.
    Dotted lines mark the reference configuration.}}
    \label{fig:2dpoisson-sens-boundary}
\end{figure}
\refB{
We first analyze the local sensitivity around the reference case.
Figure~\ref{fig:2dpoisson-sens-boundary} shows how the boundary sampling parameters
$N_d$, $\rho_d$, $N_n$ and $\rho_n$ affect the resulting solution error,
along with the condition number of the associated buffer coefficient system.
In general, the condition numbers increase with both the number of sample points
and the radius multipliers.
This is the expected behavior of a radial basis system:
adding sample points brings the basis centers closer together,
and raising the multiplier widens each basis function,
so in either case neighboring basis functions overlap more
and become increasingly linearly dependent,
driving the coefficient system toward singularity~\cite{schaback1995error}.
The resulting solution error, however, seems not so correlated with
the condition number.
While the condition numbers vary almost monotonically,
each error curve attains an interior minimum,
and the parameter value at that minimum is not the one that minimizes the conditioning.
}
\begin{figure}[tbh]
    \begin{tikzpicture}[
font=\small,
]
    \begin{groupplot}[
        group style={
            group name = my plots,
            group size= 3 by 2,
            horizontal sep=1.9cm,
            vertical sep=2.2cm,
        },
        height = 0.28\textwidth,
        width = 0.30\textwidth,
        ymode=log,
        tick scale binop ={\times},
        name=chung,
    ]
\pgfplotsset{set layers=standard}%

        \nextgroupplot[
            xlabel={$N_i$},
            ylabel={Relative $L_2$ error},
            xmin=3, xmax=11, xtick={4,6,8,10},
            ymin=1e-3, ymax=1e0,
        ]

            \addplot [gray, dotted, mark=none, line width=0.8]
            coordinates {(8, 1e-8) (8, 1e8)};

            \addplot [blue, line width=1.0, mark=*, mark size=1.8pt,
                      mark options={fill=white,}]
            table [x index=0, y index=1,
            ]{data_2dpoisson_sensitivity_slice_n_itf.txt};

        \nextgroupplot[
            xlabel={$N_i$},
            ylabel={$\mathrm{cond}(\bA_l)$},
            xmin=3, xmax=11, xtick={4,6,8,10},
            ymin=1e1, ymax=1e3,
        ]

            \addplot [gray, dotted, mark=none, line width=0.8]
            coordinates {(8, 1e-8) (8, 1e8)};

            \addplot [blue, line width=1.0, mark=*, mark size=1.8pt,
                      mark options={fill=white,}]
            table [x index=0, y index=2,
            ]{data_2dpoisson_sensitivity_slice_n_itf.txt};

        \nextgroupplot[
            xlabel={$N_i$},
            ylabel={$\mathrm{cond}(\bA_r)$},
            xmin=3, xmax=11, xtick={4,6,8,10},
            ymin=1e1, ymax=1e3,
        ]

            \addplot [gray, dotted, mark=none, line width=0.8]
            coordinates {(8, 1e-8) (8, 1e8)};

            \addplot [blue, line width=1.0, mark=*, mark size=1.8pt,
                      mark options={fill=white,}]
            table [x index=0, y index=3,
            ]{data_2dpoisson_sensitivity_slice_n_itf.txt};

        \nextgroupplot[
            xlabel={$\rho_i$},
            ylabel={Relative $L_2$ error},
            xmin=0.5, xmax=1.5, xtick={0.6,0.8,1.0,1.2,1.4},
            ymin=1e-3, ymax=1e0,
        ]

            \addplot [gray, dotted, mark=none, line width=0.8]
            coordinates {(0.6, 1e-8) (0.6, 1e8)};

            \addplot [blue, line width=1.0, mark=*, mark size=1.8pt,
                      mark options={fill=white,}]
            table [x index=0, y index=1,
            ]{data_2dpoisson_sensitivity_slice_rb_itf_factor.txt};

        \nextgroupplot[
            xlabel={$\rho_i$},
            ylabel={$\mathrm{cond}(\bA_l)$},
            xmin=0.5, xmax=1.5, xtick={0.6,0.8,1.0,1.2,1.4},
            ymin=1e1, ymax=1e3,
        ]

            \addplot [gray, dotted, mark=none, line width=0.8]
            coordinates {(0.6, 1e-8) (0.6, 1e8)};

            \addplot [blue, line width=1.0, mark=*, mark size=1.8pt,
                      mark options={fill=white,}]
            table [x index=0, y index=2,
            ]{data_2dpoisson_sensitivity_slice_rb_itf_factor.txt};

        \nextgroupplot[
            xlabel={$\rho_i$},
            ylabel={$\mathrm{cond}(\bA_r)$},
            xmin=0.5, xmax=1.5, xtick={0.6,0.8,1.0,1.2,1.4},
            ymin=1e1, ymax=1e3,
        ]

            \addplot [gray, dotted, mark=none, line width=0.8]
            coordinates {(0.6, 1e-8) (0.6, 1e8)};

            \addplot [blue, line width=1.0, mark=*, mark size=1.8pt,
                      mark options={fill=white,}]
            table [x index=0, y index=3,
            ]{data_2dpoisson_sensitivity_slice_rb_itf_factor.txt};

  \end{groupplot}
\node[below = 1.2cm of my plots c1r1.south west, anchor=west] {(a)};
\node[below = 1.2cm of my plots c2r1.south west, anchor=west] {(b)};
\node[below = 1.2cm of my plots c3r1.south west, anchor=west] {(c)};
\node[below = 1.2cm of my plots c1r2.south west, anchor=west] {(d)};
\node[below = 1.2cm of my plots c2r2.south west, anchor=west] {(e)};
\node[below = 1.2cm of my plots c3r2.south west, anchor=west] {(f)};
\end{tikzpicture}
%
    \caption{\refB{A local sensitivity study for the interface sampling parameters:
    (a) relative $L_2$ error, (b) $\mathrm{cond}(\bA_l)$ and (c) $\mathrm{cond}(\bA_r)$
    against the number of interface sample points $N_i$; and
    (d) relative $L_2$ error, (e) $\mathrm{cond}(\bA_l)$ and (f) $\mathrm{cond}(\bA_r)$
    against the radius multiplier $\rho_i$.
    Dotted lines mark the reference configuration.}}
    \label{fig:2dpoisson-sens-itf}
\end{figure}
\refB{
Similar trends are observed for the interface sampling parameters
$N_i$ and $\rho_i$, as shown in Figure~\ref{fig:2dpoisson-sens-itf}.
Since the interface is shared by the two subdomains,
these parameters enter both coefficient systems.
Widening the interface basis functions raises the condition number of both,
whereas adding interface sample points raises that of $\bA_l$
while leaving that of $\bA_r$ nearly unchanged.
The error, in turn, is strongly correlated with $\rho_i$, growing monotonically by
almost two orders of magnitude across the range examined,
but is again non-monotone in $N_i$.
}
\par
\begin{figure}[tbph]
    \begin{tikzpicture}[
]
    \begin{groupplot}[
        group style={
            group name = my plots,
            group size= 2 by 1,
            horizontal sep=2.2cm,
        },
        height = 0.32\textwidth,
        width = 0.42\textwidth,
        xmode=log,
        ymode=log,
        xmin=1e1, xmax=1e6,
        xlabel={$\mathrm{cond}(\bA)$},
        tick scale binop ={\times},
        axis on top=true,
        name=chung,
    ]
\pgfplotsset{set layers=standard}%

        \nextgroupplot[
            ylabel={$\cJ$},
            ymin=1e-5, ymax=1e2,
        ]

            \edef\imagepath{%
                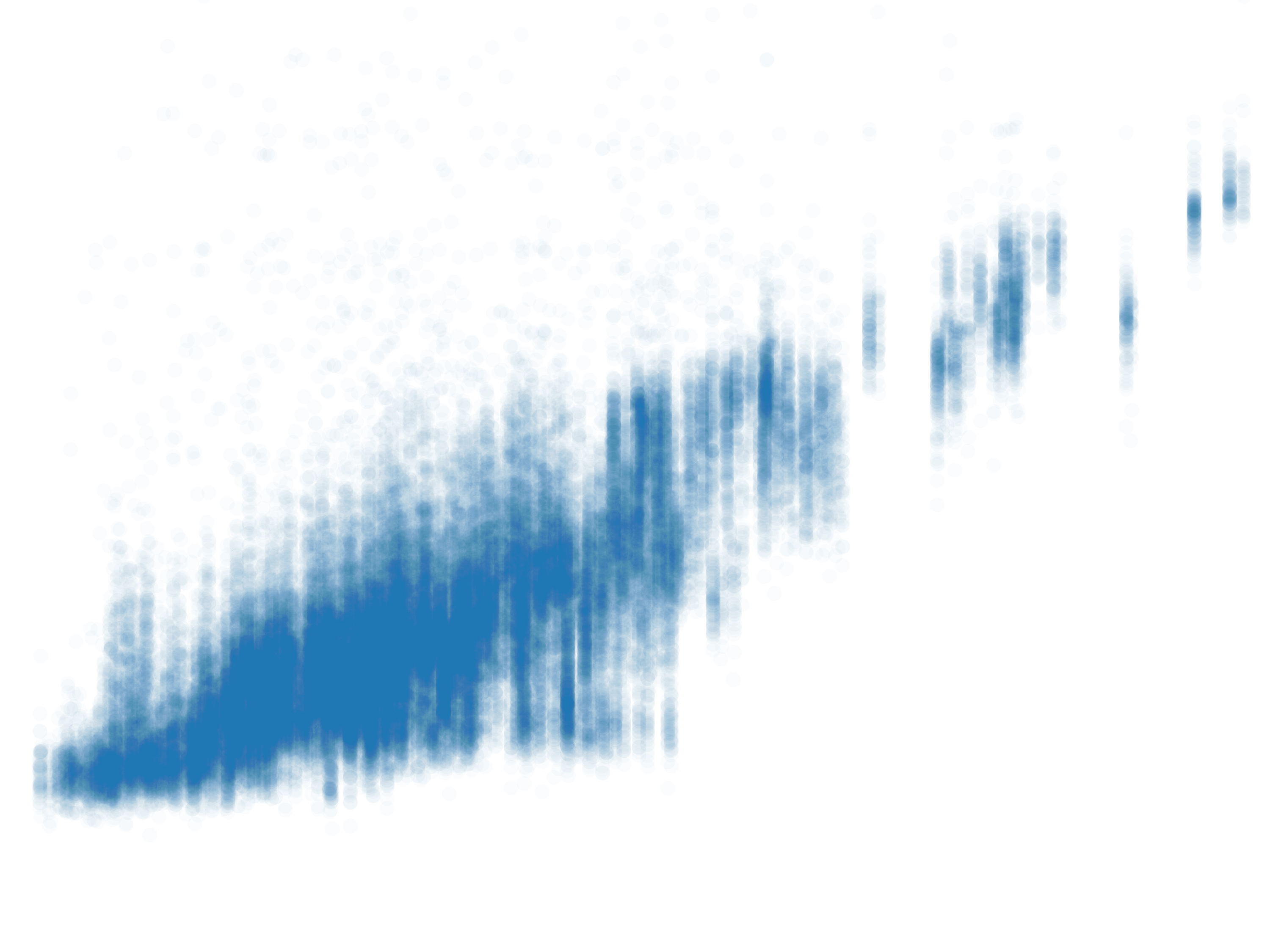}
            \addplot graphics[
                xmin=1e1, xmax=1e6,
                ymin=1e-5, ymax=1e2,
            ]{\imagepath};

        \nextgroupplot[
            ylabel={Relative $L_2$ error},
            ymin=1e-3, ymax=1e1,
        ]

            \edef\imagepath{%
                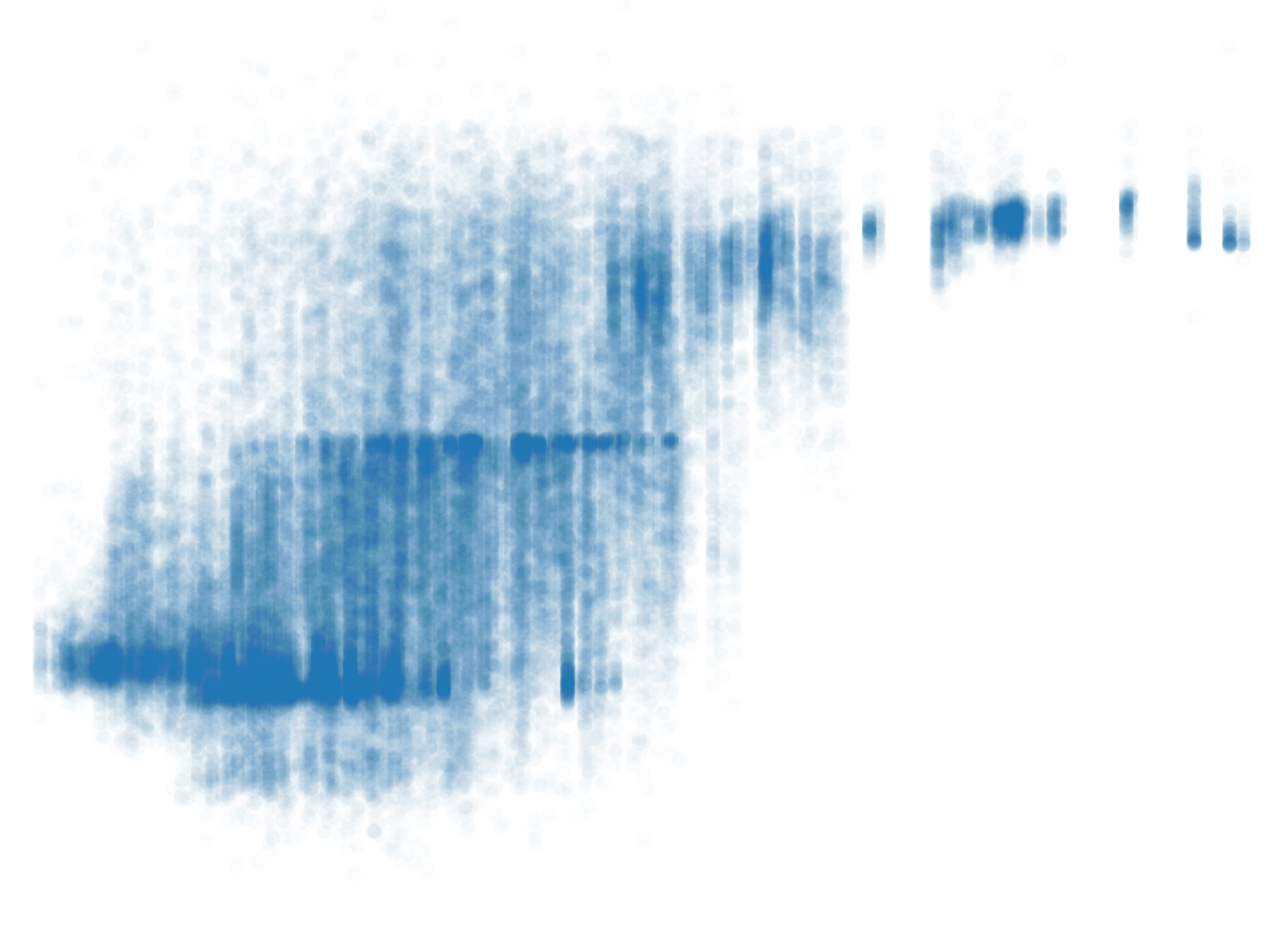}
            \addplot graphics[
                xmin=1e1, xmax=1e6,
                ymin=1e-3, ymax=1e1,
            ]{\imagepath};

  \end{groupplot}
\node[below = 1.2cm of my plots c1r1.south west, anchor=west] {(a)};
\node[below = 1.2cm of my plots c2r1.south west, anchor=west] {(b)};
\end{tikzpicture}
%
    \caption{\refB{
    A global sensitivity study of (a) lowest training loss attained and (b) relative $L_2$ error
    with respect to the conditioning of the combined buffer
    system $\bA = \mathrm{diag}\left(\bA_l, \bA_r\right)$.
    One marker per parameter set.}}
    \label{fig:2dpoisson-sens-cond}
\end{figure}
\refB{
The overall picture changes once we consider the global trend over all parameter ranges.
Figure~\ref{fig:2dpoisson-sens-cond}~(a) shows a scatter plot of the training loss
against the conditioning of the combined buffer coefficient system
$\bA = \mathrm{diag}\left(\bA_l, \bA_r\right)$,
which is the system the two block solves present to the optimizer as a whole.
Unlike the previous local sensitivity studies,
a strong correlation between the training loss and $\mathrm{cond}(\bA)$ is observed.
This does not mean, however, that the parameters minimizing $\mathrm{cond}(\bA)$
are the ones to use.
The best-conditioned systems are those with the fewest sample points and the narrowest
basis functions, which enforce the boundary and interface conditions at correspondingly
few locations with a buffer function that decays quickly between them,
so the conditions are met only loosely away from the sample points.
Figure~\ref{fig:2dpoisson-sens-cond}~(b) shows the corresponding solution error
against the same conditioning.
While the relative $L_2$ error follows the training loss closely,
the lowest condition numbers do not yield the lowest solution error.
The lower envelope of the scatter instead attains its minimum around
$\mathrm{cond}(\bA)\sim 2\times10^2$,
where the best parameter set reaches $4.8\times10^{-3}$,
against $1.2\times10^{-2}$ at the best-conditioned end of the sweep.
This presumably reflects a balance between the ability of the buffer functions to
represent the boundary and interface conditions and the conditioning of the system that
determines their coefficients.
While this condition number may serve as an empirical starting point for choosing the
boundary and interface hyperparameters,
we note that such an ``optimal'' value is likely geometry- and problem-specific,
and establishing its generality is beyond the scope of this study.
}
\par
\begin{figure}[tbph]
    \begin{tikzpicture}[
font=\small,
]
    \begin{groupplot}[
        group style={
            group name = my plots,
            group size= 4 by 4,
            x descriptions at=edge bottom,
            y descriptions at=edge left,
            horizontal sep=0.5cm,
            vertical sep=0.6cm,
        },
        height = 0.20\textwidth,
        width = 0.28\textwidth,
        xmode=log,
        xmin=1e-3, xmax=1e1,
        xtick={1e-2, 1e0},
        ymin=0, ymax=0.45,
        ytick={0, 0.2, 0.4},
        tick scale binop ={\times},
        tick label style={font=\scriptsize},
        title style={font=\small, yshift=-0.8ex},
        every axis plot/.append style={
            const plot,
            no markers,
            draw=blue!60!black,
            line width=0.6,
        },
        name=chung,
    ]
\pgfplotsset{set layers=standard}%

        \nextgroupplot[title={$N_i=4$}]
            \addplot table [x index=0, y index=1,
            ]{data_2dpoisson_sensitivity_hist_cells_rel_l2.txt};
        \nextgroupplot[title={$N_i=6$}]
            \addplot table [x index=0, y index=2,
            ]{data_2dpoisson_sensitivity_hist_cells_rel_l2.txt};
        \nextgroupplot[title={$N_i=8$}]
            \addplot table [x index=0, y index=3,
            ]{data_2dpoisson_sensitivity_hist_cells_rel_l2.txt};
        \nextgroupplot[title={$N_i=10$}]
            \addplot table [x index=0, y index=4,
            ]{data_2dpoisson_sensitivity_hist_cells_rel_l2.txt};

        \nextgroupplot
            \addplot table [x index=0, y index=5,
            ]{data_2dpoisson_sensitivity_hist_cells_rel_l2.txt};
        \nextgroupplot
            \addplot table [x index=0, y index=6,
            ]{data_2dpoisson_sensitivity_hist_cells_rel_l2.txt};
        \nextgroupplot
            \addplot table [x index=0, y index=7,
            ]{data_2dpoisson_sensitivity_hist_cells_rel_l2.txt};
        \nextgroupplot
            \addplot table [x index=0, y index=8,
            ]{data_2dpoisson_sensitivity_hist_cells_rel_l2.txt};

        \nextgroupplot
            \addplot table [x index=0, y index=9,
            ]{data_2dpoisson_sensitivity_hist_cells_rel_l2.txt};
        \nextgroupplot
            \addplot table [x index=0, y index=10,
            ]{data_2dpoisson_sensitivity_hist_cells_rel_l2.txt};
        \nextgroupplot
            \addplot table [x index=0, y index=11,
            ]{data_2dpoisson_sensitivity_hist_cells_rel_l2.txt};
        \nextgroupplot
            \addplot table [x index=0, y index=12,
            ]{data_2dpoisson_sensitivity_hist_cells_rel_l2.txt};

        \nextgroupplot
            \addplot table [x index=0, y index=13,
            ]{data_2dpoisson_sensitivity_hist_cells_rel_l2.txt};
        \nextgroupplot
            \addplot table [x index=0, y index=14,
            ]{data_2dpoisson_sensitivity_hist_cells_rel_l2.txt};
        \nextgroupplot
            \addplot table [x index=0, y index=15,
            ]{data_2dpoisson_sensitivity_hist_cells_rel_l2.txt};
        \nextgroupplot
            \addplot table [x index=0, y index=16,
            ]{data_2dpoisson_sensitivity_hist_cells_rel_l2.txt};

  \end{groupplot}
\node[below = 1.0cm of my plots c2r4.south east, anchor=north]
    {Relative $L_2$ error};
\node[rotate=90, anchor=center] at ([xshift=-1.15cm] my plots c1r2.south west)
    {Fraction of runs};
\node[rotate=-90, anchor=center] at ([xshift=0.45cm] my plots c4r1.east) {$N_d=2$};
\node[rotate=-90, anchor=center] at ([xshift=0.45cm] my plots c4r2.east) {$N_d=4$};
\node[rotate=-90, anchor=center] at ([xshift=0.45cm] my plots c4r3.east) {$N_d=6$};
\node[rotate=-90, anchor=center] at ([xshift=0.45cm] my plots c4r4.east) {$N_d=8$};
\end{tikzpicture}
%
    \caption{\refB{Distribution of the relative $L_2$ error over the global sweep,
    grouped by the number of Dirichlet sample points $N_d$ and interface sample points
    $N_i$. Each panel collects $4{,}500$ configurations,
    and all panels share the same axes.}}
    \label{fig:2dpoisson-sens-cells}
\end{figure}
\refB{
Interestingly, Figure~\ref{fig:2dpoisson-sens-cond}~(b) shows three distinct groups of
solution errors, located at relative $L_2$ errors of around
$1.8\times10^{-2}$, $1.4\times10^{-1}$ and $1.1$, roughly a decade apart.
It is not entirely clear how these basins are formed,
and they are not associated with any single hyperparameter.
We found them to be mostly associated with the interplay between the number of Dirichlet
sample points $N_d$ and the number of interface sample points $N_i$.
Figure~\ref{fig:2dpoisson-sens-cells} shows how the error distribution changes
as $N_d$ and $N_i$ increase.
At the lowest $N_d=2$ and $N_i=4$, the error peaks at $1.8\times10^{-2}$.
As $N_i$ increases to $N_i=10$, the peak ``transitions'' to $1.4\times10^{-1}$.
Then, with $N_i=10$ fixed and $N_d$ increasing,
the peak transitions once more, to $1.1$.
A similar trend is observed along the path that increases $N_d$ first and then $N_i$,
and in the remaining rows and columns.
We conjecture that both parameters act on $\bA_r$, which determines the value of the
buffer function rather than its slope, so that a degraded solve there displaces the
solution level itself.
The Neumann rows, in contrast, prescribe only normal derivatives,
and the solution error accordingly varies far less with $N_n$:
its median stays between $5.4\times10^{-2}$ and $8.4\times10^{-2}$ over the whole range,
against the tenfold variation seen for both $N_d$ and $N_i$.
}
\par
\begin{figure}[tbph]
    \begin{tikzpicture}[
]
    \begin{groupplot}[
        group style={
            group name = my plots,
            group size= 2 by 2,
            y descriptions at=edge left,
            horizontal sep=0.8cm,
            vertical sep=2.2cm,
        },
        height = 0.30\textwidth,
        width = 0.42\textwidth,
        xmode=log,
        ymin=0,
        ylabel={Fraction of runs},
        tick scale binop ={\times},
        scaled y ticks=false,
        yticklabel style={/pgf/number format/fixed},
        every axis plot/.append style={
            const plot,
            no markers,
            line width=0.9,
        },
        legend style={
            draw=none, fill=none,
            at={(rel axis cs: 0.5, 1.02)},
            anchor=south west,
            legend cell align={left},
            legend columns=3,
            nodes={scale=0.9},
            /tikz/every even column/.append style={column sep=0.4cm},
        },
        name=chung,
    ]
\pgfplotsset{set layers=standard}%

        \nextgroupplot[
            xlabel={$\cJ$},
            xmin=1e-5, xmax=1e2,
            ymax=0.09,
        ]

            \addplot [blue, dashed]
            table [x index=0, y index=1,
            ]{data_2dpoisson_sensitivity_hist_ratio_best_loss.txt};
            \addplot [black]
            table [x index=0, y index=2,
            ]{data_2dpoisson_sensitivity_hist_ratio_best_loss.txt};
            \addplot [red, dash dot]
            table [x index=0, y index=3,
            ]{data_2dpoisson_sensitivity_hist_ratio_best_loss.txt};

            \legend{$\gamma=0.1$, $\gamma=1$, $\gamma=10$};

        \nextgroupplot[
            xlabel={$\cJ$},
            xmin=1e-5, xmax=1e2,
            ymax=0.09,
        ]

            \addplot [blue, dashed]
            table [x index=0, y index=4,
            ]{data_2dpoisson_sensitivity_hist_ratio_best_loss.txt};
            \addplot [black]
            table [x index=0, y index=5,
            ]{data_2dpoisson_sensitivity_hist_ratio_best_loss.txt};
            \addplot [red, dash dot]
            table [x index=0, y index=6,
            ]{data_2dpoisson_sensitivity_hist_ratio_best_loss.txt};

        \nextgroupplot[
            xlabel={Relative $L_2$ error},
            xmin=1e-3, xmax=1e1,
            ymax=0.12,
        ]

            \addplot [blue, dashed]
            table [x index=0, y index=1,
            ]{data_2dpoisson_sensitivity_hist_ratio_rel_l2.txt};
            \addplot [black]
            table [x index=0, y index=2,
            ]{data_2dpoisson_sensitivity_hist_ratio_rel_l2.txt};
            \addplot [red, dash dot]
            table [x index=0, y index=3,
            ]{data_2dpoisson_sensitivity_hist_ratio_rel_l2.txt};

        \nextgroupplot[
            xlabel={Relative $L_2$ error},
            xmin=1e-3, xmax=1e1,
            ymax=0.12,
        ]

            \addplot [blue, dashed]
            table [x index=0, y index=4,
            ]{data_2dpoisson_sensitivity_hist_ratio_rel_l2.txt};
            \addplot [black]
            table [x index=0, y index=5,
            ]{data_2dpoisson_sensitivity_hist_ratio_rel_l2.txt};
            \addplot [red, dash dot]
            table [x index=0, y index=6,
            ]{data_2dpoisson_sensitivity_hist_ratio_rel_l2.txt};

  \end{groupplot}
\node[below = 1.5cm of my plots c1r1.south west, anchor=west] {(a) $\gamma_0$};
\node[below = 1.5cm of my plots c2r1.south west, anchor=west] {(b) $\gamma_1$};
\node[below = 1.5cm of my plots c1r2.south west, anchor=west] {(c) $\gamma_0$};
\node[below = 1.5cm of my plots c2r2.south west, anchor=west] {(d) $\gamma_1$};
\end{tikzpicture}
%
    \caption{\refA{Distribution of the training loss (a, b) and of the relative $L_2$
    error (c, d) over the global sweep, grouped by the interface split ratios
    $\gamma_0$ and $\gamma_1$.}}
    \label{fig:2dpoisson-sens-ratios}
\end{figure}
\refA{
Figure~\ref{fig:2dpoisson-sens-ratios} shows the training loss and the solution error
distribution over the global parameter sweep,
with respect to the interface split ratios $\gamma_0$ and $\gamma_1$.
Little sensitivity is observed to either,
with the distributions retaining essentially the same shape in both cases.
This is somewhat expected in that $\gamma_0$ and $\gamma_1$ change only the right-hand
side of the buffer system, leaving $\bA$ and hence its conditioning untouched.
This shows that the buffer approach is not strongly tied to a particular division of
the interface conditions between the two subdomains.
The even split $\gamma_0=\gamma_1=1$ used throughout this study is a sound default,
while the freedom to choose otherwise remains available for problems in which the two
subdomains are less comparable.
}


\bibliographystyle{elsarticle-num} 
\bibliography{references}

%
%
%
%
\end{document}